\documentclass[reqno,a4paper]{amsart}
\usepackage{amssymb,amsthm}
\usepackage{ifthen}
\usepackage{datetime}

\usepackage[linktocpage=true]{hyperref}
\usepackage{graphicx,color}

\usepackage[normalem]{ulem}

\title[Hidden symmetries and decay for the wave equation on Kerr]{Hidden symmetries and decay for the wave equation on the Kerr spacetime}
\author[L. Andersson]{Lars Andersson${}^\dagger$} \email{laan@aei.mpg.de}
\address{Albert Einstein Institute, Am M\"uhlenberg 1, D-14476 Potsdam,
  Germany 
}

\author[P. Blue]{Pieter Blue${}^\ddagger$} \email{P.Blue@ed.ac.uk}
\address{The School of Mathematics and the Maxwell Institute,  University of Edinburgh,James Clerk Maxwell Building, 
The King's Buildings, 
Mayfield Road, 
Edinburgh, 
Scotland EH9 3JZ,UK}

\date{\today \ {\em File:\jobname{.tex} at \currenttime}}

\newcommand{\blue}[1]{{\color{blue} #1}}

\newcounter{mnotecount}[section]

\newcounter{mymnotecount}[section]
\renewcommand{\themymnotecount}{\thesection.\arabic{mymnotecount}}
\newcommand{\mymnote}[1]{\protect{\stepcounter{mymnotecount}}${\raisebox{0.5\baselineskip}[0pt]{\makebox[0pt][c]{\color{magenta}{\tiny\em$\bullet$\themymnotecount}}}}$\marginpar{\raggedright\tiny\em$\!\!\!\!\!\!\,\bullet$\themymnotecount: \blue{#1}}\ignorespaces}
\renewcommand{\mymnote}[1]{}

\newcommand{\newpageForSubsection}{}

\newtheorem{theorem}{Theorem}
\newtheorem{lemma}[theorem]{Lemma}
\newtheorem{definition}[theorem]{Definition}
\newtheorem{remark}[theorem]{Remark}
\newtheorem{corollary}[theorem]{Corollary}
\numberwithin{equation}{section} 
\numberwithin{theorem}{section}
\newcounter{step}
\newenvironment{multistep}{\setcounter{step}{0}}{}
\newcommand{\newstep}[1]{\addtocounter{step}{1}{\bf Step \thestep: #1} }
\newcommand{\Naturals}{\mathbb{N}}

\newcommand{\Reals}{\mathbb R}

\newcommand{\Id}{\text{Id}}
\newcommand{\Indicator}{\mathbf{1}}
\newcommand{\supp}{\text{supp}}

\newcommand{\KDelta}{\Delta}
\newcommand{\KSigma}{\Sigma}
\newcommand{\KPi}{\Pi}
\newcommand{\rp}{r_+}

\newcommand{\CurlyR}{\mathcal{R}}

\newcommand{\Dt}{\frac{\di}{\di t}}
\newcommand{\dr}{\partial_r}

\newcommand{\dt}{\partial_t}
\newcommand{\dtheta}{\partial_\theta}
\newcommand{\dphi}{\partial_\phi}
\newcommand{\dx}{\partial_x}
\newcommand{\dAng}{{\nabla\!\!\!\!/}}
\newcommand{\Li}{\Theta_i}
\newcommand{\LsubThree}{\Theta_{3}}

\newcommand{\di}{\mathrm{d}} %{\text{d}}
\newcommand{\diNormal}{\di \eta}
\newcommand{\diFourNatural}{\gVol\di^4 x}
\newcommand{\diTwo}{\di^2\mu}
\newcommand{\diThree}{\di^3\mu}
\newcommand{\diFour}{\di^4\mu}

\newcommand{\Svol}{\mu}
\newcommand{\gVol}{\sqrt{|\gMetric|}}

\newcommand{\vecfont}[1]{\mathbf{#1}}
\newcommand{\tensorfont}[1]{#1}
\newcommand{\STwovecfont}[1]{\underline{\overline{\vecfont{#1}}}}
\newcommand{\STwofnfont}[1]{\underline{\overline{#1}}}

\newcommand{\solu}{\psi}
\newcommand{\gMetric}{\tensorfont{g}}
\newcommand{\vecX}{\vecfont{X}}
\newcommand{\vecY}{\vecfont{Y}}
\newcommand{\formxi}{\xi}
\newcommand{\MMTTscalar}{q}

\newcommand{\gWave}{\square}

\newcommand{\GenEnergyGeodesic}[1]{e_{#1}}
\newcommand{\rorbit}{r_{o}} %aka \rRandom
\newcommand{\vecMGeodesic}{A}
\newcommand{\fnMrGeodesic}{\mathcal{F}}
\newcommand{\fnMpGeodesic}{q_{\text{reduced}}}
\newcommand{\GeodesicL}{\mathbf{L}}

\newcommand{\gMetricGeneral}{\tensorfont{g}}
\newcommand{\StressEnergyGeneral}{\tensorfont{T}}
\newcommand{\GenMomentum}[1]{P_{#1}}
\newcommand{\GenEnergy}[1]{E_{#1}}
\newcommand{\GenEnergyOrder}[2]{\GenEnergy{#1,#2}}

\newcommand{\SymGeneraln}[1]{\mathbb{S}_{#1}}
\newcommand{\CQA}{S}
\newcommand{\ua}{{\underline{a}}}
\newcommand{\ub}{{\underline{b}}}
\newcommand{\uc}{{\underline{c}}}
\newcommand{\Lie}{\mathcal{L}}

\newcommand{\HypersurfaceGeneral}{\Sigma}
\newcommand{\STwovecX}{\STwovecfont{X}}
\newcommand{\STwoMMTTscalar}{\STwofnfont{q}}

\newcommand{\Geodesic}{\gamma}
\newcommand{\GeodesicVelocity}{\dot{\gamma}}

\newcommand{\GeodesicEnergy}{\boldsymbol{e}}
\newcommand{\GeodesicLz}{\boldsymbol{\ell_z}}
\newcommand{\GeodesicQ}{\boldsymbol{q}}
\newcommand{\TensorQ}{Q}
\newcommand{\KCarter}{\boldsymbol{k}}
\newcommand{\TensorKCarter}{K}
\newcommand{\nq}{{n_\OpQ}}
\newcommand{\nt}{{n_t}}
\newcommand{\np}{{n_\phi}}

\newcommand{\TensorK}{K}

\newcommand{\VFTrueSym}{\mathbb{S}_1}
\newcommand{\SOSym}{\SymGeneraln{2}}
\newcommand{\OpQ}{Q}

\newcommand{\pairXq}{(\vecX,\MMTTscalar)}

\newcommand{\hst}[1]{\Sigma_{#1}}

\newcommand{\vecNormalSigmat}{n_{\hst{t}}}
\newcommand{\omegaperp}{\omega_\perp}

\newcommand{\vecTperp}{\vecfont{T}_{\perp}}
\newcommand{\omegaH}{\omega_H}

\newcommand{\OpGeneral}{X}
\newcommand{\OpSetGeneral}{\mathbb{X}}
\newcommand{\normPtwiseTn}[2]{|#2|_{#1}}
\newcommand{\SLap}{{\Delta \!\!\!\! /}}
\newcommand{\OpL}{\mathcal{L}}
\newcommand{\OpLEpsilon}[1]{\OpL_{#1}}
\newcommand{\normPtwiseTnHomo}[2]{\normPtwiseTn{#1,1}{#2}}
\newcommand{\normPtwiseTnEpsilon}[3]{|#3|_{#1,#2}}

\newcommand{\VFAngAlg}{\mathbb{O}_1}
\newcommand{\VFEffSym}{\mathbb{T}_1}

\newcommand{\vecTBlend}{\vecfont{T}_\chi}
\newcommand{\fnBlend}{\chi}
\newcommand{\fnAtBlendLocation}{\Indicator_{\supp\fnBlend'}}

\newcommand{\CurlyRTilde}{\tilde{\mathcal{R}}}
\newcommand{\DiffCurlyRTilde}{{\tilde{\mathcal{R}}'}{}}
\newcommand{\DiffCurlyRTTilde}{\tilde{\tilde{\mathcal{R}}}'}
\newcommand{\DDiffCurlyRTTilde}{\tilde{\tilde{\mathcal{R}}}''}
\newcommand{\vecM}{{\vecfont{A}^{\ua\ub}}}            % Morawetz vector field
\newcommand{\vecMr}{\mathcal{F}^{\ua\ub}}   
\newcommand{\fnM}{{q^{\ua\ub}}}             % MMTT correction
\newcommand{\fnMp}{{q_{\text{reduced}}^{\ua\ub}}}
\newcommand{\pairM}{\STwovecfont{A}}
\newcommand{\vecMone}{{\vecfont{A}^\ua}}
\newcommand{\vecMrone}{\mathcal{F}^{\ua}}
\newcommand{\fnMone}{q^\ua}
\newcommand{\fnMpone}{q_{\text{reduced}}^\ua}
\newcommand{\CurlyA}{\mathcal{A}}
\newcommand{\CurlyAone}{\mathcal{A}}
\newcommand{\CurlyU}{\mathcal{U}}
\newcommand{\CurlyUfour}{\mathcal{U}}
\newcommand{\CurlyV}{\mathcal{V}}
\newcommand{\CurlyVone}{\mathcal{V}}

\newcommand{\MorBoundaryTermA}[1]{B_{#1;\text{I}}}

\newcommand{\newcommandMG}[3]{\newcommand{#1}{#3}} % A notation
\newcommandMG{\fnMna} {{f_{1}}} {z}
\newcommandMG{\fnMnb} {{f_{2}}} {w}
\newcommandMG{\fnMca} {f_{1,1}}   {z_1}
\newcommandMG{\fnMcb} {f_{2,1}}   {w_1}
\newcommandMG{\fnMda} {f_{1,2}}   {z_2}
\newcommandMG{\fnMdb} {f_{2,2}}   {w_2}
\newcommand{\MorBoundaryTermB}[1]{B_{#1;\text{II}}}
\newcommand{\normPtwiseTnEpsilondtTurnedOn}[2]{|#2|_{#1,\epsilonMorawetzTurnOndtSquared}}

\newcommand{\fnHardyToEstimate}{\phi}
\newcommand{\fnHardyRedToEstimate}{\varphi}
\newcommand{\PotlGeneraldrrCoeff}{A}
\newcommand{\PotlGeneralLOCoeff}{V}
\newcommand{\fnHardyODESolu}{u}
\newcommand{\fnHardyODEAux}{f}
\newcommand{\fnHardyRedODESolu}{v}
\newcommand{\PotlRedGeneralLOCoeff}{W}
\newcommand{\fnHardyRRedODESolu}{\tilde{v}}
\newcommand{\fnHardyRRRedODESolu}{\Tilde\psi}
\newcommand{\PotlRedGeneralLOCoeffPert}{\bar{W}}
\newcommand{\aAlt}{\tilde{a}}

\newcommand{\localiseAwayFromPhotonOrbits}{\Indicator_{r\not\eqsim 3M}}
\newcommand{\EnergyThree}{\GenEnergyOrder{\vecTBlend}{3}[\solu]}
\newcommand{\slabST}[2]{[#1,#2]\times(\rp,\infty)\times S^2}
\newcommand{\Lz}{L_z}
\newcommand{\vecMclassical}{\vecfont{A}}
\newcommand{\fnMclassical}{q}
\newcommand{\fnMpclassical}{q_{\text{reduced}}}
\newcommand{\fnMncClassical}{f}
\newcommand{\pairClassical}{(\vecMclassical,\fnMclassical)}
\newcommand{\CurlyALzZero}{\mathcal{A}}
\newcommand{\CurlyULzZero}{\mathcal{U}}
\newcommand{\CurlyVLzZero}{\mathcal{V}}

\newcommand{\ModelEnergyThree}{E_{\text{model},3}}

\renewcommand{\vecfont}[1]{#1}
\renewcommand{\tensorfont}[1]{#1}
\renewcommand{\STwovecfont}[1]{\mathbf{#1}}
\renewcommand{\STwofnfont}[1]{\mathbf{#1}}

\newcommand{\newepsilon}[3]{\newcommand{#1}{#2}}
 
\newepsilon{\ConsantIntroUniformBound}{C}{C_1}
\newepsilon{\epsilonSlowRotationIntroUniformBound}{\bar{a}}{\bar{a}_1}
\newepsilon{\ConstantIntroDecay}{C}{C_2}
\newepsilon{\epsilonSlowRotationIntroDecay}{\bar{a}}{\bar{a}_2}
\newepsilon{\ConsantIntroUniformBoundsolu}{C}{C_3}
\newepsilon{\epsilonSlowRotationIntroUniformBoundSolu}{\bar{a}}{\bar{a}_3}
\newepsilon{\ConstantSTwoSobolev}{C}{C_{S}}
\newepsilon{\ConstantLEpsilonL}{C}{C_{\OpL}}
\newepsilon{\epsilonInvBlendLocation}{r_{\fnBlend}}{r_{\fnBlend}}
\newepsilon{\epsilonSlowRotationNearEnergyPositive}{\bar{a}}{\bar{a}_{T}}
\newepsilon{\epsilonMorawetzSlowRotation}{\bar{a}}{\bar{a}_{\text{Morawetz}}}
\newepsilon{\epsilonMorawetzPhotonWidth}{\bar{r}}{\bar{r}_{\text{Morawetz}}}
\newepsilon{\epsilonMorawetzTurnOndtSquared}{\epsilon_{\dt^2}}{\epsilon_{\dt^2}}
\newepsilon{\epsilonMorawetzTurnOndtSquaredUpperBound}{\overline{\epsilonMorawetzTurnOndtSquared}}{\epsilon_{\dt^2,\text{upper bound}}}
\newepsilon{\epsilonInverseErrorInMorawetzBeforeHardy}{C}{C_\text{Morawetz error}}
\newepsilon{\epsilonInvIntegratedMorawetzHomo}{C_1}{C_{\text{Homogeneous Morawetz}}}
\newepsilon{\epsilonInvIntegratedMorawetzHomoB}{C_2}{C_{\text{Homogeneous Morawetz,2}}}
\newepsilon{\epsilonHardySlowRotation}{\bar{a}}{\bar{a}_\text{Hardy}}
\newcommand{\epsilonHardyExtraFactorOnPotential}{\epsilon_{\text{\normalfont{Hardy}}}}
\newcommand{\epsilonHardyExtraFactorOnPotentialB}{\epsilon_{\text{Hardy,1}}}
\newcommand{\epsilonHardyExtraFactorOnPotentialC}{\epsilon_{\text{Hardy,2}}}
\newcommand{\epsilonHardySlowRotationInPert}{\bar{a}_\text{Hardy,3}}
\newcommand{\epsilonHardyExtraFactorOnPotentialInPert}{\epsilon_{\text{Hardy,3}}}
\newepsilon{\epsilonInvIntegratedMorawetz}{C}{C_{\text{Morawetz}}}
\newepsilon{\epsilonBoundedEnergySlowRotation}{\bar{a}}{\bar{a}_E}
\newepsilon{\epsilonInvBoundedEnergy}{C}{C_E}
\newepsilon{\epsilonInvBoundedEnergyInner}{C'}{C_E'}
\newepsilon{\epsilonSlowRotationDominantEnergy}{\bar{a}}{\bar{a}_{DEC}}
\newepsilon{\ConstantDominantEnergy}{C}{C_{DEC}}
\newepsilon{\epsilonSlowRotationRelationBetweenEnergies}{\bar{a}}{\bar{a}_{\tilde{E}}}
\newepsilon{\ConstantRelationBetweenEnergies}{C}{C_{\tilde{E}}}
\newepsilon{\epsilonSlowRotationKTimelike}{\bar{a}}{\bar{a}_{TL}}
\newepsilon{\epsilonSlowRotationKEnergy}{\bar{a}}{\bar{a}_{KE}}

\newepsilon{\epsilonSlowRotationLightConeLocalise}{\bar{a}}{\bar{a}_{LC}}
\newepsilon{\ConstantLightConeLocalisation}{C}{C_{LC}}

\newepsilon{\epsilonSlowRotationInK}{\bar{a}}{\bar{a}_K}
\newepsilon{\ConstantInK}{C}{C_K}
\newepsilon{\ConstantInKexp}{C'}{C_K'}
\newepsilon{\epsilonSlowRotationOtherSfcs}{\bar{a}}{\bar{a}_{OS}}
\newepsilon{\ConstantOtherSfcs}{C}{C_{OS}}
\newepsilon{\ConstantOtherSfcsexp}{C'}{C_{OS}'}

\newepsilon{\epsilonSlowRotationStationaryDecay}{\bar{a}}{\bar{a}_{\text{SD}}}
\newepsilon{\ConstantStationaryDecayexp}{C'}{C_{SD}'}

\newepsilon{\epsilonSlowRotationND}{\bar{a}}{\bar{a}_{ND}}
\newepsilon{\ConstantND}{C}{C_{\text{ND}}}
\newepsilon{\ConstantNDexp}{C'}{C_{\text{ND}}'}
\newepsilon{\epsilonSlowRotationFarAlmostNullEnergy}{\bar{a}}{\bar{a}_{FANE}}
\newepsilon{\ConstantFarAlmostNullEnergy}{C}{C_{FANE}}
\newepsilon{\ConstantFarAlmostNullEnergyexp}{C'}{C_{FANE}'}
\newepsilon{\epsilonSlowRotationFarDecay}{\bar{a}}{\bar{a}_{FD}}
\newepsilon{\ConstantFarDecay}{C}{C_{FD}}
\newepsilon{\ConstantFarDecayExp}{C'}{C_{FD}'}

\newepsilon{\epsilonAzimuthal}{\epsilon_{\Lz=0}}

\newcommand{\rPSLzZero}{r_{\Lz=0}}

\newcommand{\tk}{\breve{t}}
\newcommand{\rk}{\breve{r}}
\newcommand{\hk}{\breve{\theta}}
\newcommand{\pk}{\breve{\phi}}
\newcommand{\dtk}{\partial_{\tk}}
\newcommand{\drk}{\partial_{\rk}}
\newcommand{\dhk}{\partial_{\hk}}
\newcommand{\dpk}{\partial_{\pk}}

\newcommand{\hstk}[1]{\breve{\Sigma}_{#1}}
\newcommand{\hsrk}[1]{\mathcal{H}_{#1}}
\newcommand{\Omegak}[1]{\Omega_{#1}}

\newepsilon{\epsilonInverseAppendix}{\bar{a}}

\newcommand{\OmegaH}{\Omega_H}

\newcommand{\opQT}{\breve{\mathcal{Q}}}

\newcommand{\vecTk}{\breve{\vecfont{T}}}
\newcommand{\epsilonNH}{\epsilon_{NH}}
\newcommand{\rNH}{r_{NH}}
\newcommand{\chiNH}{\chi_{NH}}

\newcommand{\f}{f}
\newcommand{\h}{h}

\newcommand{\normalhst}[1]{n_{\hst{#1}}}
\newcommand{\ONBasis}{\hat{X}}

\newcommand{\TimeAndAngularDerivativesTwoPtThreea}[1]{B_{\eqref{eq:ExplicitLepsilonLExpansionControlsTwoDerivatives}}[#1]}

\newcommand{\MorBoundaryTermBa}[1]{B_{#1;\text{IIa}}}

\begin{document}

\begin{abstract}
  Energy and decay estimates for the wave equation on the exterior
  region of slowly rotating Kerr spacetimes are proved. The method
  used is a generalisation of the vector-field method, which allows
  the use of higher-order symmetry operators. In particular, our
  method makes use of the second-order Carter operator, which is a
  hidden symmetry in the sense that it does not correspond to a
  Killing symmetry of the spacetime.
\end{abstract}

\maketitle

%\tableofcontents

\section{Introduction}
In this paper we prove boundedness and 
integrated energy decay for solutions of the covariant wave equation 
\begin{align*}
 \nabla^\alpha \nabla_\alpha \solu = 0
\end{align*}
in the exterior region of the slowly rotating Kerr spacetime. 
In Boyer-Lindquist coordinates $(x^\alpha) = (t,r,\theta,\phi)$, the exterior region is given by $\Reals\times(\rp,\infty)\times S^2$ with the Lorentz metric 
\begin{align}
\gMetric_{\mu\nu} \di x^\mu \di x^\nu \nonumber ={}&{}
-\left(1-\frac{2Mr}{\KSigma}\right)\di t^2 -\frac{4Mra\sin^2\theta}{\KSigma}
\di t\di\phi 
\nonumber\\ {}&{}+
\frac{\KPi\sin^2\theta}{\KSigma} \di\phi^2  +\frac{\KSigma}{\KDelta}\di r^2
 +\KSigma \di\theta^2 ,
\label{eq:KerrMetric}
\end{align}
where $\rp=M+\sqrt{M^2-a^2}$ and
\begin{align*}
\KDelta= r^2-2Mr +a^2 , \quad  \KSigma= r^2+a^2\cos^2\theta , \quad  \KPi=
(r^2+a^2)^2 -a^2\sin^2\theta\KDelta .
\end{align*}
For $0\leq |a|\leq M$, the Kerr family of metrics describe an asymptotically flat, stationary and axi-symmetric solution of the vacuum Einstein equations, containing a  rotating black
hole with mass $M$ and angular momentum $Ma$, and with horizon
located at $r=\rp$. The Schwarzschild spacetime is the subcase with
$a=0$. We will take $M>0$ fixed and will study the slowly rotating
case, $|a|\ll M$. The exterior region is globally hyperbolic, with the
surfaces of constant $t$, $\hst{t}$, as Cauchy surfaces. Thus, the
wave equation is well posed in the exterior region, even though the
Kerr spacetime can be extended. We consider initial data on the
hypersurface $\hst{0}$.

The isometry group of the Kerr spacetime is 2-dimensional, generated
by the stationary vector field $\dt$ which is timelike near spatial
infinity, and the axial rotation vector field $\dphi$. In a general 4-dimensional spacetime with a 2-dimensional isometry group, one may expect that there are only three constants of the motion for geodesics, and that the geodesic motion is chaotic. However, the fourth conserved quantity for geodesics in the Kerr spacetime discovered by Carter 
\cite{Carter}
allows the geodesic equations to be integrated, and the geometry of the Kerr spacetime to be completely analysed. The Killing tensor associated to Carter's constant \cite{PenroseWalker:1970} can be used to construct a second-order symmetry operator (i.e. an operator which takes solutions to solutions)  for the wave equation on the Kerr spacetime \cite{Carter:KillingTensor}. 
This symmetry operator, which is a hidden symmetry in the sense that it is not reducible to first-order symmetries, will play a central role in this paper. 

The Kerr black-hole spacetime is expected to be the unique,
stationary, asymptotically flat, vacuum spacetime containing a
nondegenerate Killing horizon
\cite{AlexakisKlainermanIonescu:Uniqueness}.  Further, from
considerations including the weak cosmic censorship conjecture, the
asymptotic limit of the evolution of asymptotically flat, vacuum data
in general relativity is expected to be decomposable into regions,
each of which approaches a Kerr black hole.  An important step towards
establishing the validity of this scenario is to prove the black-hole stability
conjecture, i.e.\ to show that vacuum spacetimes evolving
from data which represent a small perturbation of Kerr initial data with $|a| < M$ 
asymptotically approach a Kerr solution. 

During the past 15 years, there has been a considerable amount of work by several groups towards constructing uniformly bounded energies and proving Morawetz estimates for solutions of the wave equation on black-hole spacetimes. This activity has been motivated by that fact that proving boundedness and decay in time for solutions to the scalar wave equation on the asymptotically flat exterior of the Kerr spacetime is an important model problem for the full black-hole stability problem. 

The basic mechanism for decay of waves 
is by dispersion,
% of energy to asymptotic regions, 
an effect which manifests itself by decay of local energy. The integrated local energy estimate of Morawetz, which captures this effect, was first proved for waves propagating in the exterior of an
obstacle in Minkowski space \cite{morawetz:1961:MR0132908}. 
Both the multiplier method, which was used in the original proof of the Morawetz estimate, and its generalisation, the vector-field method of Klainerman \cite{Klainerman:LorentzGenerators},
provide flexible tools to construct 
%higher-order 
energy estimates for solutions of the wave equation.  

Perhaps the most important application of the vector-field method to date is the monumental proof of the nonlinear stability of Minkowski space \cite{ChristodoulouKlainerman:MinkowskiStability}, one of the central results in general relativity. A partial result had previously been proved using the conformal method \cite{Friedrich:MinkowskiStability}. 
More recently, a simpler proof of the nonlinear stability of Minkowski space 
has been developed \cite{LindbladRodnianski}; however, it also makes use
of the vector-field method. 

%In this paper, we introduce a generalisation of the vector-field method 
%of Klainerman to prove a Morawetz estimate for solutions of the wave
%equation on the exterior of a slowly rotating Kerr black hole. This
%estimate allows us to construct a uniformly bounded energy purely in
%terms of local differential operators.  Our generalisation allows one
%to construct energy identities related not only to first-order Killing
%symmetries but also to the second-order symmetry operator
%corresponding to Carter’s constant.

The fundamental difficulty in proving the existence of a conserved,
positive-definite energy and a Morawetz estimate in the exterior of a
rotating black hole is that neither claim is true without some
adjustment to accommodate the effects of superradiance and trapping.
There is no globally timelike Killing vector field on the exterior of
the Kerr black hole, and hence there is no exactly conserved positive
definite energy for the wave equation. Wave packets which enter the
vicinity of the black hole can leave with higher energy than they had
upon entering \cite{Starobinskii:1973JETP...37...28S}. This phenomenon is called superradiance.  Further, the
Kerr geometry exhibits trapping in the form of null geodesics,
corresponding to light rays, which orbit the black hole. The trapping
in Kerr is complicated in the sense that the orbiting null geodesics
fill an open region of spacetime.  Added to these difficulties is the
fact, mentioned above, that the isometry group of the Kerr geometry is
only 2-dimensional, which should be compared to the Schwarzschild
spacetime with a 4-dimensional isometry group and the flat Minkowski
spacetime with a 10-dimensional group of isometries.

The available results for the wave equation on the Kerr spacetime all make use of Fourier transforms or pseudo-differential operators to overcome the difficulties related to complicated trapping, superradiance, and lack of symmetries. 
The existence of a uniformly bounded, positive-definite energy in the exterior of a slowly rotating Kerr black hole was first proved in \cite{DafermosRodnianski:KerrEnergyBound}, and the Morawetz estimate was first proven in \cite{TataruTohaneanu}. 

In the context of this paper, the most serious difficulty posed by the Kerr geometry is the existence of trapped null geodesics, which orbit the black hole. 
%As mentioned above, the trapped orbits fill an open region in phase space. However, the trapped set has positive codimension in phase space, and the trapped orbits are unstable. 
From \cite{ralston:MR0254433}, one expects that it should be possible to construct solutions of the wave equation for which all but an arbitrarily small amount of energy remains arbitrarily close for an arbitrarily long period of time to a chosen null geodesic. In fact 
this holds in the general Lorentzian setting, not merely one where the metric has a product structure, cf. \cite{sbierski:2013arXiv1311.2477S}. By taking the length of time to be large compared to the constant %larger than the constant 
in any 
%putative 
supposed Morawetz estimate with a spacetime energy density which does not degenerate at the trapped set, one can construct a counterexample to it.
The resolution of this problem is to allow the Morawetz estimate to degenerate at the orbiting null geodesics.  

In the Schwarzschild subclass, where the rotation speed of the black hole vanishes, 
the orbits are all located at Schwarzschild radius $r=3M$, and are unstable. This allows a Morawetz estimate to be proven using a radial vector field 
$\vecfont{A}=\mathcal{F}\partial_r$ where $\mathcal{F}$, 
see \cite{BlueSoffer:ODE,BlueSoffer:LongPaper,BlueSterbenz,DafermosRodnianski:RedShiftSchwarzschild}, 
is continuously differentiable and changes sign at $r=3M$, so that $\vecfont{A}$ points away from the trapped orbits, with $\mathcal{F}$ vanishing to first order there.
Away from the trapped null orbits, the bulk term in the resulting Morawetz estimate is a nondegenerate quadratic expression in the field and its first derivatives. However, the terms involving time and angular derivatives effectively contain a factor $\mathcal{F}^2$ which vanishes quadratically at the trapped set. 

Outside a rotating Kerr black hole, the orbits have a significantly more complicated structure and fill an open set in spacetime. 
%exist on a range of $r$ value. 
Nonetheless, they remain unstable and the trapped set has nonzero codimension in phase space. 
%Further, the special features of the Kerr geometry allow the trapped orbits to be completely analyzed.  
%can be treated by, an albeit subtle form of, separation of variables. 
The primary difficulty then is to construct a smooth vector field which points away from the orbiting null geodesics and vanishes linearly there.
%the tangent space of geodesics, 
It is relatively easy to construct a function on phase space which vanishes linearly on the orbiting null geodesics. The key idea of \cite{DafermosRodnianski:KerrEnergyBound,TataruTohaneanu} is to start with such a function $\mathcal{F}$, to replace the phase space coordinates for null geodesics by 
their conjugate Fourier or spectral variables, and to use the the resulting $\mathcal{F}$ in $\vecfont{A}$.

In our approach, we replace the \emph{conserved quantities for null geodesics} (which may be viewed as coordinates on phase space) by \emph{partial differential symmetry operators} of up to 
%\mymnote{LA: comment that this is similar to the ``separation of variables'' methods, but NOT the ``wild'' Fourier-in-all-directions approaches}
second order. 
This 
%is the fundamentally new idea introduced in this work, and 
%leads to a generalisation
allows us to introduce a generalisation 
of the vector-field method which allows the use of not only Killing
symmetries but also the hidden symmetry corresponding to Carter's constant 
in the construction of 
suitable generalisations of Noetherian 
currents for the analysis of Lagrangian field
equations. The generalised vector-field method allows us, in contrast to other recent work on the
wave equation on Kerr, to carry out our proof of uniform boundedness
and integrated energy decay exclusively in physical space, using only
the coordinate functions and differential operators. This technique
eliminates the need for methods involving separation of variables or
Fourier analysis.  The suitability of the classical vector-field method for nonlinear problems, which was demonstrated in the proof of the nonlinear stability of Minkowski space, partly motivates, in view of the black hole stability problem, our work on generalising the vector-field method to deal with the linear wave equation in the Kerr spacetime. 

In our proof of the Morawetz estimate, we use a radial vector field
with a coefficient function $\mathcal{F}$, which takes values in the algebra generated by the symmetry operators of second order. 
This function $\mathcal F$ is constructed from the radial derivatives of coefficients appearing in the wave equation. 
%In particular, this construction makes direct use of the hidden symmetry of the Kerr spacetime which manifests itself in the existence of Carter constant and its corresponding symmetry operator. 
As in the Schwarzschild case discussed above, there is a quadratic degeneracy in terms involving time and angular derivatives, which again occurs where the corresponding terms in $\mathcal F$ vanish. 
%A more detailed discussion can be found later in this introduction, see section \ref{SS:Strategy}. An explicit form of the bulk term in the Morawetz estimate is given in \eqref{eq:BigImportantNotNumberedMorawetz}. The quadratically degenerate term occurs on the first line in this equation. 
This is discussed in more detail following formula \eqref{eq:BigImportantNotNumberedMorawetz}. It is in the first line of this formula that the quadratic degeneracy can be most clearly seen.
We remark that using Fourier techniques, it is possible to strengthen the local energy estimate at the trapped set. In \cite{MarzuolaMetcalfeTataruTohaneanu} a local energy estimate with logarithmic losses at the trapped set was proved for the Schwarzschild case, and in \cite{tohaneanu:strichartz:kerr:MR2846348} for Kerr. Further, it is possible to prove a local energy estimate controlling a fractional energy norm with a fractional loss of derivatives, but which is uniform at the trapped set, cf. \cite{andersson:blue:maxwell:2013arXiv1310.2664A,andersson:blue:nicolas:MR3021792,BlueSoffer:LongPaper}.

The solution to the problem posed by the lack of an exactly conserved
positive definite energy, and the related superradiance phenomenon, is
to construct a vector field that is globally timelike, approaches the
generator of time-translations $\partial_t$ in the asymptotically flat
region, and for which a Morawetz estimate can be used to control the
change, from one time to another, in the energy associated with the
new vector field. Because of the dominant energy condition, this
vector field provides a positive-definite energy, which, using the
Morawetz estimate, can be bounded uniformly in terms of its initial
value, although it is not conserved. This is the method used in this
paper and \cite{TataruTohaneanu}. In
\cite{DafermosRodnianski:KerrEnergyBound}, solutions to the wave
equation are decomposed into various frequency regimes, and, to
construct a uniformly bounded energy, it is not necessary to prove a
Morawetz estimate in all frequency regimes. 
%This avoidance of a global Morawetz estimate may have some advantages in approaching the rapidly rotating case \cite{Dafermos:Rodnianski:2010maph.conf..421D,shlapentokh-rothman:2013arXiv1302.6902S}. 

During the preparation of this paper, uniform energy and Morawetz
estimates have been obtained for the full range $|a|<M$
\cite{dafermos:rodnianski:Shlapentokh-Rothman:2014arXiv1402.7034}. This
builds on
\cite{shlapentokh-rothman:QuantitativeModeStability}
as well as
\cite{DafermosRodnianski:KerrEnergyBound} and related works.

\newpageForSubsection
%%%%%%%%%%%%%%%%%%%%%%%%%%%%%%%%%%%%%%%%%%%%%%%%%%%%%%%%%%%%%%%%%%%%%%%%%%
\subsection{Hidden symmetries and null geodesics of the Kerr spacetime}
\label{SS:SymNull} 
We begin this subsection with a brief discussion of conserved quantities for null geodesics and of symmetry operators for the wave equation. We then review how conserved quantities are used to analyse the null geodesics in the Kerr spacetime. Finally, we review the close connections in the Kerr spacetime between the analysis of the null geodesics and of the wave equation, particularly its hidden symmetries.

For an affinely parametrised geodesic $\Geodesic^\alpha$ with velocity $\GeodesicVelocity^\alpha$, the mass squared,
%Hamiltonian $\GeodesicH$ given by $\GeodesicH^2 = - \gMetric_{\alpha\beta}
%\GeodesicVelocity^\alpha \GeodesicVelocity^\beta$ 
$-\gMetric_{\alpha\beta}
\GeodesicVelocity^\alpha \GeodesicVelocity^\beta$, 
is a constant of the
motion. For null geodesics, the mass is zero. 
For a coordinate, e.g.{} $t$, let $\GeodesicVelocity_t$ be the $t$-momentum given by 
$\GeodesicVelocity_t = \gMetric_{\alpha\beta} \GeodesicVelocity^\alpha
(\partial_t)^\beta $.
In the Kerr
spacetime, we have the Killing fields $\dt$ and $\dphi$, and the associated momenta are then constants of the motion. For a null geodesic $\Geodesic^\alpha$, we define the energy and the azimuthal angular momentum to be $\GeodesicEnergy=-\dot\gamma_t$ and
$\GeodesicLz=-\dot\gamma_\phi$ respectively.\footnote{The sign for $\GeodesicEnergy$ is chosen so that it is positive for future-directed null geodesics for sufficiently large $r$. The sign for $\GeodesicLz$ is chosen so that the sign convention is consistent with that of $\GeodesicEnergy$.}

A symmetric 2-tensor $\TensorK_{\alpha\beta}$ is called a  Killing tensor if $\nabla_{(\alpha} \TensorK_{\beta\gamma)} = 0$, cf. \cite{PenroseWalker:1970}.  
If $\TensorK_{\alpha\beta}$ is a Killing tensor, then for any affinely parametrised geodesic $\Geodesic^\alpha$, the quantity $\TensorK_{\alpha\beta} \GeodesicVelocity^{\alpha}
\GeodesicVelocity^{\beta}$ is a conserved quantity. 
As shown by Carter \cite{Carter:KillingTensor}, if $\TensorK_{\alpha\beta}$ is a Killing tensor in a vacuum spacetime, 
then the second-order operator $\nabla_\alpha \TensorK^{\alpha\beta} \nabla_\beta$ satisfies 
$[\nabla_\alpha \TensorK^{\alpha\beta} \nabla_\beta,\nabla^\gamma\nabla_\gamma] = 0$, i.e.{} $\nabla_\alpha \TensorK^{\alpha\beta} \nabla_\beta$ commutes with the d'Alembertian. 

We take a symmetry operator for the wave equation, $\nabla^\alpha\nabla_\alpha\solu=0$, to be a differential operator that takes solutions to solutions. Recall that if $\vecX$ is a Killing field, then the
operator $\Lie_\vecX$ generated by Lie differentiation with respect to
$\vecX$ is a symmetry operator. From the previous paragraph, if $\TensorK^{\alpha\beta}$ is a Killing tensor, then $\nabla_\alpha \TensorK^{\alpha\beta} \nabla_\beta$ is also a symmetry operator. The set of symmetry operators 
for the wave equation is closed
under scalar multiplication, addition, and composition, and each
symmetry operator has a well-defined order as a differential operator.
Thus, the set of symmetry operators forms a graded algebra. 

A hidden symmetry is defined to be a symmetry operator which is not in
the algebra generated by the Killing vector fields and the
d'Alembertian. In Minkowski spacetime, since the
Delong-Takeuchi-Thompson inequality is saturated, there are no hidden
symmetries \cite{ChanuEtAl}. Similarly, in the Schwarzschild
spacetime, there are no hidden symmetries \cite{Caviglia}. In
contrast, the Kerr spacetime admits a Killing 2-tensor
\cite{PenroseWalker:1970}, $\TensorK^{\alpha\beta}$, which generates
both Carter's constant \cite{Carter}
$\KCarter=\TensorKCarter_{\alpha\beta}\dot\gamma^\alpha\dot\gamma^\beta$
and a second-order symmetry operator; this gives a hidden symmetry in
the Kerr spacetime.

In the Kerr spacetime, Carter's constant provides a fourth constant of the motion, in addition to the mass, energy, and azimuthal angular momentum. Here we shall be interested in the expression 
\begin{align}
%\GeodesicQ = \GeodesicMomentum_{\theta}^2 +\frac{\cos^2\theta}{\sin^2\theta}\GeodesicMomentum_\phi^2 +a^2\sin^2\theta\GeodesicMomentum_t^2 .
\GeodesicQ = \GeodesicVelocity_{\theta}^2 +\frac{\cos^2\theta}{\sin^2\theta}\GeodesicVelocity_\phi^2 +a^2\sin^2\theta\GeodesicVelocity_t^2 .
\label{eq:IntroDefGeodesicQ}
\end{align}
%which can be written in the form
%\mnote{LA: introduced macro {\tt $\backslash$TensorQ} }
%$\GeodesicQ = \TensorQ_{\alpha\beta}
%\GeodesicVelocity^\alpha \GeodesicVelocity^\beta$. 
%The symmetric tensor $\TensorQ_{\alpha\beta}$ given here is a conformal Killing tensor.
The quantity $\GeodesicQ$ is closely related to the commonly used form of Carter's constant $\KCarter$, see appendix 
\ref{S:CarterOperatorAndQ} for details.  
For null geodesics, we have 
$\GeodesicQ = \KCarter - 2a \GeodesicEnergy \GeodesicLz -\GeodesicLz^2$.
Any
linear combination of $\GeodesicEnergy^2$,
$\GeodesicEnergy\GeodesicLz$, and $\GeodesicLz^2$ can be added to
$\GeodesicQ$ to give an alternate choice for the fourth constant of
the motion for null geodesics. 
The form we have chosen is uncommon,
but useful for our purposes because it has nonnegative coefficients. 

%As was demonstrated by Carter, 
The presence of the extra conserved
quantity allows one to 
%separate 
integrate
the equations of geodesic motion.  Of
most interest to us is the equation for the $r$-coordinate of null
geodesics\cite{FrolovNovikov},
\begin{align}
\KSigma^2 \left ( \frac{\di r}{\di \lambda} \right )^2 ={}&{}
-\CurlyR(r;M,a;\GeodesicEnergy,\GeodesicLz,\GeodesicQ) , 
\label{eq:IntroNullGeodesicsr}\\
\intertext{where $\lambda$ is the affine parameter of the null geodesic and }
\CurlyR(r;M,a;\GeodesicEnergy,\GeodesicLz,\GeodesicQ)
={}&{} -(r^2+a^2)^2\GeodesicEnergy^2 -4aMr\GeodesicEnergy\GeodesicLz
+(\KDelta-a^2)\GeodesicLz^2 +\KDelta\GeodesicQ . 
\label{eq:RR:consquant} 
\end{align} 
One finds that there are null geodesics for which the $r$ coordinate
is constant. We refer to these as orbiting 
null geodesics. The $r$-values for which this is possible are the solutions to the equations
\begin{align}
\CurlyR ={}&{} 0,\quad \partial_r \CurlyR = 0.
\label{eq:RotatingGeodesicCondition}
\end{align} 
The corresponding
null geodesics are unstable, which with our conventions corresponds to
$\partial_r^2 \CurlyR < 0$. The orbiting null geodesics are
the only ones which neither go to nor have come from infinity or the horizon. 
There are
other null geodesics that fail to fall into the black hole or 
escape to infinity, but the $r$ coordinate along these asymptotically
approaches (towards the future) an $r$ value for which there is an
orbiting null geodesic.

In the Schwarzschild case, i.e.\ for $a=0$, there are only orbiting
null geodesics on the sphere at $r=3M$, which is called the photon
sphere. For nonzero $a$, the orbiting null geodesics fill up an open
region in spacetime. 
As $a \to 0$, this region tends to
$r=3M$. There are many descriptions of the Kerr spacetime and its
geodesics, including \cite{Bardeen,FrolovNovikov,Teo}.

In Boyer-Lindquist coordinates, the d'Alembertian $\gWave =
\nabla^\alpha \nabla_\alpha$  takes the form 
\begin{align}
\gWave = \frac{1}{\KSigma} \left(\dr\KDelta\dr +\frac1\KDelta\CurlyR(r;M,a;\dt,\dphi,\OpQ)\right) ,
\label{eq:KerrWaveBLCoord}
\end{align}
where $\CurlyR$ is given by \eqref{eq:RR:consquant} with the conserved
quantities $\GeodesicEnergy,\GeodesicLz,\GeodesicQ$ replaced by their corresponding operators $\dt, \dphi$, and the second-order 
%Carter 
operator 
$\OpQ$ given by
\begin{equation}\label{eq:OpQdef}
\OpQ = \frac1{\sin\theta}\dtheta\sin\theta\dtheta
+\frac{\cos^2\theta}{\sin^2\theta}\dphi^2 +a^2\sin^2\theta \dt^2 .
\end{equation} 
In appendix \ref{S:CarterOperatorAndQ}, we explain the relationship between  $\OpQ$ and $\nabla_\alpha \TensorK^{\alpha\beta}\nabla_\beta$. The operator $\CurlyR$ is given by 
\begin{equation} \label{eq:OpCurlyR}
\CurlyR(r;M,a;\dt,\dphi,\OpQ)
= -(r^2+a^2)^2\dt^2 -4aMr\dt\dphi
+(\KDelta-a^2)\dphi^2 +\KDelta\OpQ . 
\end{equation}
We have used some unconventional 
sign conventions in defining $\CurlyR$ to avoid having to use 
factors of $i$ when replacing the constants of motion by differential operators.  
It is clear from the above that $\partial_t$, $\partial_\phi$, and
$\OpQ$ are symmetry operators for the wave equation on Kerr. 
In fact, we see from \eqref{eq:KerrWaveBLCoord} that the operator $\OpQ$ \emph{commutes} with the operator $\KSigma\gWave$, and in particular, that $\OpQ$ is a symmetry operator for $\gWave$. 
The operators $\nabla_\alpha\TensorK^{\alpha\beta}\nabla_\beta$ and $\OpQ$ are both hidden symmetries. 
%Since the Carter Killing tensor, as well as the conformal Killing tensor $\TensorQ_{\alpha\beta}$ cannot be written as a linear combination of tensor products of Killing fields and the metric, the operators $\OpK$ and $\OpQ$ are hidden symmetries. 
We remark that the operator $\OpQ$ is closely related to the angular operator in the spin-0 Teukolsky system, cf. \cite{Teukolsky}.

We denote 
the set of order-$n$ generators of the symmetry algebra 
generated by $\dt$, $\dphi$, and $\OpQ$ by
\begin{align}
\SymGeneraln{n}= \{ \dt^\nt \dphi^\np \OpQ^\nq | \nt +\np +2\nq = n; \nt,\np,\nq\in\Naturals \} .
\label{eq:KerrHigherSymmetries}
\end{align}
In particular, 
$$
\SymGeneraln{0}=\{\Id\} , \quad \SymGeneraln{1}=\{\dt,\dphi\} . 
%\SymGeneraln{2}={}&{}\{\dt^2,\dt\dphi,\dphi^2,\OpQ\} .
$$
Of particular importance in our analysis will be the set of second-order symmetry operators, 
\begin{align*}
\SOSym =\{ \dt^2, \dt\dphi, \dphi^2, \OpQ \}=\{\CQA_\ua \} ,
\end{align*}
and 
underlined
%roman 
indices always refer to the index in this set. 
Higher-order pointwise norms are defined 
in terms of $\SymGeneraln{n}$
by 
\begin{equation}\label{eq:ptwise} 
\normPtwiseTn{n}{\solu}^2=
\sum_{j=0}^n\sum_{\CQA\in\SymGeneraln{j}} |\CQA\solu|^2. 
\end{equation} 
The function $\CurlyR$ is polynomial in its last three arguments, so
$\CurlyR(r;M,a;\dt,\dphi,\OpQ)$ is well defined. Furthermore, it can
be written as a linear combination of the second-order symmetries with
coefficients which are rational in $r$, $M$, and $a$.
Such linear combinations of second-order symmetry operators play a crucial role in the analysis of this paper. Having introduced underlined subscript indices for the second-order symmetry operators, we can introduce underlined superscript indices for the rational functions in $r$, $M$, and $a$ that are the coefficients. Thus, using the standard Einstein index convention in the underlined indices, we can write these linear combinations as 
\begin{align*}
\CurlyR(r;M,a;\dt,\dphi,\OpQ)={}&{} \CurlyR^\ua \CQA_\ua .
\end{align*}

\newpageForSubsection
%%%%%%%%%%%%%%%%%%%%%%%%%%%%%%%%%%%%%%%%%%%%%%%%%%%%%%%%%%%%%%%%%%%%%%%%%%
\subsection{Statement of results}
\label{SS:Statement} 
We now state our main results and briefly compare them with previous results.
In formulating our estimates, we shall make use of the following model energy, 
\begin{align*} 
\ModelEnergyThree[\solu](\hst{t})\\
%=\sum_{j=0}^2\sum_{\CQA\in\SymGeneraln{2}}\int_{\hst{t}}{}&{} \left(\frac{(r^2+a^2)^2}{\KDelta}|\dt\CQA\solu|^2 +\KDelta|\dr\CQA\solu|^2+|\dtheta\CQA\solu|^2+\frac{1}{\sin^2\theta}|\dphi\CQA\solu|^2\right)\diThree ,
=\int_{\hst{t}}{}&{} \left(\frac{(r^2+a^2)^2}{\KDelta}|\dt\solu|_2^2 +\KDelta|\dr\solu|_2^2+|\dtheta\solu|_2^2+\frac{1}{\sin^2\theta}|\dphi\solu|_2^2\right)\diThree , 
\end{align*}
where $|\cdot|_2$ is the second-order point-wise norm introduced in \eqref{eq:ptwise} above, and 
$
\diThree=\sin\theta\di r\di\theta\di\phi  
$
is a reference volume element on the Cauchy slice $\hst{t}$.

As discussed above, the main contribution of this paper is a new method to prove the
following results, which had previously been known  
  from \cite{DafermosRodnianski:KerrEnergyBound}
and \cite{TataruTohaneanu}:  
\begin{theorem}[Uniformly bounded, positive energy]
\label{Thm:IntroUniformEnergyBound}
Given $M>0$, there are positive constants $\ConsantIntroUniformBound$
and $\epsilonSlowRotationIntroUniformBound$, such that, if
$|a|\leq\epsilonSlowRotationIntroUniformBound$ and
$\solu:\Reals\times(\rp,\infty)\times S^2\rightarrow\Reals$ is a
solution of the wave equation, $\nabla^\alpha\nabla_\alpha\solu=0$,
then $\forall t$
\begin{align*}
\ModelEnergyThree(\hst{t}) \leq \ConsantIntroUniformBound \ModelEnergyThree(\hst{0}) .
%\GenEnergyOrder{\vecTBlend}{3}[\solu](\hst{t}) \leq \ConsantIntroUniformBound \GenEnergyOrder{\vecTBlend}{3}[\solu](\hst{0}) .
\end{align*}
\end{theorem}

\begin{theorem}[Morawetz estimate]
\label{Thm:IntroMorawetz}
Given $M>0$, there are positive constants
$\epsilonMorawetzSlowRotation$, 
$\epsilonMorawetzPhotonWidth$, $\epsilonInvIntegratedMorawetz$ and a function $\localiseAwayFromPhotonOrbits$ which is identically one for $|r-3M|>\epsilonMorawetzPhotonWidth$ and zero otherwise, 
such that for all $|a|\leq\epsilonSlowRotationIntroDecay$ and all smooth
$\solu$ solving the wave equation, 
$\nabla^\alpha\nabla_\alpha\solu=0$,
%$\gWave \solu = 0$
\begin{align*}
\int_{-\infty}^{\infty} \int_{\rp}^{\infty} \int_{S^2}
 \bigg( \left(\frac{\KDelta^2}{r^4}\right) \normPtwiseTn{2}{\dr \solu}^2 +r^{-2}\normPtwiseTn{2}{\solu}^2  
+ \localiseAwayFromPhotonOrbits \frac{1}{r} \left(\normPtwiseTn{2}{\dt \solu}^2+\normPtwiseTn{2}{\dAng \solu}^2\right) \bigg)
\diFour  \\
%\nonumber\\
\leq\epsilonInvIntegratedMorawetz \ModelEnergyThree(\hst{0}) ,
%\label{eq:IntegratedMorawetz} 
\end{align*}
where $\dAng$ is the angular gradient in Boyer-Lindquist coordinates and $\diFour = \diThree \di t$.
\end{theorem}
Theorem \ref{Thm:IntroUniformEnergyBound} is the conclusion of section
\ref{S:BoundedEnergy}, cf.{} theorem
\ref{Thm:UniformEnergyBound}. Theorem \ref{Thm:IntroMorawetz} follows from lemma
\ref{Lemma:IntegratedMorawetz} and the uniform bound in theorem
\ref{Thm:IntroUniformEnergyBound}. We remark that Lemma \ref{Lemma:IntegratedMorawetz} gives an estimate of the integrated Morawetz density appearing in Theorem 1.2 in terms of the initial and final energy which may be of independent interest. 

The degeneracy near $r=3M$ in the Morawetz estimate of
theorem \ref{Thm:IntroMorawetz} is not optimal. In fact, a sharper
estimate follows from lemma \ref{Lemma:Morawetz}. 
However, some
degeneracy in this region is unavoidable, due to the existence in the
Kerr spacetime of null geodesics which orbit at fixed $r$
values. This degeneracy is discussed further in section
\ref{SS:Strategy}. 

The energy $\ModelEnergyThree$ is both ad hoc and degenerate as $r\rightarrow\rp$. In subsection \ref{SS:TBlend}, we relate it to a geometrically defined energy $\GenEnergyOrder{\vecTBlend}{3}$. The degeneracy as $r\rightarrow\rp$ is concealed by the Boyer-Lindquist coordinates and can even appear, in the coefficient of $|\dr\solu|^2$, as a divergence. In appendix \ref{S:DRVector}, we use a nondegenerate coordinate system and apply the
ideas of
\cite{DafermosRodnianski:RedShiftSchwarzschild,DafermosRodnianski:LectureNotes},
to overcome the degeneracy in the energy. (There is a similar degeneracy in theorem
\ref{Thm:IntroMorawetz}, which is also removed. See equation
\eqref{eq:NonDegenerateMorawetzBulk} and the subsequent discussion.) On $\hst{t}$, there is a nondegenerate energy $\GenEnergy{\normalhst{t},3}$  
which is equivalent to the sum of the $L^2$ norms of all derivatives of order $1$ through $3$, with derivatives taken in spatial directions tangential to $\hst{0}$ and in the timelike direction orthogonal to this surface. This allows us to control nondegenerate, third-order Sobolev norms on a new foliation, from which we can obtain the following $L^\infty$ estimate:

\begin{corollary}[Uniformly bounded solution]
\label{Cor:IntroUniformSoluBound}
Given $M>0$, there are positive constants
$\ConsantIntroUniformBoundsolu$ and
$\epsilonSlowRotationIntroUniformBoundSolu$, and a nonnegative
quadratic form $\GenEnergy{\normalhst{0},3}$ on $\hst{0}$, such
that, if $|a|<\epsilonSlowRotationIntroUniformBoundSolu$ and
$\solu:\Reals\times(\rp,\infty)\times S^2\rightarrow\Reals$ is a
solution of the wave equation, $\nabla^\alpha\nabla_\alpha\solu=0$,
then $\forall (t,r,\omega)\in\Reals\times(\rp,\infty)\times S^2$
\begin{align*}
|\solu(t,r,\omega)|
\leq \ConsantIntroUniformBoundsolu \GenEnergyOrder{\normalhst{0}}{3}[\solu](\hst{0})^{1/2} .
\end{align*}
\end{corollary}

\begin{remark} 
Corollary
  \ref{Cor:IntroUniformSoluBound} and an analogue of theorem
  \ref{Thm:IntroUniformEnergyBound} were first proven in
  \cite{DafermosRodnianski:KerrEnergyBound}, and an analogue of
  theorem \ref{Thm:IntroMorawetz} was first proven in
  \cite{TataruTohaneanu}. Both works deal directly with
  $\GenEnergy{\normalhst{t}}$ without estimating
  $\GenEnergy{\vecTBlend}$. Since our focus is to control the
  influence of the orbiting null geodesics, which are relatively far
  from $r=\rp$, we find it convenient to work with the weaker norms
  $\GenEnergy{\vecTBlend}^{1/2}$ which are dominated by
  $\GenEnergy{\normalhst{t}}^{1/2}$.  Appendix \ref{S:DRVector}
  provides the details in removing the degeneracy at
  $r\rightarrow\rp$. 
\end{remark} 

The quadratic forms $\ModelEnergyThree(\hst{0})$ and
$\GenEnergyOrder{\normalhst{0}}{3}[\psi](\hst{0})$ 
are bounded if $\solu$
extends smoothly to the closure of the hypersurface $\hst{0}$ in the
extended spacetime and satisfy the following decay conditions. First,
as $r\rightarrow\infty$, it is sufficient that $\solu$ and its first
three derivatives (with respect to the Boyer-Lindquist coordinates)
decay like $r^{-3/2+\delta}$ for some positive $\delta$. Since this
decay rate is stated with respect to the angular derivatives, a more
geometric statement of this decay condition is that the first three
normalised derivatives in the angular directions decay at a rate with
one additional power of $r^{-1}$ for each additional angular
derivative. By ``smooth'' we mean $C^\infty$ with respect to local
coordinates in the extended spacetime.  As $r\rightarrow\rp$, this is
not the same as being smooth with respect to the Boyer-Lindquist
coordinates. The
$\GenEnergyOrder{\vecTBlend}{3}[\psi](\hst{0})^{1/2}$ and
$\GenEnergyOrder{\normalhst{0}}{3}[\psi](\hst{0})^{1/2}$ 
are $L^2$-based
weighted Sobolev norms, so these $L^\infty$ decay conditions on the
initial data are sufficient, but not necessary.

\subsection{Summary of previous results}
\label{SS:PreviousResults} 
As mentioned above, 
our work builds on previous results in the subcase of the
Schwarzschild spacetime, where $a=0$. This subcase
%
%We briefly comment on other related work. Estimates for the decay rate
%of solutions to the wave equation have been proven in the subcase of
%the Schwarzschild spacetime, where $a=0$. Birkhoff's theorem states
%that the Schwarzschild spacetime is the unique spherically symmetric,
%vacuum spacetime solution of Einstein's equation. For the coupled
%Einstein and scalar wave system, a decay rate and nonlinear stability
%of the Schwarzschild solution have been proven in the spherically
%symmetric setting \cite{DafermosRodnianski:PriceLaw}.
%
%For the wave equation without a symmetry
%assumption but on a fixed background spacetime, the case of the linear
%wave equation on the Schwarzschild spacetime 
is significantly simpler than the situation in the Kerr
spacetime, since the $\dt$ Killing vector is timelike in the entire
exterior region and generates a conserved positive energy, there is
the full $SO(3)$ group of rotation symmetries available to generate
higher energies, and the orbiting null geodesics are restricted to
$r=3M$.  The first two of these properties imply theorem
\ref{Thm:IntroUniformEnergyBound} in the $a=0$ case.  Following the
introduction of a Morawetz vector field and of the equivalent of an
almost conformal vector field to the Schwarzschild spacetime
\cite{LabaSoffer}, decay estimates for the wave equation were proven
\cite{BlueSterbenz}, proven with a weaker decay rate but less
regularity loss \cite{BlueSoffer:LongPaper}, and proven separately
with a stronger estimate near the event horizon
\cite{DafermosRodnianski:RedShiftSchwarzschild}. These were extended
to Strichartz estimates for the wave equation
\cite{MarzuolaMetcalfeTataruTohaneanu} and to decay estimates for
Maxwell's equation \cite{Blue:Maxwell} and for Einstein's equation on
ultimately Schwarzschildean spacetimes \cite{Holzegel}. All of these
works relied upon a Morawetz estimate. The Morawetz vector field which
made these estimates possible was centred about the orbiting geodesics
at $r=3M$.

This construction of a classical Morawetz vector field fails when
$a\not=0$, since there are orbiting null geodesics filling an open set in
spacetime. Using Fourier-analytic techniques, one may define
generalised Morawetz vector fields which circumvent this problem
\cite{DafermosRodnianski:KerrEnergyBound,TataruTohaneanu}. 
These
Fourier-analytic Morawetz vector fields have coefficients that depend
on both spacetime position and on the Fourier variables
conjugate to the spacetime coordinates. One
  advantage of such an approach is that it allows the analogues of
  theorems \ref{Thm:IntroUniformEnergyBound} and
  \ref{Thm:IntroMorawetz} to be proven in $H^1$ regularity Sobolev
  spaces, instead of $H^3$ regularity Sobolev spaces. At the end
of the introduction, 
we further compare the method of those works with
the present paper. The techniques in this paper remain the
only ones to provide a Morawetz or integrated energy estimate without
using a Fourier variable conjugate to time.

%The techniques used in these papers may 
%reasonably be called Fourier-analytic, pseudo-differential, microlocal,
%or phase-space techniques, since the Fourier operators represent
%coordinates in momentum space in contrast to spacetime coordinates in
%physical or configuration space. These results include a form of weak
%decay, since the Morawetz estimate implies integrability of the local
%energy. 
%
A number of stronger decay estimates have been obtained by
building upon a Morawetz estimate in the slowly rotating
case. These include a pointwise decay estimate with rate
$t^{-1+\epsilon}$ \cite{DafermosRodnianski:LectureNotes}, a pointwise
in time decay with rate $t^{-3/2+\epsilon}$
\cite{Luk:Kerr:improveddecay:2010arXiv1009.0671L} (see also
\cite{Dafermos:Rodnianski:2010maph.conf..421D}), a Strichartz estimate
\cite{tohaneanu:strichartz:kerr:MR2846348}, 
%\cite{Tohaneanu:Kerr:Strichartz:2009arXiv0910.1545T}, and pointwise
decay with rate $t^{-3}$ \cite{Tataru:Price}. In these statements
$\epsilon$ refers to a continuous function of $a$ which vanishes as
$|a|\searrow0$. These pointwise in time decay estimates hold in
stationary regions, where $r$ is bounded by $\rp<r_1<r<r_2<\infty$,
but, in all cases, there are corresponding decay estimates along the
boundary of the exterior region, i.e.\ the event horizon $r=\rp$ and
null infinity $r\rightarrow\infty$. The recent result covering the
full range $|a|<M$ also includes a pointwise decay estimate with rate
$t^{-3/2+\epsilon}$ 
\cite{dafermos:rodnianski:Shlapentokh-Rothman:2014arXiv1402.7034}.

There have also been several lines of work that do not make use of a
Morawetz estimate. Fourier-analytic vector fields were used previously
to prove Mourre estimates, which are similar to Morawetz estimates, in
the proof of scattering for the Klein-Gordon equation
\cite{Hafner:KerrKGScattering} and the Dirac equation
\cite{HafnerNicolas:KerrDiracScattering}. The complete separability of
the wave equation in the Kerr spacetime was used to derive an explicit
representation
\cite{FinsterKamranSmollerYau:KerrWaveRepresentation}. For solutions
of the form $\solu(t,r,\theta,\phi)=\solu_{\Lz}(t,r,\theta)
e^{i\Lz\phi}$ or where $\solu$ is made up of a finite number of
azimuthal modes of this form, the integral representation was used to
prove an $L^\infty_{\text{loc}}$ decay result.  Such solutions form a
dense set in the space of solutions, but, without a uniform estimate
on the decay rate, it is not possible to pass to a limit.  Decay rates
have been obtained from this separability method for solutions to the
Dirac equation \cite{FinsterKamranSmollerYau:DiracDecay} and
spherically symmetric solutions to the wave equation when $a=0$
\cite{Kronthaler}.  The decay without rate results for the wave
equation built on two earlier results. The first showed that for
$|a|\in[0,M)$, it is possible to do a full separation of variables
  \cite{Teukolsky}. The second showed that for each hypothetical
  exponentially growing mode there is a conserved, nonnegative,
  energy-like expression, and hence that there are no exponentially
  growing modes \cite{Whiting}. Recently, separation of variables
  techniques have been used to prove a uniform $t^{-3}$ decay rate (in
  stationary regions) for solutions to the wave equation on the
  Schwarzschild background
  \cite{donninger:schlag:soffer:2009arXiv0911.3179D}. The uniform
  $t^{-3}$ rate is the one conjectured by the ``summed Price law'',
  based on Price's prediction that all modes decay at a rate of
  $t^{-3}$ or faster \cite{price:1972:I,price:1972:II}.

For the coupled Einstein and
scalar wave system, a decay rate and nonlinear stability of the
Schwarzschild solution have been proven in the spherically symmetric
setting \cite{DafermosRodnianski:PriceLaw}. Birkhoff's theorem 
states that the Schwarzschild spacetime is the
unique spherically symmetric, vacuum spacetime solution of Einstein's
equation. Hence, to consider a dynamical problem, one must couple the
Einstein equation to some matter model.

\newpageForSubsection
%%%%%%%%%%%%%%%%%%%%%%%%%%%%%%%%%%%%%%%%%%%%%%%%%%%%%%%%%%%%%%%%%%%%%%%%%%%%%
\subsection{A monotonicity result for null geodesics}
%{Hidden symmetries and the vector-field method}
\label{SS:ResultForGeodesics} 

Here we illustrate the key idea of the paper by exploring a related one for null geodesics. At the heart of the Morawetz estimate is a monotonicity formula for the wave equation. Null geodesics are the characteristic curves for the wave equation and provide important insight into the behaviour of solutions of the wave equation. 

For a null geodesic, $\Geodesic(\lambda)$, we define the energy associated with  a vector field $\vecX$ and evaluated on a Cauchy hypersurface $\HypersurfaceGeneral$ to be
\begin{align*}
\GenEnergyGeodesic{\vecX}[\Geodesic](\HypersurfaceGeneral)
{}&{}= -\gMetric_{\alpha\beta}\vecX^\alpha \dot\Geodesic^\beta|_{\HypersurfaceGeneral}. 
\end{align*}
%Since $\HypersurfaceGeneral$ is a Cauchy hypersurface, the null geodesic $\Geodesic$ has a unique intersection with $\HypersurfaceGeneral$, so that $\dot\Geodesic|_{\HypersurfaceGeneral}$ and the energy are well-defined. 
The sign in the energy is introduced so that it is nonnegative if $\vecX$ and $\dot{\Geodesic}$ are both future directed and causal. 
If the spacetime is globally hyperbolic, for each value of the geodesic parameter $\lambda$ there is a unique Cauchy surface $\HypersurfaceGeneral$ for which $\Geodesic(\lambda)$ intersects $\HypersurfaceGeneral$. Differentiation with respect to $\lambda$ is equivalent to differentiation along the tangent vector. Since $\nabla_{\dot\Geodesic}\dot\Geodesic=0$ for a geodesic, integrating the derivative of the energy gives
% the following energy identity
\begin{align}
\GenEnergyGeodesic{\vecX}[\Geodesic](\HypersurfaceGeneral_2)-\GenEnergyGeodesic{\vecX}[\Geodesic](\HypersurfaceGeneral_1)
={}&{} -\int_{\lambda_1}^{\lambda_2} (\dot\Geodesic_\alpha\dot\Geodesic_\beta)(\nabla^{(\alpha}\vecX^{\beta)}) \di\lambda ,
\label{eq:deformForGeodesics}
\end{align}
where $\lambda_i$ is the unique value of $\lambda$ such that $\Geodesic(\lambda)$ is the intersection of $\Geodesic$ with $\HypersurfaceGeneral_i$. The sign arises from our choice of sign in the definition of the energy. Formula \eqref{eq:deformForGeodesics} is particularly easy to work with, if one recalls that 
\begin{align*}
\nabla^{(\alpha}\vecX^{\beta)}{}&{}= -\frac12 \Lie_{\vecX}\gMetric^{\alpha\beta} .
\end{align*}
The tensor $\nabla^{(\alpha}\vecX^{\beta)}$ is commonly called the ``deformation tensor''. In the following, unless there is room for confusion, we will drop reference to $\Geodesic$ and $\HypersurfaceGeneral$ in referring to $\GenEnergyGeodesic{\vecX}$. 

These energies for null geodesics are useful for understanding the monotonicity at the heart of the original Morawetz estimate in the Minkowski spacetime, $\Reals^{1+3}$. That estimate is proven using the radial vector field $\dr$ in $(t,r,\theta,\phi)$ coordinates. In Minkowski spacetime, null geodesics are simply straight lines, and one can consider the projection onto the spatial component in $\Reals^n$, which is also a straight line $\vec{x}(t)$ and can be parametrised by time. The projection of a null geodesic will have a constant and nonvanishing tangent, $\vec{v}(t)$. Asymptotically, the position and velocity will become aligned, so that $\lim_{t\rightarrow\infty} \vec{v}(t)\cdot\vec{x}(t)/|\vec{x}(t)| =1$, and similarly in the past, $\lim_{t\rightarrow-\infty} \vec{v}(t)\cdot\vec{x}(t)/|\vec{x}(t)| =-1$. Thus, the radial component of the velocity, $\dot\Geodesic_r=-\GenEnergyGeodesic{\dr}$, increases overall from the asymptotic past to the asymptotic future. In fact, it is not hard to show that the radial component of the velocity increases monotonically. In particular, with $\eta^{\alpha\beta}=-\dt^\alpha\dt^\beta +\dr^\alpha\dr^\beta +r^{-2}(\dtheta^\alpha\dtheta^\beta+\sin^{-2}\dphi^\alpha\dphi^\beta)$, for null geodesics that do not pass through $r=0$,\footnote{For null geodesics passing through $r=0$ in $\Reals^{1+n}$, there is a singular contribution as the radial component of the velocity instantaneously goes $-1$ to $1$}  one finds the monotonicity formula  $-(\di/\di\lambda)\GenEnergyGeodesic{\dr}=(-1/2)(\Lie_{\dr}\eta^{\alpha\beta})\dot\Geodesic_\alpha\dot\Geodesic_\beta=r^{-3}(\dot\Geodesic_\theta\dot\Geodesic_\theta+\sin^{-2}\dot\Geodesic_\phi\dot\Geodesic_\phi)\geq 0$. 

We now consider the behaviour of the radial velocity of a null geodesic
in the Kerr spacetime. If one makes the (implicitly defined) change
of variables $\di\tau/\di\lambda=\KSigma^{-1}$, then equation
\eqref{eq:IntroNullGeodesicsr} for the radial component becomes $(\di
r/\di\tau)^2=-\CurlyR(r;M,a;\GeodesicEnergy,\GeodesicLz,\GeodesicQ)$. For
fixed $(M,a)$ and $(\GeodesicEnergy,\GeodesicLz,\GeodesicQ)$, this
takes the form of the equations of motion of a particle in $1$-dimension with a potential. The roots and double roots of the potential $\CurlyR$ determine the turning points and stationary points, respectively, for the motion in the $r$ direction. The potential $-\CurlyR=((r^2+a^2)\GeodesicEnergy+a\GeodesicLz)^2 -\KDelta(\GeodesicQ+\GeodesicLz^2+2a\GeodesicEnergy\GeodesicLz)$ is always nonnegative at $r=\rp$ and, unless $\GeodesicEnergy=0$, is positive as $r\rightarrow\infty$. As we will show below, it has at most two roots counting multiplicity. 
%as $r\rightarrow\infty$ and has no roots, a double root, or exactly two (distinct, single) roots \mnote{Citation needed. Check Frolov-Novikov}. 
%% We'll only include true statements. 

By simply considering the turning points, one can use $r$ and $\dot\gamma_r$ to construct a quantity that is increasing on average. Throughout this argument, we will take $(M,a)$ and $(\GeodesicEnergy,\GeodesicLz,\GeodesicQ)$ as fixed. Consider null geodesics which came from infinity, i.e.{} for which $r\rightarrow\infty$ as $\tau\rightarrow-\infty$. For these, $\GeodesicEnergy$ is positive and $-\CurlyR$ has no roots, a single double root, or two simple roots. If the potential $-\CurlyR$ has no roots, then we can arbitrarily choose $\rorbit\in(\rp,\infty)$, so that when the geodesic falls in from infinity, the quantity $(r-\rorbit)\dot\gamma_r$ goes from negative to positive. Similarly, if there are two distinct roots, we can choose $\rorbit$ between these two roots (which are between $\rp$ and $\infty)$, in which case, before the geodesic reaches the turning point, the quantity $(r-\rorbit)\dot\gamma_r$ is negative, but that after the geodesic leaves the turning point, the quantity $(r-\rorbit)\dot\gamma_r$ is positive. Finally, in the case that there is a double root, we can define the root to be $\rorbit$, so that $(r-\rorbit)\dot\gamma_r$ is be large and negative in the far past, but that it goes to zero as the null geodesic asymptotes onto the double root of $-\CurlyR$. Using the terminology from the start of this section and taking $\vecMGeodesic=(r-\rorbit)\dr$, we write $(r-\rorbit)\dot\gamma_r$ as the energy $-\GenEnergyGeodesic{\vecMGeodesic}$. In all three cases considered, $\GenEnergyGeodesic{\vecMGeodesic}$ is decreasing overall, in the sense that the limit in the future that is less than the limit in the past. A similar analysis shows that $\GenEnergyGeodesic{\vecMGeodesic}$ is nonincreasing overall along all other null geodesics. (Along the orbiting null geodesics, it is identically zero for all $\tau$.) Thus, in all cases, we can define $\rorbit$ and $\vecMGeodesic=(r-\rorbit)\dr$ so that $\GenEnergyGeodesic{\vecMGeodesic}$ is nonincreasing overall. 

To construct a monotone quantity on each null geodesic, it remains to replace $r-\rorbit$ by a function $\fnMrGeodesic$ so that for $\vecMGeodesic=\fnMrGeodesic\dr$, the energy $\GenEnergyGeodesic{\vecMGeodesic}$ is nondecreasing for all $\tau$ and not merely nondecreasing overall. For $a\not=0$, both $\rorbit$ and $\fnMrGeodesic$ will have to depend on both the Kerr parameters $(M,a)$ and the constants of motion $(\GeodesicEnergy,\GeodesicLz,\GeodesicQ)$; the function $\fnMrGeodesic$ will also depend on $r$, but no other variables. We define $\vecMGeodesic=\fnMrGeodesic\dr$ and emphasise to the reader that this is a map from the tangent bundle to the tangent bundle, which is not the same as a standard vector field, which is a map from the manifold to the tangent bundle. To derive a monotonicity formula, we wish to choose $\fnMrGeodesic$ so that $\GenEnergyGeodesic{\vecMGeodesic}$ has a nonnegative derivative. We define the covariant derivative of $\vecMGeodesic$ by holding the values of $(\GeodesicEnergy,\GeodesicLz,\GeodesicQ)$ fixed and computing the covariant derivative as if $\vecMGeodesic$ were a regular vector field. Similarly, we define $\Lie_{\vecMGeodesic}\gMetric^{\alpha\beta}$ by fixing the values of the constants of geodesic motion. Since the constants of motion have zero derivative along null geodesics, equation \eqref{eq:deformForGeodesics} remains valid. 

We can use this to illustrate the core calculation of this paper. The Kerr metric can be written as 
\begin{align*}
\KSigma\gMetric^{\alpha\beta}={}&{} \KDelta\dr^\alpha\dr^\beta +\frac{1}{\KDelta}\CurlyR^{\alpha\beta} ,
\end{align*}
where
\begin{align}
\CurlyR^{\alpha\beta}={}&{} -(r^2+a^2)^2\dt^\alpha\dt^\beta -4aMr
\dt^{(\alpha}\dphi^{\beta)} +(\KDelta-a^2)\dphi^\alpha\dphi^\beta
+\KDelta\TensorQ^{\alpha\beta}, 
\label{eq:DefCurlyRalphabeta}\\
\TensorQ^{\alpha\beta}={}&{}\dtheta^\alpha\dtheta^\beta+\cot^2\theta\dphi^\alpha\dphi^\beta+a^2\sin^2\theta\dt^\alpha\dt^\beta .
\end{align}
The double contraction of the tensor $\CurlyR^{\alpha\beta}$ with the
tangent to a null geodesic gives the potential
$\CurlyR(r;M,a;\GeodesicEnergy,\GeodesicLz,\GeodesicQ)=\CurlyR^{\alpha\beta}\dot\gamma_\alpha\dot\gamma_\beta$. The
crucial quantity
$\Lie_{\vecMGeodesic}\gMetric^{\alpha\beta}\dot\gamma_\alpha\dot\gamma_\beta$
is calculated below in \eqref{eq:GeodesicBulk}. For now, we ignore
distracting factors of $\KSigma$, $\KDelta$, their derivatives, and
constant factors, so we can see that the most important terms in $\Lie_{\vecMGeodesic}\gMetric^{\alpha\beta}\dot\gamma_\alpha\dot\gamma_\beta$ are 
\begin{align*}
-2(\dr\fnMrGeodesic)\dot\gamma_r\dot\gamma_r +\fnMrGeodesic(\dr\CurlyR^{\alpha\beta})\dot\gamma_\alpha\dot\gamma_\beta =-2(\dr\fnMrGeodesic)\dot\gamma_r\dot\gamma_r +\fnMrGeodesic(\dr\CurlyR).
\end{align*} 
The second term in this sum will be nonnegative if
$\fnMrGeodesic=\dr\CurlyR(r;M,a;\GeodesicEnergy,\GeodesicLz,\GeodesicQ)$. Recall
that the vanishing of
$\dr\CurlyR(r;M,a;\GeodesicEnergy,\GeodesicLz,\GeodesicQ)$ is one of
the two conditions \eqref{eq:RotatingGeodesicCondition} for orbiting null geodesics. With this choice of $\fnMrGeodesic$, the instability of the null geodesic orbits ensures that, for these null geodesics, the coefficient in the first term, $-2(\dr\fnMrGeodesic)$, will be positive. We can now perform the calculations more carefully to show that this nonnegativity holds for all null geodesics. 

Since, up to reparametrisation, null geodesics are conformally invariant, it is sufficient to work with the conformally rescaled metric $\KSigma\gMetric^{\alpha\beta}$. Furthermore, since $\gamma$ is a null geodesic, for any function $\fnMpGeodesic$, we may add $\fnMpGeodesic\KSigma\gMetric^{\alpha\beta}\dot\gamma_\alpha\dot\gamma_\beta$ wherever it is convenient. Thus, the change in $\GenEnergyGeodesic{\vecMGeodesic}$ is given as the integral of
\begin{align*}
\KSigma\dot\gamma_\alpha\dot\gamma_\beta \nabla^{(\alpha}\vecMGeodesic^{\beta)}
{}&{}=\left(-\frac12\Lie_{\vecMGeodesic}(\KSigma\gMetric^{\alpha\beta}) +\fnMpGeodesic\KSigma\gMetric^{\alpha\beta}\right)\dot\gamma_\alpha\dot\gamma_\beta 
%\\
%{}&{}=\left(-\frac12\fnMrGeodesic\left(\dr\frac{\CurlyR^{\alpha\beta}}{\KDelta}\right)+\fnMpGeodesic\frac{\CurlyR^{\alpha\beta}}{\KDelta}\right)\dot\gamma_\alpha\dot\gamma_\beta \\
%{}&{}\qquad
%+\left(-\frac12\fnMrGeodesic(\dr\KDelta)+\KDelta(\dr\fnMrGeodesic)-\fnMpGeodesic\right)\dot\gamma_r\dot\gamma_r .
\end{align*}

To progress further, choices of $\fnMrGeodesic$ and $\fnMpGeodesic$
must be made. For the choices we make here, the calculations are
straightforward but lengthy. More detail is given in subsections
\ref{SS:SetUpForRadial}-\ref{SS:ChoosingTheWeights} where analogous
calculations are made for the wave equation. Let $\fnMna$ and $\fnMnb$
be smooth functions of $r$ and the Kerr parameters $(M,a)$. We
introduce the notation
\begin{align*}
\DiffCurlyRTilde{}={}&\dr\left(\frac{\fnMna}{\KDelta}\CurlyR(r;M,a;\GeodesicEnergy,\GeodesicLz,\GeodesicQ)\right),&
\DDiffCurlyRTTilde{}={}&
\left(\dr\left(\fnMnb\frac{\fnMna^{1/2}}{\KDelta^{1/2}}\DiffCurlyRTilde\right)\right), 
\end{align*}
and make the choices
\begin{align*}
\fnMrGeodesic={}&{}-\fnMna\fnMnb\DiffCurlyRTilde,&
\fnMpGeodesic={}&{}-(1/2)(\dr\fnMna)\fnMnb\DiffCurlyRTilde .
\end{align*}
In terms of these functions,
\begin{align}
\dot\gamma_\alpha\dot\gamma_\beta \nabla^{(\alpha}\vecMGeodesic^{\beta)}
{}&{}=\fnMnb (\DiffCurlyRTilde)^2 -\fnMna^{1/2}\KDelta^{3/2}\DDiffCurlyRTTilde \dot\gamma_r^2. 
\label{eq:GeodesicBulk}
\end{align}
If we now take $\fnMna=\fnMca=\KDelta(r^2+a^2)^{-2}$, then the
coefficient of $\GeodesicEnergy^2$ in $\DiffCurlyRTilde$ vanishes, and if we further take
$\fnMnb=\fnMcb=(r^2+a^2)^4/(3r^2-a^2)$, then the coefficient of
$\GeodesicEnergy\GeodesicLz$ in $\DDiffCurlyRTTilde$ also vanishes, leaving
\begin{subequations}
\begin{align}
\DiffCurlyRTilde
{}&{}= 4Ma\frac{3r^2-a^2}{(r^2+a^2)^3}\GeodesicEnergy\GeodesicLz\nonumber\\
{}&{}\quad-2\frac{r^3-3Mr-a^2r+Ma^2}{(r^2+a^2)^3}\GeodesicLz^2
-2\frac{r^3-3Mr+a^2r+Ma^2}{(r^2+a^2)^3}\GeodesicQ ,
\label{eq:GeodesicsDiffCurlyRTilde}\\
\DDiffCurlyRTTilde
{}&{}=-2\frac{3r^4+a^4}{(3r^2-a^2)^2}\GeodesicLz^2
-2\frac{(3r^2-6a^2r^2+a^4)^2}{(3r^2-a^2)^2}\GeodesicQ .
\label{eq:GeodesicsDDiffCurlyRTTilde}
\end{align}
\end{subequations}
Since $\GeodesicQ$ is nonnegative by equation \eqref{eq:IntroDefGeodesicQ}, it follows that the right-hand side of \eqref{eq:GeodesicsDDiffCurlyRTTilde} is nonpositive and that the right-hand side of equation \eqref{eq:GeodesicBulk} is nonnegative. Since equation \eqref{eq:GeodesicBulk} gives minus the rate of change, the energy $\GenEnergyGeodesic{\vecMGeodesic}$ is monotone. 

Furthermore, these calculations reveal useful information about the geodesic motion. The positivity of the term on the right-hand side of \eqref{eq:GeodesicsDDiffCurlyRTTilde} shows that $\DiffCurlyRTilde$ can have at most one root, which must be simple. In turn, this shows that $\CurlyR$ can have at most two roots, as previously asserted. 

The role of orbiting geodesics can be seen in equation
\eqref{eq:GeodesicBulk}. Along null geodesics for which $\CurlyR$ has
a double root, the double root occurs at the root of
$\DiffCurlyRTilde$, so it is convenient to think of the corresponding
value of $r$ as being $\rorbit$. In particular, this root is where
null geodesics with the given values of $\GeodesicEnergy$,
$\GeodesicLz$, and $\GeodesicQ$ orbit the black hole with a constant
value of $r$. The first term in \eqref{eq:GeodesicBulk} vanishes at
the root of $\DiffCurlyRTilde$, as it must so that
$\GenEnergyGeodesic{\vecMGeodesic}$ can be constantly zero on the
orbiting null geodesics. When $a=0$, the quantity $\DiffCurlyRTilde$
reduces to $-2(r-3M)r^{-4}(\GeodesicLz^2+\GeodesicQ)$, so that the
orbits occur at $r=3M$. The continuity in $a$ of $\DiffCurlyRTilde$
guarantees that its root converges to $3M$ as $a\rightarrow 0$ for
fixed $(\GeodesicEnergy,\GeodesicLz,\GeodesicQ)$. In subsection
\ref{SS:ChoosingTheWeights}, a slightly more complicated choice of
$\fnMna$ leads to an $\DiffCurlyRTilde$ for which the convergence of the root to $3M$ as $a\rightarrow 0$ can be made uniform in $(\GeodesicEnergy,\GeodesicLz,\GeodesicQ)$. 

Because of the geometric optics correspondence between null geodesics and solutions of the wave equation, it is natural to try to adapt the monotonicity of $\GenEnergyGeodesic{\vecMGeodesic}$ for null geodesics to a similar result for the wave equation and, in particular, to adapt the nonnegativity of the terms in equation \eqref{eq:GeodesicBulk} to help in the proof of the Morawetz estimate. In making this adaptation, one must find a replacement for the constants of motion as arguments in the weight $\fnMrGeodesic$. There are several ways in which this can be done. One approach \cite{DafermosRodnianski:KerrEnergyBound} is to use the complete separability of the wave equation; to observe that separation of variables in the $t$ and $\phi$ coordinates is equivalent to the Fourier transform; to observe that the Fourier variables conjugate to $t$ and $\phi$ can be treated as $\GeodesicEnergy$ and $\GeodesicLz$; and to treat the final separation constant, typically associated with $\theta$ but perhaps more properly thought of as the eigenvalues of the hidden symmetry $\OpQ$, as analogous to $\GeodesicQ$; to construct a monotone energy like $\GenEnergyGeodesic{\vecMGeodesic}$ at least for some values of $(\GeodesicEnergy,\GeodesicLz,\GeodesicQ)$; and then to show that the estimate on separated components can be summed to provide an estimate for general solutions. Another approach \cite{TataruTohaneanu} to treating the wave equation is to define a pseudo-differential operator with symbol given by $\vecMGeodesic$; this is possible since $\vecMGeodesic$ is a map from the tangent bundle to the tangent bundle. 

The method that we introduce in this paper uses only differential operators. Since $\vecMGeodesic$ is constructed only from the constants $\{\GeodesicEnergy^2,\GeodesicEnergy\GeodesicLz,\GeodesicLz^2,\GeodesicQ\}$, which are quadratic in $\dot\gamma$ and constructed from the conformal Killing tensors $\{\dt^\alpha\dt^\beta, \dt^{(\alpha}\dphi^{\beta)}, \dphi^\alpha\dphi^\beta,\TensorQ^{\alpha\beta}\}$, our approach is to construct an analogue which primarily uses the symmetries $\SOSym=\{\dt^2,\dt\dphi,\dphi^2,\OpQ\}$, which are second-order differential operators constructed from the same set of Killing tensors. 
%Since we want to construct a quadratic form associated with $\vecMGeodesic$, we multiply $\vecMGeodesic$ by $\GeodesicL$ so that a symmetry from $\SOSym$ can be applied to each term in the quadratic form. \mnote{Is the additional information about $\GeodesicL$ distracting?}
In \cite{DafermosRodnianski:KerrEnergyBound}, the Morawetz estimate is only proved for ``high'' frequency waves, which have a large ratio between certain frequencies corresponding to the constants of motion for null geodesics $(\GeodesicEnergy,\GeodesicLz,\GeodesicQ)$; such a decomposition is not necessary in deriving the nonnegativity for null geodesics in equation \eqref{eq:GeodesicBulk} or when proving estimates for the wave equation in \cite{TataruTohaneanu} or this paper. In summary, our approach allows us to use differential operators to construct a multiplier which treats all frequency ranges in a uniform manner and in particular gives a nonnegative bulk term at $r \sim 3M$.

\newpageForSubsection
%%%%%%%%%%%%%%%%%%%%%%%%%%%%%%%%%%%%%%%%%%%%%%%%%%%%%%%%%%%%%%%%%%%%%%%%%%%%%
\subsection{Generalising the vector-field method}
%{Hidden symmetries and the vector-field method}
\label{SS:genvect} 
In this section, we outline a generalisation of the
vector-field method which allows us to take advantage of the presence
of hidden symmetries in the Kerr spacetime. In particular, we consider
energies based on operators of order greater than one, rather than
just vector fields.
% 
%Let $\gWave = \nabla^\alpha \nabla_\alpha$.  \mnote{LA: already def} 
In the discussion here,
we consider the scalar wave equation $\gWave \solu = 0$, but much of
the discussion applies equally to general field equations derived from
a quadratic action.

The energy-momentum tensor for the wave equation is
\begin{equation} \label{eq:em-tensor} 
\StressEnergyGeneral[\solu]_{\alpha\beta}= \nabla_\alpha\solu \nabla_\beta\solu -\frac12\gMetricGeneral_{\alpha\beta}\left(\nabla_\gamma\solu\nabla^\gamma\solu \right). 
\end{equation} 
The momentum associated with a vector field $\vecX$ and the energy
associated with a vector field $\vecX$ and evaluated on a hypersurface
$\HypersurfaceGeneral$ are
\begin{subequations}
\begin{align}
\GenMomentum{\vecX}[\solu]_\alpha ={}&{} \StressEnergyGeneral[\solu]_{\alpha\beta}\vecX^\beta , 
\label{eq:momentum}\\
\GenEnergy{\vecX}[\solu](\HypersurfaceGeneral)={}&{} -\int_\HypersurfaceGeneral \GenMomentum{\vecX}[\solu]_\alpha \diNormal^\alpha ,
\label{eq:energy}
\end{align}
\end{subequations}
where $\diNormal^\alpha$ is the normal volume form on
$\HypersurfaceGeneral$. That is, for any $1$-form $\formxi$,
$\int_\Sigma \formxi_\alpha\diNormal^\alpha$ is the integral over
$\HypersurfaceGeneral$ of the $3$-form given by the Hodge dual $*
\formxi$ \cite{wald:MR757180}. When the spacetime is foliated by
surfaces of constant time, as is the case in the Kerr spacetime, we
will denote these surfaces by $\hst{t}$. In this case, we take the
normal to be past directed. Thus, the sign in the energy ensures that
the energy is nonnegative on constant $t$ surfaces if $\vecX$ is
future directed and timelike.  In the following, unless there is room
for confusion, we will drop reference to $\solu$ in the notation for
momentum and energy. 

The energy momentum tensor (\ref{eq:em-tensor}) satisfies the dominant
energy condition, and hence for $\vecX$ future-oriented and timelike,
the energy induced on a hypersurface with a past-oriented timelike
normal (i.e.\ a spacelike hypersurface) is positive definite. The
energy conservation law takes the form
\begin{align}
\GenEnergy{\vecX}(\HypersurfaceGeneral_2)-\GenEnergy{\vecX}(\HypersurfaceGeneral_1)
={}&{} -\int_\Omega \left(\nabla_\alpha\GenMomentum{\vecX}^\alpha  \right) \diFourNatural \nonumber\\
={}&{} -\int_\Omega  \StressEnergyGeneral_{\alpha\beta}\nabla^{(\alpha}\vecX^{\beta)}   \diFourNatural , 
\label{eq:deform}
\end{align}
where $\Omega$ is the region enclosed between $\HypersurfaceGeneral_1$
and $\HypersurfaceGeneral_2$. This is often referred to as the
deformation formula. The sign in the right-hand side arises from Stokes'
theorem and our sign conventions. Energy estimates are often performed
by controlling the bulk (also called deformation) terms $\nabla_\alpha
\GenMomentum{\vecX}^\alpha$. However, for Morawetz estimates (e.g.{}
inequality \eqref{eq:Morawetz} below), one makes use of the sign of the bulk
term itself; this is similar to the derivation of the monotonicity
formula for null geodesics in subsection \ref{SS:ResultForGeodesics}.

Recall that the wave equation is a Lagrangian field equation. Taking this into account, equation \eqref{eq:deform} with a momentum 1-form as in \eqref{eq:momentum} is simply a restatement of Noether's theorem. We will now consider generalisations of the deformation formula, involving momentum 1-forms, and energies, which are not derived directly from a deformation of a Lagrangian. These generalisations include the addition of lower-order correction terms, which is a familiar feature of the multiplier method, as well as the introduction of higher-order conservation laws defined in terms of symmetry operators of the field equation. The existence of symmetry operators is closely related to separability properties for field equations and has been studied for a long time, see e.g.{} \cite{Koornwinder,miller:MR0460751, MR1240056}
and references therein. However, the application of these ideas in the context of the vector-field method and, in particular, with nonsymmetries as in our proof of the Morawetz estimate is new.

By estimating higher-order energies one may, via Sobolev estimates,
obtain pointwise control of the fields. Higher-order energies may be
defined by using symmetries. Given, for $0\leq i\leq n$, a
collection of order-$i$ differential operators, $\SymGeneraln{i}$,
we can define the higher-order energy (of order $n+1$) for a vector
field $\vecX$ to be
\begin{align}
\GenEnergyOrder{\vecX}{n+1}[\solu](\HypersurfaceGeneral)
=\sum_{i=0}^n \sum_{\CQA\in\SymGeneraln{i}} \GenEnergy{\vecX}[\CQA \solu](\HypersurfaceGeneral) .
\label{eq:DefHigherEnergy}
\end{align}

Since the energy-momentum tensor is quadratic in $\solu$, we can
define a bilinear form of it 
by
\begin{align*}
\StressEnergyGeneral[\solu_1,\solu_2]_{\alpha\beta}
={}&{}\frac14\left(\StressEnergyGeneral[\solu_1+\solu_2]_{\alpha\beta} -\StressEnergyGeneral[\solu_1-\solu_2]_{\alpha\beta} \right) . 
\end{align*}
It is convenient to define an indexed version of the bilinear 
energy momentum, with respect to a set of symmetry operators $\{\CQA_\ua\}$, by  
\begin{align*}
\StressEnergyGeneral[\solu]_{\ua\ub\alpha\beta}
={}&{}\StressEnergyGeneral[\CQA_\ua\solu,\CQA_\ub\solu]_{\alpha\beta} .
\end{align*}
Given a double-indexed collection of vector fields, $\STwovecX=\{\vecX^{\ua\ub}\}$, we define the associated generalised momentum and energy to be 
\begin{align*}
\GenMomentum{\STwovecX}[\solu]_\alpha 
={}&{}\StressEnergyGeneral[\solu]_{\ua\ub\alpha\beta} \vecX^{\ua\ub\beta} ,\\
\GenEnergy{\STwovecX}[\solu](\HypersurfaceGeneral)
={}&{} -\int_\HypersurfaceGeneral \GenMomentum{\STwovecX}[\solu]_\alpha \diNormal^\alpha .
\end{align*}

In practice it is convenient to consider momenta with lower-order
terms, designed to improve certain deformation terms in
$\nabla_\alpha\GenMomentum{\vecX}^\alpha$.  For a scalar function,
$\MMTTscalar$ (\cite{MarzuolaMetcalfeTataruTohaneanu}, but previously appearing in
\cite{DafermosRodnianski:MorawetzWithoutHarmonics}), or a double-indexed collection of
functions, $\STwoMMTTscalar=\{\MMTTscalar^{\ua\ub}\}$, the associated momenta are defined
to be 
\begin{align*}
\GenMomentum{\MMTTscalar}[\solu]_\alpha
={}&{} \MMTTscalar(\nabla_\alpha\solu)\solu -\frac12(\partial_\alpha\MMTTscalar)\solu^2 , \\
\GenMomentum{\STwoMMTTscalar}[\solu]_\alpha
={}&{} \MMTTscalar^{\ua\ub}(\nabla_\alpha\CQA_\ua\solu)\CQA_\ub\solu -\frac12(\partial_\alpha\MMTTscalar^{\ua\ub})(\CQA_\ua\solu)(\CQA_\ub\solu) . 
\end{align*}
For a pair consisting of a vector field and a scalar function,
$(\vecX,\MMTTscalar)$, the associated momentum is defined to be the
sum of the momenta associated with the vector field and the scalar.
For a pair of collections,
$\STwovecX=(\{\vecX^{\ua\ub}\},\{\MMTTscalar^{\ua\ub}\})$, again the
momentum is defined to be the sum of the momenta.  In all cases, the
energy on a hypersurface is given by the negative of the flux, defined
with respect to the momentum vector, through the hypersurface.  From
the above, we have the following version of the deformation formula,
\begin{align}
\GenEnergy{\STwovecX}(\HypersurfaceGeneral_2)-\GenEnergy{\STwovecX}(\HypersurfaceGeneral_1)
={}&{} -\int_\Omega \left(\nabla_\alpha\GenMomentum{\STwovecX}^\alpha  \right) \diFourNatural .
\label{eq:deform:twovec}
\end{align}

It is important to point out, as we show in lemma \ref{Lemma:DivOfP},
that the deformation terms for the generalised momenta are
computationally not much more difficult to handle than the classical
ones.  As for the classical momenta and energies, in defining the
generalised vector fields, momenta, and energies as outlined above,
one is interested in getting positive definiteness of the energies or
bulk terms. Here, an additional subtlety arises. Namely, in the
Morawetz estimate presented in equation \eqref{eq:Morawetz}, one
achieves positive definiteness only modulo boundary terms. We generate
these boundary terms when we integrate by parts to use the formal
self-adjointness of the second-order symmetry operators. These
boundary terms can then be controlled by the energy. The presence of
these boundary terms is a completely new feature compared to the
classical energies and momenta.

\subsection{Strategy of proof}
\label{SS:Strategy}
Recall from earlier in the introduction that there are three major
problems in the Kerr spacetime:
\begin{enumerate}
\item{No positive, conserved energy:} There is no timelike, Killing
  vector. In particular, the vector field $\dt$, which is Killing, is only
  timelike outside the ergosphere, i.e.\ for
  $r>M+\sqrt{M^2-a^2\cos^2\theta}$.
\label{Problem:NoPositiveConservedEnergy} 
\item{Lack of sufficient classical symmetries:} The higher energies generated
  by the Lie derivatives in the $\dt$ and $\dphi$ directions do not control
  enough directions to estimate Sobolev norms of the function. 
\label{Problem:InsufficientClassicalSymmetries}
\item{Complicated trapping:} There are orbiting null geodesics. These
  neither escape to null infinity nor enter the black hole. Since
  solutions to the wave equation can follow null geodesics for an
  arbitrarily long time, this presents an obstacle to
  decay. Furthermore, there are orbiting null geodesics occurring over
  a range of $r$ in the Kerr spacetime (with $|a|>0$), which makes the
  situation more complicated than in the Schwarzschild spacetime
  ($a=0$), where there are orbiting null geodesics only at $r=3M$.
\label{Problem:ComplicatedTrapping}
\end{enumerate}
%We shall compare our approach to overcoming these problems to those in
%\cite{DafermosRodnianski:KerrEnergyBound,TataruTohaneanu} at the end
%of this subsection.

To overcome problem \ref{Problem:NoPositiveConservedEnergy}, we first observe that the vector $\dt$ is timelike for sufficiently large $r$; that, if 
\begin{align*}
\omegaH=\frac{a}{\rp^2+a^2}
\end{align*} 
denotes the angular velocity of the horizon, then the vector $\dt+\omegaH\dphi$ is null on the horizon and timelike for sufficiently small $r>\rp$; that the regions where $\dt$ and $\dt+\omegaH\dphi$ are timelike overlap when $|a|$ is sufficiently small; and that both $\dt$ and $\dt+\omegaH\dphi$ are Killing. Thus, if we let
\begin{align}
\vecTBlend = \dt +\fnBlend \omegaH \dphi ,
\label{eq:IntroDefTBlend} 
\end{align}
where $\fnBlend$ is identically $1$ for $r<\epsilonInvBlendLocation$ for some constant $\epsilonInvBlendLocation$, identically $0$ for $r>\epsilonInvBlendLocation +M$, and smoothly decreases on $[\epsilonInvBlendLocation,\epsilonInvBlendLocation+M]$, then, for sufficiently small $a$, this vector field will be timelike everywhere and will be Killing outside the fixed region $r\in[\epsilonInvBlendLocation,\epsilonInvBlendLocation+M]$. Thus, to prove the boundedness of the associated positive energy, it will be sufficient to control the behaviour of solutions in this fixed region through a Morawetz estimate. 

To overcome problem \ref{Problem:InsufficientClassicalSymmetries}, we note that the second-order operator $\OpQ$ is a symmetry and is a weakly elliptic 
operator. Using $\OpQ$, $\dphi^2$, and $\dt^2$ as symmetries to generate higher energies, we can control energies of the spherical Laplacian of $\solu$. These control Sobolev norms which are sufficiently strong to control $|\solu|^2$. 

To overcome problem \ref{Problem:ComplicatedTrapping}, the complicated trapping, we will adapt $\vecMGeodesic$ from
the monotonicity formula for null geodesics in subsection
\ref{SS:ResultForGeodesics}. This adaptation is possible because the
double-indexed energy momentum tensor
$\StressEnergyGeneral[\solu]_{\ua\ub\alpha\beta}$ from subsection
\eqref{SS:genvect} allows one to use the hidden symmetries in defining double-indexed sets of vectors. If we introduce 
\begin{align*}
\OpL{}&{}=\OpL^\ua\CQA_\ua =\dt^2+\dphi^2+\OpQ
\end{align*} 
to give us a weakly elliptic operator and an extra, free, 
underlined index, we can take as our collection of Morawetz vector fields
\begin{align*}
\vecM ={}&{} -\fnMna\fnMnb \DiffCurlyRTilde^{(\ua}\OpL^{\ub)} \dr ,\\
\fnM={}&{} -\frac12\fnMna \left(\dr\left( \fnMnb\DiffCurlyRTilde^{(\ua}\right)\right)\OpL^{\ub)} , \\
\pairM={}&{}(\{\vecM\},\{\fnM\}),\\
\DiffCurlyRTilde^\ua={}&{}\dr\left(\frac{\fnMna}{\KDelta}\CurlyR^{\ua}\right) , 
\end{align*}
with $\fnMna$ and $\fnMnb$ smooth, positive functions to be chosen. In
subsection \ref{SS:ChoosingTheWeights}, we choose $\fnMna$ and
$\fnMnb$ slightly differently from how they were chosen in subsection
\ref{SS:ResultForGeodesics}, so that they satisfy some additional
conditions that were not necessary there. Applying the
generalised deformation formula \eqref{eq:deform:twovec}, the
difference between the energies on one hypersurface and another is
\begin{align*}
\GenEnergy{\pairM}(\HypersurfaceGeneral_2 )-\GenEnergy{\pairM}(\HypersurfaceGeneral_1 )
={}&{} - \int \left( \nabla_\alpha\GenMomentum{\pairM}^\alpha \right)  \diFourNatural . 
\end{align*}
Ignoring several distracting details, the bulk term 
%deformation 
is of the form
\begin{align}
{}&{}\frac12 \fnMnb \DiffCurlyRTilde^\ua\DiffCurlyRTilde^\ub\OpL^{\alpha\beta} (\partial_\alpha\CQA_\ua\solu)(\partial_\beta\CQA_\ub\solu) \nonumber \\
{}&{}+\fnMna^{1/2}\KDelta^{3/2}\left(-\dr\left(\fnMnb\frac{\fnMna^{1/2}}{\KDelta^{1/2}}\DiffCurlyRTilde^{\ua}\right)\right)\OpL^\ub (\dr\CQA_\ua\solu)(\dr\CQA_\ub\solu)  \nonumber \\
{}&{}+\frac14(\dr\KDelta\dr\fnMna(\dr\fnMnb\DiffCurlyRTilde^{\ua}))\OpL^\ub (\CQA_\ua\solu)(\CQA_\ub\solu) , \label{eq:BigImportantNotNumberedMorawetz}
\end{align}
see the proof of lemma \ref{Lemma:Morawetz}. The first two terms are very similar to the terms in equation \eqref{eq:GeodesicBulk} for null geodesics. 
In the first line, one factor of $\DiffCurlyRTilde$ arises from the wave equation, and the other from our choice of the Morawetz vector field $\pairM$, which allows us to construct a perfect square to obtain positivity. The term in the second line involves two derivatives of $-\CurlyRTilde$. Near the photon orbits, 
the convexity properties of $\CurlyR$, which ensured that the orbits are unstable, ensure that this term is positive. We choose $\fnMna$ and $\fnMnb$ to get positivity away from the photon orbits. The third term is lower-order, since it involves fewer derivatives of $\solu$, and has no analogue for null geodesics.

Recall from the discussion above that when viewed as a function on
phase space, $\CurlyR$ vanishes together with $\dr \CurlyR$ at the
trapped set, and $\dr^2 \CurlyR < 0$ there.  Since $\fnMna$ and
$\fnMnb$ are positive, the equivalent statement holds for
$\CurlyRTilde$. The expression in
\eqref{eq:BigImportantNotNumberedMorawetz} is quadratic in up to
third-order derivatives of the field $\solu$.  For this reason it is
not appropriate in the context of
\eqref{eq:BigImportantNotNumberedMorawetz} to think of $\CurlyR$ as a
function on phase space, but rather to compare this expression with
the square of a nondegenerate pointwise third-order norm of
$\solu$. The vanishing of $\CurlyR$ and $\dr \CurlyR$ at the orbiting
null geodesics then is reflected in a degeneracy of the first line of
\eqref{eq:BigImportantNotNumberedMorawetz} compared to such a
third-order norm. This analysis is done in detail in section
\ref{SS:ChoosingTheWeights}, see in particular the proof of lemma
\ref{Lemma:Morawetz}.

For small $|a|$, with $v$ denoting terms of the form $\CQA_\ua\solu$, and with our choices of $\fnMna$ and $\fnMnb$, the sum of the second and third terms in \eqref{eq:BigImportantNotNumberedMorawetz}
is of the form
\begin{align}
M\left( \frac{\KDelta^2}{r^2(r^2+a^2)} (\dr v)^2 +\frac{9r^2-46Mr +54M^2}{6 r^4} v^2 \right)
\label{eq:IntroNeedHardy}
\end{align}
with small perturbations on the coefficients. The coefficient on $v^2$ is positive outside a compact interval in $(\rp,\infty)$. As shown in \cite{BlueSoffer:ODE}, it is sufficient to prove a Hardy estimate which bounds the quadratic form in \eqref{eq:IntroNeedHardy} from below by a sum of positive weights times $(\dr v)^2$ and $v^2$. 
This provides an estimate on the third line in \eqref{eq:BigImportantNotNumberedMorawetz}, which has no analogue for null geodesics. 

The positive terms arising from the deformation of $\pairM$ dominate
the deformation terms (with extra derivatives) arising in the failure
of $\vecTBlend$ to be Killing. (In fact, at this stage, only the
second and third derivatives of $\solu$ are controlled, whereas the
deformation terms from the third-order $\vecTBlend$ energy also
involve the first derivatives and undifferentiated factors of
$\solu$. To handle this, a classical vector field is also introduced
to prove a Morawetz estimate that controls the lower-order terms.) On
the other hand, the energy associated with $\pairM$ is dominated by
the (third-order) energy associated with $\vecTBlend$. Since there is
a factor of $|a|$ on the $\vecTBlend$ deformation terms, we have a
small parameter, which allows us to close the boot-strap argument in
which the $\vecTBlend$ energy is controlled by the integral of its
deformation, which is controlled by the integral of the $\pairM$
deformation, which is controlled by the $\pairM$ energy, which is
finally controlled by the $\vecTBlend$ energy. This allows us to
establish theorem \ref{Thm:IntroUniformEnergyBound}.

A similar argument can be used to show that for null geodesics, there is a uniform bound on the positive energy $\GenEnergyGeodesic{\GeodesicL^2\vecTBlend}[\Geodesic]$, where $\GeodesicL=\GeodesicEnergy^2+\GeodesicLz^2+\GeodesicQ$. The method is essentially the same as for the wave equation. With $\vecMGeodesic=\GeodesicL\fnMrGeodesic\dr$, with $\fnMna$ and $\fnMnb$ chosen as in subsection \ref{SS:ChoosingTheWeights}, the positive terms in the deformation dominate the deformation terms arising from the failure of $\vecTBlend$ to be Killing. Similarly, the $\vecTBlend$ energy (with two extra factors of $\GeodesicL$) dominates the $\vecMGeodesic$ energy. 

The small $|a|$ condition which we impose is significantly stronger
than the condition that $|a|\leq M$ which implies the existence of a
black hole and which might be ideally imposed. There are several
fundamental and technical reasons for this small $|a|$
condition. Perhaps most importantly, the construction of $\vecTBlend$
relies on there being a region where both $\dt+\fnBlend\omegaH\dphi$
and $\dt$ are timelike in which to perform the blending. When $|a|$ is
sufficiently large, but still smaller than $M$, there is no such
overlapping region, so this particular construction fails. In
addition, we use the assumption on the smallness of $|a|$ to close the
bounded $\vecTBlend$ energy argument.  If $|a|$ is not small relative
to the absolute constants appearing in that estimate, it would not be
possible to close the boot-strap argument. A clear technical obstacle
is that, in the proof of the Morawetz estimate, we perturb the Hardy
estimate in \eqref{eq:IntroNeedHardy}. If $|a|$ were too large, the
perturbation argument would fail, and our numerical investigation
suggests that when $|a|$ is larger than about $.9 M$, there are no
longer positive solutions of the associated ODE, which we use to prove
the estimate. These obstacles are the most fundamental obstacles to
extending the range of $|a|$, but there are also numerous other,
technical estimates in which we have made use of the smallness of
$|a|$.

Having summarised our method, we will now compare it with methods used
in recent, related work. Recently, others have constructed a bounded
energy \cite{DafermosRodnianski:KerrEnergyBound,TataruTohaneanu}. To make a comparison, we point to
several features which they share but which are different from those
in our approach. 

To overcome problem \ref{Problem:NoPositiveConservedEnergy}, we use
$\vecTBlend$,
 which becomes null on the event horizon.  Thus, the energy we control
has a weight which vanishes linearly at $r=\rp$. The other works use a
different timelike vector field, which includes some of the
horizon-penetrating vector field, first introduced in
\cite{DafermosRodnianski:RedShiftSchwarzschild}. This is denoted $Y$
\cite{DafermosRodnianski:KerrEnergyBound} or $X_2$
\cite{TataruTohaneanu}. We add the contribution from such a vector
field as a separate step in appendix \ref{S:DRVector}.

To overcome problem \ref{Problem:InsufficientClassicalSymmetries},
neither \cite{DafermosRodnianski:KerrEnergyBound} nor
\cite{TataruTohaneanu} use $\OpQ$ to generate higher energies. Away
from the event horizon, they use the symmetries $\dt^2$, $\dt\dphi$,
and $\dphi^2$ and the fact that $\solu$ satisfies the wave equation.
Near the event horizon, they generate higher energies using $\dt$ and
a horizon-penetrating, radial vector field (e.g. $Y$ in
\cite{DafermosRodnianski:KerrEnergyBound}). This is possible because
of a favourable sign in the error terms arising from the failure of
the radial vector field to be a symmetry.

Finally, to overcome problem \ref{Problem:ComplicatedTrapping}, both
of \cite{DafermosRodnianski:KerrEnergyBound,TataruTohaneanu} use a pseudo-differential Morawetz multiplier, as explained in subsection \ref{SS:ResultForGeodesics}. We have avoided these in favour of local differential operators. 

Less importantly, both avoid surfaces of constant $t$ in favour of
surfaces and coordinates which go through the event horizon. Since
vector-field arguments can be deformed from one surface to another,
this is a minor difference; however, the lower-order coefficients in the momenta,
$\MMTTscalar$, slightly complicate this. Although all known Morawetz
arguments have, in some sense, a troublesome lower-order term,
\cite{DafermosRodnianski:KerrEnergyBound,TataruTohaneanu} use a
different construction so that they can use positivity arising from
$Y$ or $X_2$, instead of the Hardy estimate we use to control the
negativity in \eqref{eq:IntroNeedHardy}.

The structure of this paper is as follows. Section
\ref{S:Preliminaries} introduces some preliminary results and further
notation. Section \ref{S:BoundedEnergy} contains the main argument of
this paper; in this section, we expand the energy associated with
$\vecTBlend$ and prove the Morawetz estimate using the
symmetry-indexed vector fields. Finally, a brief appendix reviews how
to derive nondegenerate energy estimates from the main estimates of
this paper.

%%%%%%%%%%%%%%%%%%%%%%%%%%%%%%%%%%%%%%%%%%%%%%%%%%%%%%%%%%%%%%%%%%%%%%%%%
\section{Notation and preliminaries}
\label{S:Preliminaries}
In this section, we present some more notation and basic estimates which we shall use throughout the paper. 

To begin, we note that we take $M>0$ as fixed. In a statement about
the existence of a sufficiently small bound $\bar{a}$ for which an
estimate holds for $|a|\leq\bar{a}$, it is understood that the upper
bound $\bar{a}$ depends on $M$. Similarly, in estimates, $C$ is used
to denote an absolute constant or a constant which depends only on
$M$. The notation $x\lesssim y$ means $x\leq Cy$, and the notation
$x\eqsim y$ means $x\lesssim y$ and $y\lesssim x$. All objects are smooth unless otherwise stated. 

In informal discussions, if $\OpSetGeneral$ is a set of operators,
then $\OpSetGeneral\psi$ will typically refer to expressions of the form $X\psi$ for $X\in\OpSetGeneral$. Similar notation is defined precisely in certain contexts in the remainder of this section. 

In the remainder of this paper, unless otherwise stated, Greek indices refer to components in the Boyer-Lindquist coordinates $(t,r,\theta,\phi)$. In the appendix, Greek indices refer to the Kerr coordinates denoted there, $(\tk,\rk,\hk,\pk)$.

\subsection{Proof of the general deformation formula}
In analogy with the definition of higher-order energies in equation
\eqref{eq:DefHigherEnergy}, we define, for a general set of operators
$\OpSetGeneral$, for a vector field $\vecY$, and a function $\solu$,
the higher-order momentum and energy as
\begin{align*}
\GenMomentum{\vecY}[\OpSetGeneral\solu]{}&{}=\sum_{\OpGeneral\in\OpSetGeneral} \GenMomentum{\vecY}[\OpGeneral\solu], \\
\GenEnergy{\vecY}[\OpSetGeneral\solu]{}&{}=\sum_{\OpGeneral\in\OpSetGeneral} \GenEnergy{\vecY}[\OpGeneral\solu]. 
\end{align*}

We now prove the validity of the general deformation formula \eqref{eq:deform:twovec}. In addition, we allow an additional positive function $\Omega$ to be introduced, which allows for many calculations to be simplified. 

\begin{lemma}
\label{Lemma:DivOfP}
If $U$ is an open set in a general $4$-dimensional Lorentzian
manifold, $(x^\alpha)$ is a system of coordinates on $U$, $\solu\in C^2(U)$ is a solution of $\nabla^\alpha\nabla_\alpha\solu=0$, $\vecX$ is a vector field, $\MMTTscalar$ is a function, and $\Omega$ is a positive function,
then the divergence of the associated momentum is
\begin{align}
\nabla_\alpha \GenMomentum{\pairXq}[\solu]^\alpha
{}&{}=-\frac{\Omega^2}{2}\Lie_{\vecX}(\Omega^{-2}\gMetric^{\alpha\beta})(\nabla_\alpha\solu)(\nabla_\beta\solu) \nonumber\\
{}&{}\quad+\left(-\frac{\Omega^{-2}}{2}\nabla_\alpha(\Omega^{2}\vecX^\alpha)+\MMTTscalar\right) (\nabla^\gamma\solu)(\nabla_\gamma\solu) 
-\frac12(\nabla_\alpha\nabla^\alpha\MMTTscalar)\solu^2 .
\label{eq:BasicDeformation}
\end{align}

Furthermore, if $\{\CQA_\ua\}$ is a set of symmetry operators for the wave equation, $\STwovecX=(\{\vecX^{\ua\ub}\},\{\MMTTscalar^{\ua\ub}\})$ is a pair consisting of symmetric collections of double-indexed vectors and scalars, and $\solu$ is a solution of the wave equation, then
\begin{align}
\nabla_\alpha \GenMomentum{\STwovecX}[\solu]^\alpha
{}&{}=-\frac{\Omega^2}{2}\Lie_{\vecX^{\ua\ub}}(\Omega^{-2}\gMetric^{\alpha\beta})(\nabla_\alpha\CQA_{\ua}\solu)(\nabla_\beta\CQA_{\ub}\solu) \nonumber\\
{}&{}\quad+\left(-\frac{\Omega^{-2}}{2}\nabla_\alpha(\Omega^{2}\vecX^{\ua\ub\alpha})+\MMTTscalar^{\ua\ub}\right) (\nabla^\gamma\CQA_{\ua}\solu)(\nabla_\gamma\CQA_{\ub}\solu) \nonumber\\
{}&{}\quad-\frac12(\nabla_\alpha\nabla^\alpha\MMTTscalar^{\ua\ub})(\CQA_{\ua}\solu)(\CQA_{\ub}\solu) .
\label{eq:DivOfP}
\end{align}
\end{lemma} 
\begin{proof}
By direct computation
\begin{align*}
\nabla_\alpha \GenMomentum{\pairXq}[\solu]^\alpha
{}&{}= (\nabla_\alpha \StressEnergyGeneral^{\alpha\beta})\vecX^\beta +\StressEnergyGeneral_{\alpha\beta}\nabla^{\alpha}\vecX^\beta \\
{}&{}\quad+\MMTTscalar\solu\nabla^\alpha\nabla_\alpha\solu+\MMTTscalar (\nabla^\alpha\solu)(\nabla_\alpha\solu) -\frac12(\nabla_\alpha\nabla^\alpha\MMTTscalar)\solu^2 .
\end{align*}
Since $\solu$ solves the wave equation, $\nabla_\alpha\StressEnergyGeneral^{\alpha\beta}=0$ and $\nabla_\alpha\nabla^\alpha\solu=0$. Expanding the energy-momentum tensor, one finds
\begin{align*}
\nabla_\alpha \GenMomentum{\pairXq}[\solu]^\alpha
{}&{}=  \nabla^{(\alpha}\vecX^{\beta)}(\nabla_\alpha\solu)(\nabla_\beta\solu) \\
{}&{}\quad+\left(-\frac12(\nabla_\alpha\vecX^\alpha)+\MMTTscalar\right) (\nabla^\alpha\solu)(\nabla_\alpha\solu) -\frac12(\nabla_\alpha\nabla^\alpha\MMTTscalar)\solu^2 .
\end{align*}
Since $\Omega$ is positive, there is a well-defined inverse. Inserting $1=\Omega^2\Omega^{-2}$ into the derivative of $\vecX$ and using the formula $\nabla^{(\alpha}\vecX^{\beta)}=(-1/2)\Lie_\vecX\gMetric^{\alpha\beta}$, one finds $\nabla^{(\alpha}\vecX^{\beta)}$ $=(-1/2)\Lie_\vecX\gMetric^{\alpha\beta}$ $=(-1/2)\Omega^2 \Lie_\vecX(\Omega^{-2}\gMetric^{\alpha\beta})-(1/2)\Omega^{-2}\gMetric^{\alpha\beta}(\vecX\Omega^{2})$, so
\begin{align*}
\nabla_\alpha \GenMomentum{\pairXq}[\solu]^\alpha
{}&{}=  -\frac{\Omega^2}{2}\Lie_{\vecX}(\Omega^{-2}\gMetric^{\alpha\beta})(\nabla_\alpha\solu)(\nabla_\beta\solu) \\
{}&{}\quad+\left(-\frac12\left((\Omega^{-2}(\vecX\Omega^{2})\right)+\nabla_\alpha\vecX^\alpha)+\MMTTscalar\right) (\nabla^\alpha\solu)(\nabla_\alpha\solu) -\frac12(\nabla_\alpha\nabla^\alpha\MMTTscalar)\solu^2 \\
{}&{}=  -\frac{\Omega^2}{2}\Lie_{\vecX}(\Omega^{-2}\gMetric^{\alpha\beta})(\nabla_\alpha\solu)(\nabla_\beta\solu) \\
{}&{}\quad+\left(-\frac{\Omega^{-2}}{2}\nabla_\alpha(\Omega^{2}\vecX^\alpha)+\MMTTscalar\right) (\nabla^\alpha\solu)(\nabla_\alpha\solu) -\frac12(\nabla_\alpha\nabla^\alpha\MMTTscalar)\solu^2 .
\end{align*}

The second part of the theorem follows by replacing
$(\vecX,\MMTTscalar)$, $(\nabla_\alpha\solu)(\nabla_\beta\solu)$, and
$\solu^2$ by $(\vecX^{\ua\ub},\MMTTscalar^{\ua\ub})$,
$(\nabla_\alpha\CQA_\ua\solu)(\nabla_\beta\CQA_\ub\solu)$, and
$(\CQA_\ua\solu)(\CQA_\ub\solu)$ respectively and then using the
$\ua\ub$ symmetry of $\vecX^{\ua\ub}$ and $\MMTTscalar^{\ua\ub}$.

\end{proof}

\subsection{Simplifying rescalings}
\label{SS:SigmaRescaling}
It is convenient to introduce the following reference volume forms
\begin{align*}
\diTwo ={}&{} \Svol \di\theta\di\phi , {}&{}
\Svol ={}&{} \sin\theta ,{}&{}
\diThree={}&{} \diTwo\di r, {}&{}
\diFour={}&{} \diTwo\di r \di t .
\end{align*}

It so happens that the Boyer-Lindquist coordinates allow the second-order symmetry operators to be expressed easily in terms of coordinate partial derivatives and $\Svol$
\begin{align*}
\CQA_\ua ={}&{} \frac1{\Svol} \partial_\alpha \Svol\CQA_\ua^{\alpha\beta}   \partial_\beta .
\end{align*}
All other operators built from these can be similarly expanded. For
example, the operator $\CurlyR$ defined in equation
\eqref{eq:OpCurlyR} can be written as 
\begin{align*}
\CurlyR ={}&{} \frac1{\Svol} \partial_\alpha \Svol\CurlyR^{\alpha\beta}   \partial_\beta ,
\end{align*}
where $\CurlyR^{\alpha\beta}$ is defined in
\eqref{eq:DefCurlyRalphabeta}. 
Similarly, the contravariant form of the metric can be written as
\begin{align}
\KSigma\gMetric^{\alpha\beta}
{}&{}=\KDelta \dr^\alpha\dr^\beta +\frac{1}{\KDelta}\CurlyR^{\alpha\beta} .
\label{eq:SigmaMetric}
\end{align}
This eliminates all $\theta$ dependence, except that arising through $\TensorQ^{\alpha\beta}$. 

Careful applications of the factor $\KSigma$ can be used in many applications to either eliminate $\theta$ or leave only $\Svol$. The volume element for $\gMetric_{\alpha\beta}$ in Boyer-Lindquist coordinates is given by 
\begin{align*}
\gVol {}&{}= \sqrt{-\det \gMetric_{\alpha\beta}} = \KSigma \sin\theta .
\end{align*}
Thus, divergences can be written as
\begin{align}
\KSigma\nabla_\alpha\vecX^\alpha
{}&{}=\KSigma\frac{1}{\gVol} \partial_\alpha \gVol \vecX^\alpha
=\frac{1}{\Svol}\partial_\alpha \Svol\KSigma\vecX^{\alpha} .
\label{eq:RescaledDivergence}
\end{align}
Similarly, $\KSigma\gWave$ can be written as
\begin{align*}
\KSigma \gWave = \KSigma \frac{1}{\gVol} \partial_\alpha \gVol \gMetric^{\alpha\beta} \partial_\beta = \frac{1}{\Svol} \partial_\alpha \Svol \KSigma \gMetric^{\alpha\beta} \partial_\beta .
\end{align*}
From the formula for $\KSigma\gMetric$ in equation \eqref{eq:SigmaMetric}, this eliminates all $\theta$ and $\dtheta$ terms except for those arising from $\OpQ$. 

In the deformation formulas \eqref{eq:BasicDeformation}-\eqref{eq:DivOfP}, we will make the choice
\begin{align*}
\Omega^{-2}{}&{}=\KSigma .
\end{align*}
This yields the Lie derivative of $\KSigma\gMetric^{\alpha\beta}$ and the divergence of $\KSigma^{-1}\vecX$. These can be simplified using equations \eqref{eq:SigmaMetric} and \eqref{eq:RescaledDivergence} respectively.

\newpageForSubsection
%%%%%%%%%%%%%%%%%%%%%%%%%%%%%%%%%%%%%%%%%%%%%%%%%%%%%%%%%%%%%%%%%%%%%%
\subsection{The $3+1$ decomposition}
The surfaces of constant $t$, $\hst{t}$, are spacelike since they are spanned by the spacelike vector fields $\dr$, $\dtheta$, and $\dphi$. Thus, the $1$-form $\di t$ is timelike in the exterior. Its length is $\gMetric(\di t,\di t)=\gMetric^{tt}=\frac{\KPi}{\KDelta\KSigma}$. For our purposes, it will be convenient to rescale this by $(\gMetric^{tt})^{-1}$, so that the component in the $\dt$ direction is $1$. Thus, we introduce
\begin{align*}
\vecTperp={}&{} \dt +\omegaperp\dphi , \\
\omegaperp={}&{} \frac{\gMetric^{t\phi}}{\gMetric^{tt}}=\frac{2aMr}{\KPi} ,
\end{align*}
which has length
\begin{align*}
\gMetric(\vecTperp,\vecTperp)=(\gMetric^{tt})^{-1}=-\frac{\KDelta\KSigma}{\KPi} .
\end{align*}

The vector field $\vecTperp$ is timelike in the exterior, and it extends
continuously to the event horizon and the bifurcation sphere. In fact, it
extends smoothly through the event horizon and the bifurcation
sphere.\footnote{The vector fields $\dt$ and $\dphi$ are known to extend
  smoothly through the bifurcation sphere 
\cite{ONeill}.}
This vector field extends to the null tangent vector on the event horizon and to axial rotation (with coefficient $\omegaH$) on the bifurcation sphere. 

To calculate the flux through hypersurfaces of constant $t$, one needs the normal volume element, 
\begin{align*}
\diNormal^\alpha
{}&{}=-\vecNormalSigmat^\alpha \sqrt{\gMetric_{rr}\gMetric_{\theta\theta}\gMetric_{\phi\phi}}\di r\di\theta\di\phi \\
%{}&{}=\vecTperp^\alpha (\gMetric^{tt})^{1/2}\sqrt{\gMetric_{rr}\gMetric_{\theta\theta}\gMetric_{\phi\phi}}\di r\di\theta\di\phi \\
%{}&{}=\vecTperp^\alpha \gMetric^{tt} \sqrt{-\text{det}g} \di r\di\theta\di\phi \\
%\diNormal^\alpha 
{}&{}= -\vecTperp^\alpha \frac{\KPi}{\KDelta}\di r\sin\theta\di\theta\di\phi .
\end{align*}
(This can be computed using that $\vecNormalSigmat=\vecTperp (-\gMetric(\vecTperp,\vecTperp))^{-1/2}$ $=\vecTperp (\gMetric^{tt})^{1/2}$; that, from the formula for inverting $2\times2$ matrices, $\gMetric^{tt}=\gMetric_{\phi\phi}(\gMetric_{tt}\gMetric_{\phi\phi}-\gMetric_{t\phi}^2)^{-1} $; and that, from the formula for determinants, $-\KSigma^2\sin^2\theta=\text{det}g=\gMetric_{rr}\gMetric_{\theta\theta}(\gMetric_{\phi\phi}\gMetric_{tt}-\gMetric_{t\phi}^2)=\gMetric_{rr}\gMetric_{\theta\theta}\gMetric_{\phi\phi}/\gMetric^{tt}$.)

For certain calculations, a contravariant form of the metric is more useful, in which case, we will use that, with our sign conventions
\begin{align*}
\diNormal_\alpha{}&{}=\KSigma \di t_\alpha \diThree .
\end{align*}

\newpageForSubsection
%%%%%%%%%%%%%%%%%%%%%%%%%%%%%%%%%%%%%%%%%%%%%%%%%%%%%%%%%%%%%%%%%%%%%%%%%%%
\subsection{Pointwise Norms}
\label{SS:Norms}

First, we introduce some notation for angular derivatives. We typically use $(\theta,\phi)$ for coordinates on the sphere, but occasionally use $\omega\in S^2$ to avoid coordinate singularities at the poles. We use $\dAng$ to denote the angular gradient and $\SLap$ for the Laplacian on the unit sphere. For two vectors on the sphere, there is an inner product defined using the standard metric on the unit sphere. Extending this notation to differential operators in the standard way, for a function $f$, $\SLap f=\Svol^{-1}\dAng\cdot(\Svol\dAng f)$. 

The Boyer-Lindquist coordinates induce coordinates $(\theta,\phi)$ on
the constant $(t,r)$ surfaces. This defines a diffeomorphism to the
unit sphere in $\Reals^3$ using the standard spherical
coordinates. This diffeomorphism, defined in the $(\theta,\phi)$
chart, extends smoothly to the entire sphere. This allows us to treat
$\dAng$ and $\SLap$ as operators defined in the Kerr spacetime. We use
$\Li$ for the pullback under this diffeomorphism of the rotation
vector fields about the coordinate axes. With the exception of
$\LsubThree=\dphi$, these are not symmetries in the Kerr spacetime. We
use $\VFAngAlg=\{\Li\}$ to denote the set of these rotations, and we
use $\VFEffSym$ for $\{\dt,\Li\}$.

Here and throughout the rest of the paper, we say that a vector field
has smooth angular components if, for fixed $r$ and $t$, the
contraction of the vector field with any smooth $1$-form on the sphere
produces a smooth function. Because of the coordinate singularity in
the $(\theta,\phi)$ coordinates, this does not assert that the
$\theta$ and $\phi$ components of the vector field are smooth. The
angular gradient of a smooth function has smooth angular
components. 

Given a set of differential operators, $\OpSetGeneral$, we use the notation
\begin{align*}
| \solu |_\OpSetGeneral^2 ={}&{} |\OpSetGeneral \solu |^2 =\sum_{\OpGeneral\in\OpSetGeneral} |\OpGeneral \solu|^2 .
\end{align*}
If no set is specified, simply an index, we mean
\begin{align*}
\normPtwiseTn{n}{\solu}^2 ={}&{} \sum_{i=0}^n |\SymGeneraln{i} \solu |^2 , 
\end{align*}
where $\SymGeneraln{i}$ is the set of generators of the order-$i$ symmetries given in equation \eqref{eq:KerrHigherSymmetries}. We will refer to $\normPtwiseTn{n}{\solu}$ as the order-$n$ pointwise norm of $\solu$. When the $n$ is clear from context, we will simply refer to this as the norm of $\solu$. 

\begin{lemma}[Spherical Sobolev estimate using symmetries]
\label{Lemma:STwoSobolev}
There is a constant, $\ConstantSTwoSobolev$, such that for all $(t,r)\in\Reals\times(\rp,\infty)$, if $\solu$ is sufficiently smooth that the quantity on the right is bounded, then
\begin{align*}
\sup_{(t,r)\times S^2} |\solu|^2 \leq \ConstantSTwoSobolev \int_{(t,r)\times S^2} \normPtwiseTn{2}{\solu}^2 \diTwo .
\end{align*}
\begin{proof}
Recall that we use $\Svol$ to denote $\sin\theta$ and $\SLap$ to
denote the spherical Laplacian, which takes the form
\begin{align*}
\SLap ={}&{} \frac{1}{\Svol}\dtheta\Svol\dtheta +\frac{1}{\Svol^2}\dphi^2 . 
\end{align*}
The absolute value of the spherical Laplacian of $u$ can be estimated by
\begin{align*}
|\SLap \solu|
={}&{} \left| \left(\frac{1}{\Svol}\dtheta\Svol\dtheta +\cot^2\theta\dphi^2 +\dphi^2\right) \solu \right| \\
\leq{}&{} \left| \left(\frac{1}{\Svol}\dtheta\Svol\dtheta +\cot^2\theta\dphi^2\right)\solu\right| +\left|\dphi^2 \solu \right| \\
\leq{}&{} |\OpQ \solu| + a^2\sin^2\theta|\dt^2 \solu| +|\dphi^2 \solu| \\
\lesssim{}&{} |\SOSym \solu| .
\end{align*}
By a standard, spherical, Sobolev estimate, 
\begin{align*}
|\solu|_{L^\infty(S^2)}^2 \lesssim \int_{S^2} \left(|\SLap \solu|^2 + |\solu|^2\right) \diTwo .
\end{align*}
Since the integrand on the right is bounded by $\normPtwiseTn{2}{\solu}$, the desired estimate holds with a uniform constant in $(t,r)$. 
\end{proof}
\end{lemma}

In subsection \ref{SS:ChoosingTheWeights}, we also require the following operators and the associated weaker norms. 

\begin{definition}
For $\epsilon\geq0$, let
\begin{align*}
\OpL ={}&{} \dt^2 +\OpQ +\dphi^2 , \\
\OpLEpsilon{\epsilon} ={}&{} \epsilon\dt^2 +\OpQ +\dphi^2 , 
\end{align*}
and
\begin{align*}
\normPtwiseTnEpsilon{2}{\epsilon}{\solu}^2
={}&{} \epsilon |\dt^2 \solu|^2 
+ (1+\epsilon)\left(|\dt\dtheta \solu|^2 +\frac{1}{\Svol^2}|\dt\dphi \solu|^2\right) 
+ |\SLap \solu|^2 , \\
\normPtwiseTnEpsilon{3}{\epsilon}{\solu}^2
={}&{} \epsilon^2 |\dt^3 \solu|^2 
+(2\epsilon+\epsilon^2) |\dt^2\dAng \solu|^2  
+(1+2\epsilon) |\dt\SLap \solu|^2 +|\dAng\SLap \solu|^2 . 
\end{align*}
We also introduce the homogeneous norms, generated from the previous norm by taking $\epsilon=1$, 
\begin{align*}
\normPtwiseTnHomo{n}{\solu} .
\end{align*}
\end{definition}

In these norms, there are coefficients that are not just monomial in $\epsilon$, such as the $(1+\epsilon)$ in $\normPtwiseTnEpsilon{2}{\epsilon}{\solu}$. These permit exact equality in some estimates below. 

\begin{lemma}[The $\OpLEpsilon{\epsilon}\OpL$ estimate]
\label{Lemma:LepsilonL}
There is a positive constant $\ConstantLEpsilonL$ such that, for $\epsilon\in(0,1)$, if $\psi$ is smooth, then\footnote{The index on $\TimeAndAngularDerivativesTwoPtThreea{\psi}$ refers to equation \eqref{eq:ExplicitLepsilonLExpansionControlsTwoDerivatives}.} 
\begin{align}
\left|(\OpLEpsilon{\epsilon} \psi)(\OpL \psi)
- \normPtwiseTnEpsilon{2}{\epsilon}{\solu}^2
+\frac{1}{\Svol}\partial_\alpha\left(\Svol(1+\epsilon)\TimeAndAngularDerivativesTwoPtThreea{\psi}^\alpha\right)\right|
\lesssim a^2\normPtwiseTnEpsilon{2}{1}{\psi}^2 .
\label{eq:ExplicitLepsilonLExpansionControlsTwoDerivatives}
\end{align}
where
$\TimeAndAngularDerivativesTwoPtThreea{\psi}^t=(\dt\psi)(\SLap\psi)$,
$\TimeAndAngularDerivativesTwoPtThreea{\psi}^r=0$,
$\TimeAndAngularDerivativesTwoPtThreea{\psi}^\theta=-(\dt\psi)(\dtheta\dt\psi)$,
and
$\TimeAndAngularDerivativesTwoPtThreea{\psi}^\phi=-(\dt\psi)(\dphi\dt\psi)/\Svol^2$. 
%\begin{align*}
%\TimeAndAngularDerivativesTwoPtThreea^t{}&{}= (\dt\psi)(\SLap\psi) ,{}&{}
%\TimeAndAngularDerivativesTwoPtThreea^r{}&{}= 0,{}&{}
%\TimeAndAngularDerivativesTwoPtThreea^\theta{}&{}= -(\dt\psi)(\dtheta\dt\psi),{}&{}
%\TimeAndAngularDerivativesTwoPtThreea^\phi{}&{}=
%-(\dt\psi)\frac{\dphi\dt\psi}{\Svol^2} ,\\
%\TimeAndAngularDerivativesTwoPtThreeb^t{}&{}= (\dt\psi)(\SLap\psi) ,{}&{}
%\TimeAndAngularDerivativesTwoPtThreeb^r{}&{}= 0,{}&{}
%\TimeAndAngularDerivativesTwoPtThreeb^\theta{}&{}= 0, {}&{}
%\TimeAndAngularDerivativesTwoPtThreeb^\phi{}&{}=
%-(\dt\psi)\frac{\dphi\dt\psi}{\Svol^2} .
%\end{align*} 

Furthermore, 
\begin{align}
\left| (\OpLEpsilon{\epsilon}\dt\psi)^2 +(\OpLEpsilon{\epsilon}\dAng\psi)^2-\normPtwiseTnEpsilon{3}{\epsilon}{\psi}^2 +2\epsilon\frac{1}{\Svol}\partial_\alpha(\Svol\TimeAndAngularDerivativesTwoPtThreea{\VFEffSym\psi}^\alpha) \right| 
\lesssim{}&{} a^2 \normPtwiseTnEpsilon{3}{1}{\psi}^2 ,
\label{eq:ExplicitLepsilonLExpansionEqualityThirdOrder} 
\end{align}
where $\TimeAndAngularDerivativesTwoPtThreea{\VFEffSym\psi}$ denotes
$\sum_{\OpGeneral\in\VFEffSym}\TimeAndAngularDerivativesTwoPtThreea{\OpGeneral\psi}$. 
\end{lemma}

\begin{proof}
By direct computation, 
\begin{align*}
(\OpLEpsilon{\epsilon}\psi)(\OpL\psi)
-((\epsilon\dt^2+\SLap)\psi)((\dt^2+\SLap)\psi)
{}&{}=(a^2\sin^2\theta\dt^2\psi)(\OpL\psi)\\
{}&{}\quad+((\epsilon\dt^2+\SLap)\psi)(a^2\sin^2\theta\dt^2\psi) ,
\end{align*}
so
\begin{align*}
\left|(\OpLEpsilon{\epsilon}\psi)(\OpL\psi)
-((\epsilon\dt^2+\SLap)\psi)((\dt^2+\SLap)\psi)\right|
{}&{}\lesssim a^2\normPtwiseTnEpsilon{2}{1}{\psi}^2 .
\end{align*}
We now expand $((\epsilon\dt^2+\SLap)\psi)((\dt^2+\SLap)\psi)$ as
\begin{align*}
((\epsilon\dt^2+\SLap)\psi)((\dt^2+\SLap)\psi)
{}&{}=\epsilon(\dt^2\psi)+(1+\epsilon)(\dt^2\psi)(\SLap\psi)+(\SLap\psi)^2, 
\end{align*}
and simplify the cross-term by gathering total derivatives
\begin{align*}
(\dt^2\psi)(\SLap\psi)
{}&{}=-(\dt\psi)(\SLap\dt\psi) +\dt((\dt\psi)(\SLap\psi)) \\
{}&{}=|\dAng\dt\psi|^2 +\dt((\dt\psi)(\SLap\psi))
-\frac{1}{\Svol}\dAng\cdot((\dt\psi)(\dAng\dt\psi)) \\
{}&{}=|\dAng\dt\psi|^2
+\frac{1}{\Svol}\partial_\alpha\left(\Svol\TimeAndAngularDerivativesTwoPtThreea{\psi}^\alpha\right) .
\end{align*}
Note that it was crucial to gather the total derivatives first in $t$ and then
in the angular directions, so that the desired bound on
$\TimeAndAngularDerivativesTwoPtThreea{\psi}^t$ holds. This completes
the proof of estimate \eqref{eq:ExplicitLepsilonLExpansionControlsTwoDerivatives}.

The proof of estimate
\eqref{eq:ExplicitLepsilonLExpansionEqualityThirdOrder} 
follows the same steps. First, there is the simplification
from
\begin{align*}
\left|(\OpLEpsilon{\epsilon}\dt\psi)^2 +(\OpLEpsilon{\epsilon}\dAng\psi)^2-\left(((\epsilon\dt^2+\SLap)\dt\psi)^2+((\epsilon\dt^2+\SLap)\dAng\psi)^2\right)\right|
{}&{}\lesssim a^2 \normPtwiseTnEpsilon{3}{1}{\psi}^2 .
\end{align*}
Second, the simplified term can be expanded as
\begin{align*}
((\epsilon\dt^2+\SLap)\dt\psi)^2+((\epsilon\dt^2+\SLap)\dAng\psi)^2
{}&{}=\epsilon^2(\dt^3\psi) +2\epsilon(\dt^3\psi)(\SLap\dt\psi) +(\SLap\dt\psi)^2
  \\
{}&{}\quad+\epsilon^2(\dt^2\dAng\psi)^2+2\epsilon(\dt^2\dAng\psi)(\SLap\dAng\psi)+(\SLap\dAng\psi)^2
  . 
\end{align*}
Third, the mixed factors can be written as perfect squares plus total
derivatives
\begin{align*}
(\dt^3\psi)(\SLap\dt\psi)
{}&{}=(\dt^2\dAng\psi)^2+\frac{1}{\Svol}\partial_\alpha(\Svol\TimeAndAngularDerivativesTwoPtThreea{\dt\psi}^\alpha)
  , \\
(\dt^2\dAng\psi)(\SLap\dAng\psi)
{}&{}=(\SLap\dt\psi)\sum_{\OpGeneral\in\VFAngAlg}\frac{1}{\Svol}\partial_\alpha(\Svol\TimeAndAngularDerivativesTwoPtThreea{\OpGeneral\psi}^\alpha) .
\end{align*}

\end{proof}

An important consequence of lemma \ref{Lemma:LepsilonL} is that the
$(2,\epsilon)$ norm is dominated by $(\OpLEpsilon{\epsilon} \psi)(\OpL \psi)$ plus a divergence term and small error terms. In
the divergence term, the time component satisfies
$|\TimeAndAngularDerivativesTwoPtThreea{\psi}^t|\lesssim |\dt\solu|
|\SOSym\solu|$, there is no $r$ component, and the vector field
$\TimeAndAngularDerivativesTwoPtThreea{\psi}$ has smooth angular
components.

\newpageForSubsection
%%%%%%%%%%%%%%%%%%%%%%%%%%%%%%%%%%%%%%%%%%%%%%%%%%%%%%%%%%%%%%%%%%%%%%%%%%%
\subsection{Further notation}
\label{SS:FurtherNotation}

We use the notation
\begin{align*}
f = O(r^p)
\end{align*}
to mean that there is a constant, uniformly in $a$ in some small interval of $a$ values containing $0$, such that for all $r>\rp $, $|f(r)| < C r^p$. We introduce also the notation 
\begin{align*}
f = O\left( \left(\frac{\KDelta}{r^2} \right)^q, r^p \right)
\end{align*}
to mean that there is a constant, uniformly in $a$ in some small interval of $a$ values containing $0$, such that for all $r>\rp $, 
\begin{align*}
|f(r)| < C \left(\frac{\KDelta}{r^2} \right)^q r^p .
\end{align*} 
Similarly, for functions $f$ of $t,r,\omega$, we use $f=O(r^p)$ or
$O\left( \left(\frac{\KDelta}{r^2} \right)^q, r^p \right)$ when the
same bounds hold with the condition $\forall t\in\Reals, r>\rp,
\omega\in S^2$ replacing $r>\rp$ and $f(t,r,\omega)$ replacing $f(r)$.
This measures the decay rate at $\rp $ and $\infty$. If $f$ is
continuous, this is all the information that is required to bound the
function.

For a set $X$, we use $\Indicator_X$ to denote the indicator function, which is identically one on $X$ and zero elsewhere. We define a function to be smooth on a closed interval if it is smooth in the interior and if all the derivatives are continuous up to the boundary. 

\newpageForSubsection
\section{The bounded-energy argument}
\label{S:BoundedEnergy}
In this section, we construct a bounded energy by first constructing
an almost conserved energy and then proving a Morawetz estimate to
control its growth. 

%%%%%%%%%%%%%%%%%%%%%%%%%%%%%%%%%%%%%%%%%%%%%%%%%%%%%%%%%%%%%%%%%%%%%%%%%%
% The blended energy
%%%%%%%%%%%%%%%%%%%%%%%%%%%%%%%%%%%%%%%%%%%%%%%%%%%%%%%%%%%%%%%%%%%%%%%%%%
\subsection{The blended energy}
\label{SS:TBlend}
Recall from \eqref{eq:IntroDefTBlend} that for $|a|$ sufficiently small, the vector field   
\begin{align*}
\vecTBlend = \dt +\fnBlend\omegaH \dphi  
\end{align*}
is timelike in the exterior and Killing outside the region $[\epsilonInvBlendLocation,\epsilonInvBlendLocation+M]$, since $\fnBlend$ is constant outside this region and decreases from one to zero inside this region. If we choose $\epsilonInvBlendLocation$ sufficiently large so that it corresponds to a larger value of $r$ than any orbiting null geodesic for our initial choice of small $|a|$, this property will remain true for any subsequent decrease in the range of $|a|$ we allow. For definiteness, we take $\epsilonInvBlendLocation=10M$, which is beyond the range of the orbiting null geodesics for any Kerr black hole. 

The vector field $\vecTBlend$ becomes null on the horizon, so the
associated energy degenerates there. In the following theorem, we
compare this with the energy associated with
$\vecTperp=(\KDelta\KSigma/\KPi)^{1/2}\vecNormalSigmat$ to make clear
that the rate of degeneration with respect to the normal is 
$(\KDelta/(r^2+a^2))^{1/2}$. We also provide a coordinate expression
which is useful for making estimates. The apparently singular
contribution to the energy from $\KDelta^{-1}(\vecTperp\solu)^2$ is in
fact vanishing, since the vector field $\vecTperp$ vanishes on the
bifurcation sphere at such a rate to exactly compensate for the factor
of $\KDelta^{-1}$, and in addition the form $\di r$, which appears in
$\diThree$, degenerates at a rate of $(\KDelta/(r^2+a^2))^{1/2}$ near the bifurcation sphere. 

\begin{lemma}
There is a positive $\epsilonSlowRotationNearEnergyPositive$ such that for $|a|\leq\epsilonSlowRotationNearEnergyPositive$, if $t\in\Reals$ and $\solu$ is smooth, then  $\vecTBlend$ is timelike and 
\begin{align}
\GenEnergy{\vecTperp}(t)
{}&{}\eqsim\int_{\hst{t}} \left(\frac{(r^2+a^2)^2}{\KDelta}(\vecTperp \solu)^2 
+\KDelta (\dr \solu)^2
+\OpQ^{\alpha\beta}(\partial_\alpha \solu)(\partial_\beta \solu) \right)\diThree ,\nonumber\\
{}&{}\eqsim\int_{\hst{t}}\left(\frac{(r^2+a^2)^2}{\KDelta}(\vecTperp \solu) 
+\KDelta (\dr \solu) 
+\KDelta (\dt \solu)^2 
+\sum_i |\Li \solu|^2 \right)\diThree , 
\label{eq:EnergyExpansion}\\
\GenEnergy{\vecTBlend}(t){}&{}\eqsim\GenEnergy{\vecTperp}(t) \nonumber.
\end{align}
Furthermore, if $\solu$ is a solution of the wave equation 
$\gWave \solu = 0$, 
%\eqref{eq:KerrWave}, 
then 
\begin{align}
\left|\KSigma\nabla_\alpha\GenMomentum{\vecTBlend}[\solu]^\alpha\right| 
={}&{} \KDelta\omegaH |\partial_r \fnBlend| |\partial_\phi \solu| |\partial_r \solu| .
\end{align}
\end{lemma}
\begin{proof}
Since $-\gMetric_{\alpha\beta}\vecTperp^\alpha\vecTperp^\beta=\KDelta\KSigma/\KPi$, the $\vecTperp$ energy is 
\begin{align*}
\GenEnergy{\vecTperp}
{}&{}=\int_{\hst{t}} \StressEnergyGeneral_{\alpha\beta}\vecTperp^{\alpha}\vecTperp^{\beta} \frac{\KPi}{\KDelta}\diThree ,\\
{}&{}=\int_{\hst{t}} \left((\vecTperp\solu)^2 -\frac12\gMetric(\vecTperp,\vecTperp)\gMetric^{\alpha\beta}(\partial_\alpha\solu)(\partial_\beta\solu)\right)\frac{\KPi}{\KDelta}\diThree\\
{}&{}=\int_{\hst{t}} \left(\frac{\KPi}{\KDelta}(\vecTperp\solu)^2 +\frac12\KSigma\gMetric^{\alpha\beta}(\partial_\alpha\solu)(\partial_\beta\solu)\right)\diThree .
\end{align*}
The integrand can be expanded as 
\begin{align*}
\frac{\KPi}{\KDelta}(\vecTperp\solu)^2+{}&{}\frac12\KSigma\gMetric^{\alpha\beta}(\partial_\alpha\solu)(\partial_\beta\solu)\\
={}&{} \frac12\left(\KDelta(\dr\solu)^2 +\frac{(r^2+a^2)^2}{\KDelta}(\vecTperp\solu)^2 +\OpQ^{\alpha\beta}(\partial_\alpha\solu)(\partial_\beta\solu) +(\dphi\solu)^2\right) \\
{}&{}-\frac{1}{2\KDelta}\left(4aMr-2\omegaperp(r^2+a^2)^2\right)(\dt\solu)(\dphi\solu) \\
{}&{}+\frac{1}{2\KDelta}\left(-a^2+(r^2+a^2)^2\omegaperp^2\right)(\dphi\solu^2)
-a^2\sin^2\theta(\vecTperp\solu)^2 . 
\end{align*}
Since the coefficients $4aMr-2\omegaperp(r^2+a^2)^2$ and
$-a^2+(r^2+a^2)^2\omegaperp^2$ vanish at $r=\rp$, are bounded by
factors that go uniformly to $0$ on bounded sets as $a\rightarrow0$, and grow as $r\rightarrow\infty$ no faster than $r$ and a constant respectively, for $|a|$ sufficiently small,
\begin{align*}
\GenEnergy{\vecTperp}(t)
{}&{}\eqsim\int_{\hst{t}} \left(\frac{(r^2+a^2)^2}{\KDelta}(\vecTperp \solu)^2 
+\KDelta (\dr \solu)^2
+\OpQ^{\alpha\beta}(\partial_\alpha \solu)(\partial_\beta \solu) \right)\diThree .
\end{align*}

A clearer bound on the angular derivatives can be obtained by noting that the $\dt^2$ term in $\OpQ$ has a bounded factor times $a^2$, so
\begin{align*}
\frac{(r^2+a^2)^2}{\KDelta}(\vecTBlend\solu)^2 +\OpQ^{\alpha\beta}(\partial_\alpha\solu)(\partial_\beta\solu) +(\dphi\solu)^2 
\gtrsim{}&{} (\dt\solu)^2+\OpQ^{\alpha\beta}(\partial_\alpha\solu)(\partial_\beta\solu) +(\dphi\solu)^2 \\
\gtrsim{}&{} \sum_i |\Li u|^2 .
\end{align*}

The $\vecTBlend$ energy can be estimated using the fact that $\vecTperp-\vecTBlend=(\omegaperp-\fnBlend\omegaH)\dphi$ is orthogonal to $\vecTperp$, so
\begin{align*}
\GenEnergy{\vecTperp}-\GenEnergy{\vecTBlend}
{}&{} =\int_{\hst{t}} (\omegaperp-\fnBlend\omegaH)(\dphi\solu)(\vecTperp\solu) \frac{\KPi}{\KDelta} \diThree .
\end{align*}
The coefficient $\omegaperp-\fnBlend\omegaH$ vanishes linearly at $r=\rp$, is bounded by a function that goes to zero uniformly as $a\rightarrow 0$, and goes to zero as $r\rightarrow\infty$ like $r^{-4}$, so, by a simple Cauchy-Schwarz estimate, one finds $|\GenEnergy{\vecTperp}-\GenEnergy{\vecTBlend}|\lesssim |a| \GenEnergy{\vecTperp}$, so $\GenEnergy{\vecTperp}\eqsim\GenEnergy{\vecTBlend}$. 

The divergence of the momentum can be estimated using equation
\eqref{eq:BasicDeformation}. Taking $\Omega^{-2}=\KSigma$ greatly
simplifies the terms on the right-hand side of equation
\eqref{eq:BasicDeformation}. For example, one
finds
$\Omega^{-2}\nabla_{\alpha}(\Omega^{2}\vecTBlend^\alpha)=\Svol^{-1}\partial_\alpha(\Svol\vecTBlend^\alpha)=0$. Similarly,
this choice of $\Omega$ eliminates the factor of $\KSigma$ when computing
$\Lie_{\vecTBlend}(\Omega^{-2}\gMetric^{\alpha\beta})$. Thus,
\begin{align*}
\nabla_\alpha \GenMomentum{\vecTBlend}^\alpha
{}&{}=\Omega^2 \Lie_{\vecTBlend}\left(\Omega^{-2}\gMetric^{\alpha\beta}\right)(\partial_\alpha\solu)(\partial_\beta\solu)
=\KSigma^{-1} (\dr\fnBlend)\KDelta\omegaH(\dr\solu)(\dphi\solu) .
\end{align*}
\end{proof}

Recall that we defined higher-order energies by
\begin{align*}
\GenEnergyOrder{\vecTBlend}{n+1}[\solu]={}&{} \sum_{i=0}^n \GenEnergy{\vecTBlend}[\SymGeneraln{i}\solu] ,
\end{align*}
where $\SymGeneraln{i}$ is the set of order-$i$ symmetries from \eqref{eq:KerrHigherSymmetries}. 

\begin{corollary}
\label{Corollary:EnergyGrowthByDerivatives} 
If $\solu$ is a solution of the wave equation 
$\gWave \solu = 0$, 
%\eqref{eq:KerrWave}
\begin{align}
\Dt \GenEnergyOrder{\vecTBlend}{n+1}[\solu] \leq{}&{} C\int_{\rp}^\infty \int_{S^2} \fnAtBlendLocation \normPtwiseTn{n}{\dr \solu} \normPtwiseTn{n}{\dphi \solu} r^2 \diThree ,
\label{eq:GrowthOfThirdOrderEnergy}
\end{align}
where the norms on the right are defined in subsection \ref{SS:Norms}. 
\begin{proof}
This follows from considering a symmetry operator $\CQA$ of order $i$,
applying the previous lemma to $\CQA\solu$, summing over the
operators $\CQA\in\SymGeneraln{i}$, and then summing in $i$. 
\end{proof}
\end{corollary}

\newpageForSubsection 
%%%%%%%%%%%%%%%%%%%%%%%%%%%%%%%%%%%%%%%%%%%%%%%%%%%%%%%%%%%%%%%%%%%%%%%%%%
% The local-decay Morawetz estimate
%%%%%%%%%%%%%%%%%%%%%%%%%%%%%%%%%%%%%%%%%%%%%%%%%%%%%%%%%%%%%%%%%%%%%%%%%%
%\section{The local-decay Morawetz estimate}
%\label{S:Morawetz}
\subsection{Set-up for radial vector fields and their fifth-order analogues} 
\label{SS:SetUpForRadial}

\begin{definition}
If $\fnMna$ and $\fnMnb$ are smooth functions of $r$ and the
parameters $M$ and  $a$, then we define the following single- and
double-indexed quantities
\begin{subequations}
\begin{align}
\CurlyRTilde^{\ua} ={}&{} \frac{\fnMna}{\KDelta}\CurlyR^{\ua} , \\
\DiffCurlyRTilde^{\ua}={}&{} \dr
\left(\frac{\fnMna}{\KDelta}\CurlyR^{\ua}\right) ,
\label{eq:DefDiffCurlyRTilde}\\
\DiffCurlyRTTilde^\ua={}&{}\fnMnb\frac{\fnMna^{1/2}}{\KDelta^{1/2}}\DiffCurlyRTilde^{\ua},\\
\DDiffCurlyRTTilde^{\ua}={}&{} \dr
\left(\fnMnb\frac{\fnMna^{1/2}}{\KDelta^{1/2}}\DiffCurlyRTilde^{\ua}\right)
.
\end{align}
\end{subequations}
\end{definition} 

These can be used to define a double-indexed family of vectors and
scalars which we shall use to prove a Morawetz estimate. 

\begin{definition}
Given smooth functions $\fnMna$ and $\fnMnb$ as above, 
the radial coefficients and reduced scalar functions are defined to be
\begin{align*}
\vecMrone  ={}&{} -\fnMna \fnMnb \DiffCurlyRTilde^{\ua} , {}&{}
\vecMr     ={}&{} -\fnMna \fnMnb \DiffCurlyRTilde^{(\ua}\OpL^{\ub)} , \\
\fnMpone   ={}&{} -\frac12(\dr \fnMna) \fnMnb \DiffCurlyRTilde^{\ua}  ,{}&{}
\fnMp      ={}&{} -\frac12(\dr \fnMna) \fnMnb \DiffCurlyRTilde^{(\ua}\OpL^{\ub)}  .
\end{align*}
The Morawetz vector fields and scalar functions are defined in terms
of these as
\begin{align*}
\vecMone {}&{}= \vecMrone\dr, {}&{}
\vecM={}&{} \vecMr \dr, \\
\fnMone={}&{} \frac12 \KSigma\nabla_\gamma (\KSigma^{-1}\vecMone^\gamma) -\fnMpone , {}&{}
\fnM ={}&{} \frac12 \KSigma \nabla_\gamma(\KSigma^{-1} \vecM^\gamma)
-\fnMp , 
\end{align*}

For simplicity, we introduce the following notation for the pair consisting of the previous sets of vector fields and functions, 
\begin{align*}
\pairM ={}&{} (\{\vecM\},\{\fnM\}) . 
\end{align*}
\end{definition}

\begin{lemma}
\label{Lemma:BasicMorawetzDeformation}
If $\solu$ is a solution to the wave equation 
$\gWave \solu = 0$, 
%\eqref{eq:KerrWave}, 
then the divergence of the momentum associated with these quantities is given by
\begin{align}
\KSigma\nabla_\alpha \GenMomentum{\pairM}[\solu]^\alpha
={}&{} \CurlyA^{\ua\ub}(\dr \CQA_\ua \solu)(\dr \CQA_\ub \solu)\nonumber\\
{}&{}+\CurlyUfour^{\ua\ub\alpha\beta}(\partial_\alpha\CQA_\ua \solu)(\partial_\beta \CQA_\ub \solu) \nonumber\\
{}&{}+\CurlyV^{\ua\ub}(\CQA_\ua \solu)(\CQA_\ub \solu) ,
\label{eq:DivPGenMorawetz}
\end{align}
where
\begin{align*}
\CurlyAone^{\ua} ={}&{} \fnMna^{1/2} \KDelta^{3/2}(-\DDiffCurlyRTTilde^{\ua}) ,& 
\CurlyA^{\ua\ub} ={}&{} \CurlyAone^{(\ua}\OpL^{\ub)} , \\
\CurlyU^{\ua\ub} ={}&{} \frac12\fnMnb \DiffCurlyRTilde^{\ua}\DiffCurlyRTilde^{\ub} , &
\CurlyUfour^{\ua\ub\alpha\beta}={}&{} \CurlyU^{\uc(\ua}\OpL^{\ub)}
\CQA_\uc^{\alpha\beta} ,\\
\CurlyVone^{\ua} ={}&{} \frac14 (\dr\KDelta\dr \fnMna (\dr \fnMnb\DiffCurlyRTilde^{\ua})) ,&
\CurlyV^{\ua\ub} ={}&{} \CurlyVone^{(\ua}\OpL^{\ub)} . 
\end{align*}

%\begin{align*}
%\CurlyAone^{\ua} ={}&{} \fnMna^{1/2} \KDelta^{3/2}
%(-\DDiffCurlyRTTilde^{\ua}) , 
%\\
%\CurlyVone^{\ua} ={}&{} \frac14 (\dr\KDelta\dr \fnMna (\dr \fnMnb\DiffCurlyRTil%de^{\ua})) ,\\
%\CurlyU^{\ua\ub} ={}&{} \frac12\fnMnb \DiffCurlyRTilde^{\ua}\DiffCurlyRTilde^{\ub} , \\
%\CurlyA^{\ua\ub} ={}&{} \CurlyAone^{(\ua}\OpL^{\ub)} , \\
%\CurlyV^{\ua\ub} ={}&{} \CurlyVone^{(\ua}\OpL^{\ub)} , \\
%\CurlyUfour^{\ua\ub\alpha\beta}={}&{} \CurlyU^{\uc(\ua}\OpL^{\ub)} \CQA_\uc^{\alpha\beta} .
%\end{align*}
\end{lemma}
\begin{proof}
In the formula for the divergence of the momentum, equation
\eqref{eq:DivOfP}, we choose $\Omega^{-2}=\KSigma$. Since
$\Omega^{-2}\gMetric^{\alpha\beta} = \KSigma\gMetric^{\alpha\beta} =
\KDelta\dr^\alpha\dr^\beta +\KDelta^{-1}\CurlyR^{\alpha\beta}$, this
choice of $\Omega$ eliminates $\KSigma$ when we need to compute the
Lie derivative along $\vecM$, enormously simplifying the
calculation. Furthermore, the term $\fnM$ has been chosen so that the
coefficient of
$(\nabla^\gamma\CQA_\ua\solu)(\nabla_\gamma\CQA_\ub\solu)$ is
$-\fnMp$. Thus, the divergence of the momentum is given by
\begin{align*}
\KSigma\nabla_\alpha\GenMomentum{\pairM}^\alpha
={}&{} \left( \KDelta (\dr\vecMr) -\frac12\vecMr(\dr\KDelta) \right) (\dr \CQA_\ua \solu)(\dr \CQA_\ub \solu) \\
{}&{}-\frac12\vecMr \left(\dr\left(\frac{\CurlyR^{\alpha\beta}}{\KDelta}\right)\right) (\partial_\alpha \CQA_\ua \solu)(\partial_\beta \CQA_\ub \solu) \\
{}&{}-\fnMp \KDelta (\dr \CQA_\ua \solu)(\dr \CQA_\ub \solu)
- \fnMp \frac{\CurlyR^{\alpha\beta}}{\KDelta} (\partial_\alpha \CQA_\ua \solu)(\partial_\beta \CQA_\ub \solu). \\
{}&{}-\frac12 (\KSigma\nabla^\alpha\nabla_\alpha \MMTTscalar^{\ua\ub}) (\CQA_\ua \solu)(\CQA_\ub \solu) .
\end{align*}
In the coefficient of the radial derivative terms, the part coming from the vector field can be rewritten as
\begin{align*}
\left( \KDelta (\dr\vecMr) -\frac12\vecMr(\dr\KDelta) \right) 
={}&{} \left(\dr\left(\frac{\vecMr}{\KDelta^{1/2}}\right)\right)\KDelta^{3/2} . 
\end{align*}

Expanding using the definitions of $\fnMna$, $\fnMnb$, and $\DiffCurlyRTilde$, we first note that $\MMTTscalar^{\ua\ub}=-(1/2)\dr(\fnMna\fnMnb\DiffCurlyRTilde^{\ua\ub})+(1/2)(\dr\fnMna)\fnMnb\DiffCurlyRTilde^{\ua\ub}=-(1/2)\fnMna\dr(\fnMnb\DiffCurlyRTilde^{\ua\ub})$. Thus, the divergence of the momentum is
\begin{align*}
\KSigma\nabla_\alpha\GenMomentum{\pairM}^\alpha
={}&{}-\fnMna^{1/2}\KDelta^{3/2} \left(\dr\left(\frac{\fnMna^{1/2}}{\KDelta^{1/2} }\fnMnb\DiffCurlyRTilde^{\ua}\right)\right)\OpL^\ub (\dr\CQA_\ua \solu)(\dr\CQA_\ub \solu) \\
{}&{}+\frac12\fnMnb\DiffCurlyRTilde^{\ua}\OpL^\ub\left( \dr\left( \frac{\fnMna\CurlyR^{\alpha\beta}}{\KDelta}\right)\right)(\partial_\alpha\CQA_\ua \solu)(\partial_\beta\CQA_\ub \solu) \\
{}&{}+\frac14 (\KSigma\nabla^\alpha\nabla_\alpha (\fnMna(\dr\fnMnb\DiffCurlyRTilde^{(\ua}))\OpL^{\ub)}) (\CQA_\ua \solu)(\CQA_\ub \solu) .
\end{align*}
The expression $\DiffCurlyRTilde$ was chosen so that it is exactly the derivative in the second term. Similarly, the quantity $\DDiffCurlyRTTilde$ was chosen so that it is the derivative in the first term. Thus, the total bulk term is
\begin{align*}
\KSigma\nabla_\alpha\GenMomentum{\pairM}^\alpha
={}&{} -\fnMna^{1/2} \KDelta^{3/2} \DDiffCurlyRTTilde^{\ua}\OpL^{\ub} (\dr \CQA_\ua \solu)(\dr \CQA_\ub \solu)\\
{}&{}+\frac12\fnMnb(\OpL^\ua\DiffCurlyRTilde^{\ub})\DiffCurlyRTilde^{\alpha\beta}(\partial_\alpha\CQA_\ua \solu)(\partial_\beta \CQA_\ub \solu) \\
{}&{}+\frac14 (\dr\KDelta\dr(\fnMna(\dr\fnMnb\DiffCurlyRTilde^{(\ua})))\OpL^{\ub)}  (\CQA_\ua \solu)(\CQA_\ub \solu) .
\end{align*}
Since $\DDiffCurlyRTTilde^\ua\OpL^\ub$ is contracted against a quantity which is symmetric in $\ua\ub$, it is not necessary to distinguish between $\DDiffCurlyRTTilde^\ua\OpL^\ub$ and $\DDiffCurlyRTTilde^{(\ua}\OpL^{\ub)}$. Substituting the definitions of $\CurlyA^{\ua\ub}$, $\CurlyU^{\ua\ub\alpha\beta}$, and $\CurlyV^{\ua\ub}$ gives the desired result. 
\end{proof}

%%%%%%%%%%%%%%%%%%%%%%%%%%%%%%%%%%%%%%%%%%%%%%%%%%%%%%%%%%%%%%%%%%%%%%%%%%%%
%% Rearrangements
\newpageForSubsection
\subsection{Rearrangements}
\label{SS:Rearrangements}
We rearrange the terms related to $\CurlyU$ to get a strictly positive contribution to the divergence. 

\begin{lemma}
\label{Lemma:rearrangements}
If $\solu$ is a solution to the wave equation 
$\gWave \solu = 0$, 
%\eqref{eq:KerrWave}, 
then 
\begin{align*}
\KSigma\nabla_\alpha \left( \GenMomentum{\pairM}[\solu]^\alpha +\MorBoundaryTermA{\pairM}[\solu]^\alpha \right)
={}&{} \CurlyA^{\ua\ub} (\dr\CQA_\ua \solu)(\dr\CQA_\ub \solu) \\
{}&{}+\CurlyU^{\ua\ub}\OpL^{\alpha\beta} (\partial_\alpha\CQA_\ua \solu)(\partial_\beta\CQA_\ub \solu) \\
{}&{}+\CurlyV^{\ua\ub} (\CQA_\ua \solu)(\CQA_\ub \solu) ,
\end{align*}
where $\CurlyA$, $\CurlyU$, and $\CurlyV$ are defined in lemma \ref{Lemma:BasicMorawetzDeformation} and
\begin{align*}
\KSigma\MorBoundaryTermA{\pairM}[\solu]^\alpha
={}&{}\left(\CurlyU^{\ua\ub}\OpL^{\alpha\beta}-\CurlyU^{\ua\ub\alpha\beta}\right) (\CQA_\ua \solu)(\partial_\beta\CQA_\ub \solu) .
\end{align*}
We refer to $\MorBoundaryTermA{\pairM}$ as the first boundary term. 
\end{lemma}
\begin{proof}
Starting from equation \eqref{eq:DivPGenMorawetz}, it is only the second term on the right side that needs to be manipulated. First, we rearrange the derivative term to get
\begin{align*}
\Svol\CurlyUfour^{\ua\ub\alpha\beta}(\partial_\alpha\CQA_\ua \solu)(\partial_\beta \CQA_\ub \solu) 
={}&{}\Svol\CurlyU^{\uc\ua}\OpL^\ub \CQA_\uc^{\alpha\beta}(\partial_\alpha\CQA_\ua \solu)(\partial_\beta \CQA_\ub \solu) \\
={}&{}-\CurlyU^{\uc\ua}\OpL^\ub (\CQA_\ua \solu)(\partial_\alpha\Svol\CQA_\uc^{\alpha\beta}\partial_\beta \CQA_\ub \solu) \\
{}&{}+\partial_\alpha(\Svol \CurlyU^{\uc\ua}\OpL^\ub \CQA_\uc^{\alpha\beta}(\CQA_\ua \solu)(\partial_\beta \CQA_\ub \solu)) 
.
\end{align*}
The first term on the right can be rewritten in terms of $\CQA_\uc$,
which can be commuted with $\CQA_\ub$, which in turn can be expanded
in partial derivatives as 
\begin{align*}
-\CurlyU^{\uc\ua}\OpL^\ub (\CQA_\ua \solu)(\partial_\alpha\Svol\CQA_\uc^{\alpha\beta}\partial_\beta \CQA_\ub \solu) 
%\\{}&{}+\partial_\alpha(\Svol \CurlyU^{\uc\ua}\OpL^\ub \CQA_\uc^{\alpha\beta}(\CQA_\ua \solu)(\partial_\beta \CQA_\ub \solu)) \\
={}&{}-\Svol\CurlyU^{\uc\ua}\OpL^\ub (\CQA_\ua \solu)(\CQA_\uc \CQA_\ub \solu) \\
%{}&{}+\partial_\alpha(\Svol \CurlyU^{\uc\ua}\OpL^\ub \CQA_\uc^{\alpha\beta}(\CQA_\ua u)(\partial_\beta \CQA_\ub \solu)) \\
={}&{}-\Svol\CurlyU^{\uc\ua}\OpL^\ub (\CQA_\ua \solu)(\CQA_\ub \CQA_\uc \solu) \\
%{}&{}+\partial_\alpha(\Svol \CurlyU^{\uc\ua}\OpL^\ub \CQA_\uc^{\alpha\beta}(\CQA_\ua \solu)(\partial_\beta \CQA_\ub \solu)) \\
={}&{}-\CurlyU^{\uc\ua}\OpL^\ub (\CQA_\ua \solu)(\partial_\alpha\Svol\CQA_\ub^{\alpha\beta}\partial_\beta \CQA_\uc \solu) .
%\\
%{}&{}+\partial_\alpha(\Svol \CurlyU^{\uc\ua}\OpL^\ub \CQA_\uc^{\alpha\beta}(\CQA_\ua \solu)(\partial_\beta \CQA_\ub \solu)) \\
\end{align*}
We can substitute this into the previous calculation, rearrange a derivative in the new expression, reindex, and use the symmetry of $\CurlyU^{\ua\ub}$  to conclude that 
\begin{align*}
\Svol\CurlyUfour^{\ua\ub\alpha\beta}(\partial_\alpha\CQA_\ua \solu)(\partial_\beta \CQA_\ub \solu)
={}&{}\Svol\CurlyU^{\uc\ua}\OpL^\ub\CQA_\ub^{\alpha\beta} (\partial_\alpha\CQA_\ua \solu)(\partial_\beta \CQA_\uc \solu) \\
{}&{}-\partial_\alpha(\Svol\CurlyU^{\uc\ua}\OpL^\ub\CQA_\ub^{\alpha\beta} (\CQA_\ua \solu)(\partial_\beta \CQA_\uc \solu))\\
{}&{}+\partial_\alpha(\Svol \CurlyU^{\uc\ua}\OpL^\ub \CQA_\uc^{\alpha\beta}(\CQA_\ua \solu)(\partial_\beta \CQA_\ub \solu)) \\
%={}&{}\Svol\CurlyU^{\uc\ua}\OpL^{\alpha\beta} (\partial_\alpha\CQA_\ua \solu)(\partial_\beta \CQA_\uc \solu) \\
%{}&{}-\partial_\alpha(\Svol\CurlyU^{\uc\ua}\OpL^{\alpha\beta} (\CQA_\ua \solu)(\partial_\beta \CQA_\uc \solu))+\partial_\alpha(\Svol \CurlyU^{\uc\ua}\OpL^\ub \CQA_\uc^{\alpha\beta}(\CQA_\ua \solu)(\partial_\beta \CQA_\ub \solu)) .
%\end{align*}
%Reindexing and using the symmetry of $\CurlyU^{\ua\ub}$, we have
%\begin{align*}
%\Svol\CurlyUfour^{\ua\ub\alpha\beta}(\partial_\alpha\CQA_\ua u)(\partial_\beta \CQA_\ub u)
={}&{}\Svol \CurlyU^{\ua\ub}\OpL^{\alpha\beta}(\partial_\alpha\CQA_\ua \solu)(\partial_\beta \CQA_\ub \solu)\\
{}&{}-\partial_\alpha\left(\Svol(\CurlyU^{\ua\ub}\OpL^{\alpha\beta}-\CurlyUfour^{\ua\ub\alpha\beta})(\CQA_\ua \solu)(\partial_\beta\CQA_\ub \solu)\right) .
\end{align*}
Applying the definition of $\MorBoundaryTermA{\pairM}$ gives the desired result. 
\end{proof}

%%%%%%%%%%%%%%%%%%%%%%%%%%%%%%%%%%%%%%%%%%%%%%%%%%%%%%%%%%%%%%%%%%%%%%%
\newpageForSubsection
\subsection{Choosing the weights}
\label{SS:ChoosingTheWeights}

In this section, we choose the weights $\fnMna$ and $\fnMnb$ to ensure
the positivity of the highest-order terms in the right-hand side of the estimate in the previous lemma, lemma \ref{Lemma:rearrangements}. 

\begin{definition}
Given a positive value for the parameter $\epsilonMorawetzTurnOndtSquared$, we use the following weights to define the Morawetz vector field, 
\begin{align*}
\fnMna={}&{} \fnMca \fnMda ,{}&{}
\fnMnb={}&{} \fnMcb \fnMdb ,\\
\fnMca={}&{} \frac{\KDelta}{(r^2+a^2)^2} ,{}&{} 
\fnMcb={}&{} \frac{(r^2+a^2)^4}{3r^2-a^2} , \\
\fnMda ={}&{} 1 -\epsilonMorawetzTurnOndtSquared\left(\frac{\KDelta}{(r^2+a^2)^2}\right) ,{}&{}
\fnMdb ={}&{} \frac{1}{2r} . 
\end{align*}
\end{definition}

\begin{remark}
The goal in choosing the various weight functions is to obtain
nonnegativity for $\CurlyA^{\ua\ub} (\dr\CQA_\ua \solu)(\dr\CQA_\ub
\solu)$ and $\CurlyU^{\ua\ub}\OpL^{\alpha\beta}
(\partial_\alpha\CQA_\ua \solu)(\partial_\beta\CQA_\ub \solu)$, with
possible degeneracy in the $\CurlyU$ term near $r=3M$. As explained in
the introduction, the guiding principle is to introduce operators
$\DiffCurlyRTilde=\DiffCurlyRTilde(r;M,a;\dt,\dphi,\OpQ)$ and
$\DDiffCurlyRTTilde=\DDiffCurlyRTTilde(r;M,a;\dt,\dphi,\OpQ)$ that are
analogues of the corresponding quantities for null geodesics,
$\DiffCurlyRTilde(r;M,a;\GeodesicEnergy,\GeodesicLz,\GeodesicQ)$ and
$\DDiffCurlyRTTilde(r;M,a;\GeodesicEnergy,\GeodesicLz,\GeodesicQ)$. In
particular, one expects that
$\DiffCurlyRTilde=\DiffCurlyRTilde(r;M,a;\dt,\dphi,\OpQ)$ should be
(weakly spacetime) elliptic when
$\DiffCurlyRTilde(r;M,a;\GeodesicEnergy,\GeodesicLz,\GeodesicQ)$ is
nonnegative, and similarly for $\DDiffCurlyRTTilde$. Analogously, we
think of $\CurlyA$ as the product of the positive quantities $\OpL$
and $-\DDiffCurlyRTTilde$, and we think of $\CurlyU$ as
$\DiffCurlyRTilde^2$.

Thus, the weight functions are chosen so that, for any null geodesic,
the quantity
$-\DDiffCurlyRTTilde(r;M,a;\GeodesicEnergy,\GeodesicLz,\GeodesicQ)$ is
nonnegative and the quantity
$\DiffCurlyRTilde(r;M,a;\GeodesicEnergy,\GeodesicLz,\GeodesicQ)$
can vanish only at a value of $r$ near $3M$. For $|a|\ll M$, the
orbiting null geodesics are near $r=3M$. On orbiting null geodesics,
$\DiffCurlyRTilde(r;M,a;\GeodesicEnergy,\GeodesicLz,\GeodesicQ)$
vanishes and
$-\DDiffCurlyRTTilde(r;M,a;\GeodesicEnergy,\GeodesicLz,\GeodesicQ)$ is
positive. The functions $\fnMna$ and $\fnMnb$ are chosen so that on
any null geodesic, the quantity
$\DiffCurlyRTilde(r;M,a;\GeodesicEnergy,\GeodesicLz,\GeodesicQ)$
vanishes only in a neighbourhood of $r=3M$ and the quantity
$-\DDiffCurlyRTTilde(r;M,a;\GeodesicEnergy,\GeodesicLz,\GeodesicQ)$ is
positive. 

We have chosen the weights so that the following properties hold:
\begin{enumerate}
\item The definition of $\DiffCurlyRTilde$ in equation
  \eqref{eq:DefDiffCurlyRTilde} is made so that $\CurlyU$ takes the form
  as $\fnMnb\DiffCurlyRTilde^2$ in lemma \ref{Lemma:BasicMorawetzDeformation}. 
\item $\epsilonMorawetzTurnOndtSquared^2$ is the coefficient of
  $\GeodesicEnergy^2$ in $\DiffCurlyRTilde(r;M,a;\GeodesicEnergy,\GeodesicLz,\GeodesicQ)$ and $\DDiffCurlyRTTilde(r;M,a;\GeodesicEnergy,\GeodesicLz,\GeodesicQ)$, and hence of $\dt^2$ in the operators
  $\DiffCurlyRTilde$ and $\DDiffCurlyRTTilde$. 
\item $\fnMca$ is such that, if $\fnMda$ had been equal to $1$, which
  corresponds to $\epsilonMorawetzTurnOndtSquared=0$, then the
  coefficient of $\GeodesicEnergy^2$ in $\DiffCurlyRTilde(r;M,a;\GeodesicEnergy,\GeodesicLz,\GeodesicQ)$, and hence
  of $\dt^2$ in $\DiffCurlyRTilde$, would be zero, as in equation
  \eqref{eq:GeodesicsDiffCurlyRTilde}. 
\item $\fnMda$ is such that, if $\epsilonMorawetzTurnOndtSquared>0$,
  then the coefficient of
  $\epsilonMorawetzTurnOndtSquared\GeodesicEnergy^2$ in $\DiffCurlyRTilde(r;M,a;\GeodesicEnergy,\GeodesicLz,\GeodesicQ)$, and hence of $\epsilonMorawetzTurnOndtSquared\dt^2$ in $\DiffCurlyRTilde$, is nonnegative and a perturbation (in $\epsilonMorawetzTurnOndtSquared$) of the coefficient of $\OpQ$.
\item $\fnMcb$ is such that, if $\fnMda$ and $\fnMdb$ had both been
  equal to $1$, then the coefficient of $\GeodesicEnergy\GeodesicLz$
  in
  $\DDiffCurlyRTTilde(r;M,a;\GeodesicEnergy,\GeodesicLz,\GeodesicQ)$,
  and hence of $\dt\dphi$ in $\DDiffCurlyRTTilde$, would vanish, as in
  equation \eqref{eq:GeodesicsDDiffCurlyRTTilde}.
\item $\fnMdb$ is such that 
\begin{enumerate}
\item  $\DDiffCurlyRTTilde(r;M,a;\GeodesicEnergy,\GeodesicLz,\GeodesicQ)$ is positive
  everywhere and 
\item $(\fnMna\fnMnb
  \DiffCurlyRTilde(r;M,a;\GeodesicEnergy,\GeodesicLz,\GeodesicQ))^2\gMetric(\dr,\dr)\lesssim (\GeodesicEnergy^2+\GeodesicLz^2+\GeodesicQ)^2
  \gMetric(\vecTBlend,\vecTBlend)$. \label{condition6.2}
\end{enumerate}
In particular, from the dominant energy condition, condition
\ref{condition6.2} allows us to show that
$\GenEnergy{\pairM}[\solu]\lesssim
\GenEnergyOrder{\vecTBlend}{3}[\solu]$. Once the form $\fnMdb=Cr^{-1}$ was
chosen, the factor of $C=1/2$ was chosen so that, when $a=0$ and $\epsilonMorawetzTurnOndtSquared=0$, the coefficient of $\GeodesicLz^2+\GeodesicQ$ in
$\CurlyA$ is equal to $1$.
\end{enumerate}

The factors $\DiffCurlyRTilde$, $\fnMca$, $\fnMcb$, and $\fnMda$ are uniquely defined by the above properties. In contrast, the factor $\fnMdb$ is both overdetermined, since we have chosen it to satisfy two conditions that are not a priori obviously compatible, and underdetermined, since it so happens that there are many functions that allow these two conditions to be satisfied. 
\end{remark}

The statement and proof of the following lemma make use of the norms given in subsection \ref{SS:Norms}

\begin{lemma}
\label{Lemma:Morawetz}
There are positive constants $\epsilonMorawetzSlowRotation$,
$\epsilonMorawetzTurnOndtSquaredUpperBound$, and
$\epsilonInverseErrorInMorawetzBeforeHardy$ such that if $|a|\leq
\epsilonMorawetzSlowRotation$ and
$0<\epsilonMorawetzTurnOndtSquared\leq\epsilonMorawetzTurnOndtSquaredUpperBound$
and $\solu$ is a solution to the wave equation  
$\gWave \solu = 0$ 
%\eqref{eq:KerrWave}, 
then 
\begin{align}
\KSigma\nabla_\alpha{}&{}\left(\left(\GenMomentum{\pairM}[\solu]^\alpha +\MorBoundaryTermA{\pairM}[\solu]^\alpha +\MorBoundaryTermB{\pairM}[\solu]^\alpha  \right) \right) \nonumber \\
\geq{}&{}
M\frac{\KDelta^2}{r^2(r^2+a^2)} \normPtwiseTnEpsilon{2}{\epsilonMorawetzTurnOndtSquared}{\dr \solu}^2
+\frac16 \frac{9Mr^2-46M^2r+54M^3}{r^4} \normPtwiseTnEpsilon{2}{\epsilonMorawetzTurnOndtSquared}{\solu}^2 \nonumber\\
{}&{}+\frac{1}{4r}\frac{(r^2+a^2)^4}{3r^2-a^2}\DiffCurlyRTilde^\ua\DiffCurlyRTilde^\ub\OpL^{\alpha\beta}(\partial_\alpha\CQA_\ua \solu)(\partial_\beta\CQA_\ub \solu)\nonumber\\
{}&{}-\epsilonInverseErrorInMorawetzBeforeHardy\frac{\KDelta^2}{r^2(r^2+a^2)} ((|a|+\epsilonMorawetzTurnOndtSquared^2)\normPtwiseTnHomo{2}{\dr \solu}^2 +\epsilonMorawetzTurnOndtSquared\normPtwiseTnEpsilon{2}{a^2}{\dr \solu}^2)\nonumber\\
{}&{}-\epsilonInverseErrorInMorawetzBeforeHardy\frac1{r^2} (|a|\normPtwiseTnHomo{2}{\solu}^2 +\epsilonMorawetzTurnOndtSquared\normPtwiseTnEpsilon{2}{a^2}{\solu}^2) ,
\label{eq:Morawetz}
\end{align}
where $\DiffCurlyRTilde$ is defined by equation \eqref{eq:DefDiffCurlyRTilde} satisfies 
\begin{align*}
\DiffCurlyRTilde
={}&{}-2(r-3M)r^{-4} \OpLEpsilon{\epsilonMorawetzTurnOndtSquared} \\
{}&{} +aO(r^{-4})\dphi\dt +a^2O(r^{-5})\OpQ +a^2 O(r^{-5})\dphi^2\\
{}&{}+\epsilonMorawetzTurnOndtSquared a^2 O(r^{-5}) \dt^2  +\epsilonMorawetzTurnOndtSquared O(r^{-5})\OpQ +\epsilonMorawetzTurnOndtSquared O(r^{-5})\dphi^2, 
\end{align*} 
and where the $\MorBoundaryTermB{\pairM}[\solu]^\alpha$ satisfy
\begin{align*}
|\KSigma\MorBoundaryTermB{\pairM}[\solu]^t|
\lesssim{}&{}  \frac{\KDelta^2}{r^2(r^2+a^2)}|\dr\dt \solu| 
\sum_\ua |\dr\CQA_\ua \solu| 
+\frac{1}{r^2}|\dt \solu| \sum_\ua|\CQA_\ua \solu|, \\
\MorBoundaryTermB{\pairM}[\solu]^r={}&{}0.  
\end{align*}
and the angular components of $\MorBoundaryTermB{\pairM}$ are smooth. 

\begin{remark}
In the applications of this lemma, it will always be the case that the regions of integration have boundary which are level sets of $t$ or $r$, so that the angular components of $\MorBoundaryTermB{\pairM}[\solu]$ do not contribute to the boundary integral. 
\end{remark}

\begin{proof}
From lemma \ref{Lemma:rearrangements}, there are three terms to control, the $\CurlyU$, $\CurlyA$, and $\CurlyV$ terms. 

\begin{multistep}
\newstep{The $\CurlyU$ term.}
The $\CurlyU$ term can be expanded using the definition in lemma \ref{Lemma:BasicMorawetzDeformation} as
\begin{align*}
\frac12\CurlyU^{\ua\ub}\OpL^{\alpha\beta}(\partial_\alpha\CQA_\ua \solu)(\partial_\beta\CQA_\ub \solu)
={}&{} \frac14 \fnMnb \DiffCurlyRTilde^\ua\DiffCurlyRTilde^\ub \OpL^{\alpha\beta}(\partial_\alpha\CQA_\ua \solu)(\partial_\beta\CQA_\ub \solu), 
\end{align*}
so it is sufficient to calculate $\DiffCurlyRTilde$. With our choice of $\fnMna$ and $\fnMnb$, this is
\begin{align*}
\DiffCurlyRTilde
={}&{} -\epsilonMorawetzTurnOndtSquared ( 2(r-3M) r^{-4} +a^2O(r^{-5}) )\dt^2 \\
{}&{}+aM O(r^{-4}) \dphi\dt\\
{}&{}- ( 2(r-3M)r^{-4} +a^2O(r^{-5}) +\epsilonMorawetzTurnOndtSquared O(r^{-5}) )\OpQ \\
{}&{}- ( 2(r-3M)r^{-4} +a^2O(r^{-5}) +\epsilonMorawetzTurnOndtSquared O(r^{-5}) )\dphi^2 . 
\end{align*}

\newstep{The $\CurlyA$ term.}
The $\CurlyA$ term is
\begin{align*}
\CurlyAone^\ua\OpL^\ub (\dr\CQA_\ua \solu)(\dr\CQA_\ub\solu)
={}&{} \frac{\KDelta^2}{r^2+a^2} (-\DDiffCurlyRTTilde^{\ua})\OpL^\ub (\dr\CQA_\ua \solu)(\dr\CQA_\ub\solu). 
\end{align*}
With our choices of $\fnMna$ and $\fnMnb$, we find: 
\begin{align*}
-\DDiffCurlyRTTilde
={}&{} M\epsilonMorawetzTurnOndtSquared ( r^{-2} +a^2O(r^{-3})+\epsilonMorawetzTurnOndtSquared O(r^{-3}) ) \dt^2 \\
{}&{}+aMO(r^{-2}) \dphi\dt \\
{}&{}+M( r^{-2} +a^2O(r^{-3}) +\epsilonMorawetzTurnOndtSquared O(r^{-3}) )\OpQ \\
{}&{}+M( r^{-2} +a^2O(r^{-3}) +\epsilonMorawetzTurnOndtSquared O(r^{-3}) )\dphi^2 .
\end{align*}

We are interested in this because the operator $-\DDiffCurlyRTTilde$ is very close to $\frac{M}{r^2}\OpLEpsilon{\epsilonMorawetzTurnOndtSquared}$ in the sense that
\begin{align*}
\left|\left((-\DDiffCurlyRTTilde)-\frac{M}{r^2}\OpLEpsilon{\epsilonMorawetzTurnOndtSquared}\right)\dr \solu\right|
={}&{} a^2 O(r^{-3}) |\VFEffSym^2 \dr \solu| \\
{}&{}+ a O(r^{-2}) |\VFEffSym^2 \dr \solu| \\
{}&{}+ \epsilonMorawetzTurnOndtSquared^2 O(r^{-3}) |\dt^2 \dr \solu| \\
{}&{}+ \epsilonMorawetzTurnOndtSquared O(r^{-3}) |\VFAngAlg^2 \dr \solu|
+\epsilonMorawetzTurnOndtSquared a^2 O(r^{-3}) |\sin^2\theta \dt^2 \dr \solu| .
\end{align*}
Thus, 
\begin{align*}
\left|(\OpL\dr\solu)\left((-\DDiffCurlyRTTilde)-\frac{M}{r^2}\OpLEpsilon{\epsilonMorawetzTurnOndtSquared}\right)\dr \solu\right|
{}&{}\leq (|a|+\epsilonMorawetzTurnOndtSquared^2)O(r^{-2})\normPtwiseTn{2}{\dr\solu}^2
+\epsilonMorawetzTurnOndtSquared O(r^{-2}) \normPtwiseTnEpsilon{2}{a^2}{\dr\solu}^2 .
\end{align*}

Since $\OpL$ and $\OpLEpsilon{\epsilonMorawetzTurnOndtSquared}$ commute with functions of $r$, we can apply lemma \ref{Lemma:LepsilonL} to $\dr\solu$, to get
\begin{align*}
(\OpLEpsilon{\epsilonMorawetzTurnOndtSquared}\dr\solu)(\OpL\dr\solu)
\geq{}&{} \normPtwiseTnEpsilon{2}{\epsilonMorawetzTurnOndtSquared}{\dr\solu}^2 \\
{}&{}+\frac{1+\epsilon}{\Svol}\partial_\alpha\left(\Svol\TimeAndAngularDerivativesTwoPtThreea{\dr\solu}^\alpha\right)
+a^2|\VFEffSym\dr\psi|^2 .
\end{align*}
The divergence terms consist only of time and angular derivatives
which are exactly those coming from lemma \ref{Lemma:LepsilonL}. Thus,
we may multiply the equation by $\KDelta^2r^{-2}(r^2+a^2)^{-1}$ and
move this factor inside the divergence term. The terms from the
angular derivatives are smooth, and the terms from the time derivative
are then of the form
\begin{align*}
M\frac{\KDelta^2}{r^2(r^2+a^2)} \dt( (\dt \dr\solu)(\SLap \dr\solu)  .
\end{align*}
Thus, we only need to control contributions from these terms when they appear as boundary terms on hypersurfaces of constant $t$. They are controlled by
\begin{align}
M\frac{\KDelta^2}{r^2(r^2+a^2)} |\dt \dr\solu| |\SLap \dr\solu)| 
\lesssim{}&{} M\frac{\KDelta^2}{r^2(r^2+a^2)} |\dt \dr\solu| \sum_\ua |\CQA_\ua\dr\solu| . 
\label{eq:MorawetzBoundaryErrorFromIBPA}
\end{align}

Thus, 
\begin{align*}
\CurlyA^{\ua\ub}(\CQA_\ua\dr \solu)(\CQA_\ub\dr \solu)
\geq{}&{} M\frac{\KDelta^2}{r^2(r^2+a^2)}\normPtwiseTnEpsilondtTurnedOn{2}{\dr \solu}^2 \\
{}&{}-\epsilonInverseErrorInMorawetzBeforeHardy\frac{\KDelta^2}{r^2(r^2+a^2)} (|a|\normPtwiseTnHomo{2}{\dr \solu}^2 +\epsilonMorawetzTurnOndtSquared\normPtwiseTnEpsilon{2}{a^2}{\dr \solu}^2 +\epsilonMorawetzTurnOndtSquared^2 \normPtwiseTnHomo{2}{\dr \solu}^2)\\
{}&{}+\frac{1}{\Svol}\partial_\alpha\left(\Svol\MorBoundaryTermBa{\pairM}[\solu]^\alpha\right)
\end{align*}
with $\MorBoundaryTermBa{\pairM}$ satisfying the properties given in
the statement of this lemma.

%%%%%%%%%%%%%%%%%%%%%%%%%%%%%%%%%%%%%%%%%%%%%%%%%%%%%%%%%%%%%%%%%%%%%%%%%%%
\newstep{The $\CurlyV$ term.}
By direct computation, the $\CurlyV$ term is given by
\begin{align*}
\CurlyV^\ua\CQA_\ua 
={}&{}\frac{1}{4}\dr\KDelta\dr\fnMone\CQA_\ua\\
={}&{}\left(\epsilonMorawetzTurnOndtSquared\frac16(9Mr^{-2}-46M^2r^{-3}+54M^3r^{-4} +aO(r^{-4}) )\right)\dt^2 \\
{}&{}+aO(r^{-4}) \dphi\dt\\
{}&{}+ \left( \frac16(9Mr^{-2}-46M^2r^{-3}+54M^3r^{-4} +aO(r^{-4}) +\epsilonMorawetzTurnOndtSquared O(r^{-4}) )\right)\OpQ\\
{}&{}+ \left( \frac16(9Mr^{-2}-46M^2r^{-3}+54M^3r^{-4} +aO(r^{-4}) +\epsilonMorawetzTurnOndtSquared O(r^{-4}) )\right)\dphi^2 , 
\\
={}&{}\frac16(9Mr^{-2} -46M^2r^{-3} +54M^3r^{-4})\OpLEpsilon{\epsilonMorawetzTurnOndtSquared} \\
{}&{}+|a|O(r^{-4})\SOSym +\epsilonMorawetzTurnOndtSquared O(r^{-4})\OpQ +\epsilonMorawetzTurnOndtSquared O(r^{-4})\dphi^2 ,
\end{align*}
where we have used $O(r^{-4})\SOSym$ to denote terms of the form $O(r^{-4})\CQA_\ua$ with $\CQA_\ua\in\SOSym$. 

Applying the estimate in lemma \ref{Lemma:LepsilonL}, we find 
\begin{align*}
\CurlyV^{\ua\ub}(\CQA_\ua\solu)(\CQA_\ub\solu)
={}&{} \frac16 (9Mr^{-2}-46M^2r^{-3}+54M^3r^{-4})(\OpLEpsilon{\epsilonMorawetzTurnOndtSquared} \solu) (\OpL \solu) \\
{}&{}+O(r^{-2}) (|a|  \normPtwiseTnHomo{2}{\solu}^2+\epsilonMorawetzTurnOndtSquared \normPtwiseTnEpsilon{2}{a^2}{\solu}^2 ) \\
\geq{}&{} \frac16 (9Mr^{-2}-46M^2r^{-3}+54M^3r^{-4}) \normPtwiseTnEpsilon{2}{\epsilonMorawetzTurnOndtSquared}{\solu}^2  \\
{}&{}+O(r^{-2}) (|a|  \normPtwiseTnHomo{2}{\solu}^2+\epsilonMorawetzTurnOndtSquared \normPtwiseTnEpsilon{2}{a^2}{\solu}^2 ) \\
{}&{}+O(r^{-2})\frac{1+\epsilonMorawetzTurnOndtSquared}{\Svol}\partial_\alpha\left(\Svol\TimeAndAngularDerivativesTwoPtThreea{\solu}^\alpha\right) .
\end{align*}
Again, the divergence terms come from the application of lemma
\ref{Lemma:LepsilonL}, so that there is no radial derivative, the terms from the angular derivatives are smooth, and the terms from the time derivative give a contribution of the form
\begin{align}
\frac{C}{r^2} |\dt\solu| |\SLap\solu| 
\leq{}&{}\frac{C}{r^2} |\dt\solu| \sum_\ua|\CQA_\ua\solu| .
\label{eq:MorawetzBoundaryErrorFromIBPV}
\end{align}

The time and angular derivative terms arising in this step and the previous one are combined into $\MorBoundaryTermB{\pairM}$ and are controlled by \eqref{eq:MorawetzBoundaryErrorFromIBPA}-\eqref{eq:MorawetzBoundaryErrorFromIBPV}. 
\end{multistep}
\end{proof}
\end{lemma}

%%%%%%%%%%%%%%%%%%%%%%%%%%%%%%%%%%%%%%%%%%%%%%%%%%%%%%%%%%%%%%%%%%%%%%%
\begin{lemma}[Controlling the boundary terms]
\label{Lemma:controlingTheBoundaryTerms}
If $\solu$ is sufficiently smooth, satisfies $\gWave\solu=0$, and has initial data which decays sufficiently rapidly at infinity, then
\begin{align*}
\left|\GenEnergy{\pairM}[\solu]\right| 
+\left|\int_{\hst{t}} \MorBoundaryTermA{\pairM}^\alpha \diNormal_\alpha\right|
+\left|\int_{\hst{t}} \MorBoundaryTermB{\pairM}^\alpha \diNormal_\alpha\right|
{}&{}\leq C\GenEnergy{\vecTBlend}[\solu],
\end{align*}
and
\begin{align*}
\lim_{r\rightarrow \rp} \GenMomentum{\pairM}[\solu]^r ={}&{} 0 \\
%\lim_{r\rightarrow \rp} \MorBoundaryTermA{\pairM}^r ={}&{} 0 \\
\lim_{r\rightarrow \infty} |\GenMomentum{\pairM}[\solu]^r| ={}&{} 0 \\
%\lim_{r\rightarrow \infty} \MorBoundaryTermA{\pairM}^r ={}&{} 0 . 
\end{align*}
Here, by ``sufficiently rapidly'', we mean that the $\vecTBlend$ energies of $\solu$, $\VFTrueSym\solu$, and $\SOSym\solu$ are convergent integrals and that $\solu$ goes to zero as $r\rightarrow\infty$.
\end{lemma} 
(For the decay hypothesis, it is sufficient that $\lim_{r\rightarrow\infty}\solu=0$, $\lim_{r\rightarrow\infty}r\dr\solu=0$, $\lim_{r\rightarrow\infty}r\dt\solu=0$, and the same estimates hold for $\VFTrueSym\solu$ and $\SOSym\solu$. The convergence of the energies implies that these limits are valid, at least along a subsequence.)
\begin{proof}
We begin by noting that, from the simple Hardy estimate 
\begin{align}
\int_0^\infty |\solu|^2 \di x \lesssim {}&{} \int_0^\infty x^2 |\dx \solu|^2 \di x ,
\label{eq:BasicOneDHardyWithxSquared}
\end{align}
with $x=r-\rp $, one finds $\int_{\hst{t}} |\solu|^2\diThree$ by $\GenEnergy{\vecTBlend}(t)$. We will refer to this as the basic Hardy estimate. We use $\SOSym \solu$ to denote a term which can be bounded in absolute value by $|\SOSym \solu|$.

By direct computation, 
\begin{align*}
\GenEnergy{\pairM}
{}&{}=-\int_{\hst{t}} \left(\GenMomentum{\pairM}\right)_\alpha\vecTperp^\alpha\frac{\KPi}{\KDelta} \diThree ,\\
\left|\GenEnergy{\pairM}\right|
{}&{}\leq C \int_{\hst{t}} \left(|\vecTperp\CQA_\ua\solu| |\vecMr| |\dr\CQA_\ub\solu|\frac{\KPi}{\KDelta} + |\vecTperp\CQA_\ua\solu| |\fnM| |\CQA_\ub\solu|\frac{\KPi}{\KDelta} \right)\diThree \\
{}&{}\leq C\int_{\hst{t}} \left(\frac{\KPi}{\KDelta}|\vecTperp\SOSym\solu|^2 +\frac{\KPi}{\Delta}|\vecMr|^2|\dr\SOSym\solu|^2 +\frac{\KPi}{\KDelta}|\fnM|^2|\SOSym\solu|^2 \right)\diThree .
\end{align*}
Since $\KPi/\KDelta$, $(r^2+a^2)^2/\KDelta$, and $r^4/\KDelta$ are all uniformly equivalent, since $\vecMr$ is bounded by a multiple of $\KDelta r^{-2}$, and since $|\fnM|$ is bounded by a multiple of $\KDelta r^{-3}$, it follows from estimate \eqref{eq:EnergyExpansion} that $|\GenEnergy{\pairM}|\leq C \GenEnergyOrder{\vecTBlend}{3}$. 

Using that $\diNormal_\alpha =-\KSigma \di t_\alpha \diThree$, and that each $\CurlyU^{\ua\ub}=O(r^{-1})$, it follows that
\begin{align*}
\left|\int_{\hst{t}} \MorBoundaryTermA{\pairM}^\alpha \diNormal_\alpha \right|
{}&{}\leq \int_{\hst{t}} |\CurlyU^{\ua\ub}\OpL^{\alpha t}-\CurlyUfour^{\ua\ub\alpha t}| |\CQA_\ua\partial_\alpha\solu| |\CQA_{\ub}\solu| \diThree \\
{}&{}\leq C \int_{\hst{t}} r^{-1} |\SOSym \solu| |\VFEffSym \SOSym \solu| \diThree \\
{}&{}\leq C \GenEnergyOrder{\vecTBlend}{3} .
\end{align*}
Similarly, since the $t$ component of the second boundary term was partially estimated in the statement of lemma \ref{Lemma:Morawetz}, it follows that
\begin{align*}
{}&{}\left|\int_{\hst{t}} \MorBoundaryTermB{\pairM}^\alpha \diNormal_\alpha\right|\\
{}&{}\leq C\int_{\hst{t}} \left(O((\KDelta/r^2)^2,1) |\dr\dt \solu| |\dr\SOSym \solu| 
+O(r^{-2}) |\dt \solu||\SOSym \solu| \right)\diThree \\
{}&{}\leq C\int_{\hst{t}} \left(O((\KDelta/r^2)^2,1) (|\dr\SOSym \solu|^2+ |\dr\dt \solu|) +O(r^{-2}) (|\SOSym\solu|^2 +|\dt \solu|^2)\right) \diThree \\
{}&{}\leq C\GenEnergyOrder {\vecTBlend}{3} .
\end{align*}

The limits at $\rp$ and $\infty$ are easily evaluated. The radial component of the momentum consists of bounded functions times a power of $\KDelta$, so they vanish at $\rp$. For $r$ large, from the calculation at the start of this proof, we know that $|\GenMomentum{\pairM}[\solu]^r|\lesssim O(r^2)|\partial_r\SOSym\solu|^2+|\SOSym\solu|^2$. Since solutions of the wave equation have finite speed of propagation, if the initial data falls off sufficiently rapidly, then so will the solution at any later time, hence $|\GenMomentum{\pairM}[\solu]^r|$ will tend to zero as $r\rightarrow\infty$. 
\end{proof}

Note that it is not necessary to estimate the limits of the radial components boundary terms, $\MorBoundaryTermA{\pairM}^r$ and $\MorBoundaryTermB{\pairM}^r$, since these components are identically zero. 

%%%%%%%%%%%%%%%%%%%%%%%%%%%%%%%%%%%%%%%%%%%%%%%%%%%%%%%%%%%%%%%%%%%%%%%%%
%%%%%%%%%%%%%%%%%%%%%%%%%%%%%%%%%%%%%%%%%%%%%%%%%%%%%%%%%%%%%%%%%%%%%%%%%
\newpageForSubsection
\subsection{The Hardy estimate}
In \eqref{eq:Morawetz}, the coefficient of
$\normPtwiseTnEpsilondtTurnedOn{2}{\solu}^2$ is positive except in a
compact range of $r$ values. The purpose of this section is to prove a
Hardy estimate which allows us to get a globally positive coefficient
for $\normPtwiseTnEpsilondtTurnedOn{2}{\solu}^2$ by using the
positivity of the term involving
$\normPtwiseTnEpsilondtTurnedOn{2}{\dr\solu}^2$. The proof is a bit
technical and can be omitted it on a first reading, since the proof is independent of the rest of the Morawetz estimate. 

%%%%%%%%%%%%%%%%%%%%%%%%%%%%%%%%%%%%%%%%%%%%%%%%%%%%%%%%%%%%%%%%%%%%%%%%%%
\begin{lemma}
There exist positive $\epsilonHardySlowRotation$ and $\epsilonHardyExtraFactorOnPotential$ such that 
if $|a|\leq\epsilonHardySlowRotation$, then 
for any smooth function $\fnHardyToEstimate$ on $[\rp,\infty)\times S^2$ which 
is bounded on $[\rp,\infty)$, 
\begin{align}
\int_{\rp}^\infty {}&{}\left( \frac{\KDelta^2}{r^2(r^2+a^2)} (\dr \fnHardyToEstimate)^2 
+  \frac16 \frac{9r^2-46Mr+54M^2}{r^4} \fnHardyToEstimate^2 \right) \di r \nonumber\\
\geq{}&{} \epsilonHardyExtraFactorOnPotential \int_{\rp}^\infty  \frac{\KDelta^2}{r^2(r^2+a^2)} (\dr \fnHardyToEstimate)^2 
+\frac{1}{r^2} \fnHardyToEstimate^2 \di r .
\label{eq:WeightedHardy}
\end{align}
\begin{proof}

The proof consists of several parts. The early parts of this proof follow the method of \cite{BlueSoffer:ODE}. First, we will demonstrate that it is sufficient to find a positive solution to an associated ODE (ordinary differential equation). Second, we rewrite the estimate and ODE in terms of a new function, $\fnHardyRedToEstimate$. Third, we will construct an explicit solution for the new ODE when $a=0$ and $\epsilonHardyExtraFactorOnPotential=0$. Fourth, we will argue that the construction of the explicit solution can be perturbed to cover nonzero $a$ and $\epsilonHardyExtraFactorOnPotential$, which will give a perturbed estimate for $\fnHardyRedToEstimate$. Fifth, we will show that this gives the estimate for the original function $\fnHardyToEstimate$. Finally, we will show that boundary conditions for the ODE do not place restrictions on the function $\fnHardyToEstimate$. 

\begin{multistep}
%%%%%%%%%%%%%%%%%%%%%%%%%%%%%%%%%%%%%%%%%%%%%%%%%%%%%%%%%%%%%%%%%%%%%%%%%%%%%
\newstep{Find a positive solution to the associated ODE.} 
We wish to show that if, for smooth, nonnegative $A$ and smooth $V$, the ODE
\begin{align*}
-\dr \PotlGeneraldrrCoeff \dr \fnHardyODESolu +\PotlGeneralLOCoeff \fnHardyODESolu = 0 ,
\end{align*}
has a smooth, positive solution $\fnHardyODESolu$ on $[r_0,\infty]$, 
then for any smooth function $\fnHardyToEstimate$ on $[r_0,\infty]$, there is the estimate
\begin{align}
\int_{r_0}^\infty \PotlGeneraldrrCoeff (\dr \fnHardyToEstimate)^2 + \PotlGeneralLOCoeff \fnHardyToEstimate^2 \di r \geq 0 ,
\label{eq:AbstractHardy}
\end{align}
as long as
\begin{align}
\fnHardyToEstimate^2 \PotlGeneraldrrCoeff \frac{\dr\fnHardyODESolu}{\fnHardyODESolu}
\label{eq:HardyBoundaryCondition}
\end{align}
vanishes at $r_0$ and $\infty$. Recall that a function is smooth on a closed interval if it is smooth on the interior and all derivatives have a limit at the boundary. 

Since $\fnHardyODESolu$ is positive, for any smooth $\fnHardyToEstimate$, we can define $\fnHardyODEAux=\fnHardyToEstimate/\fnHardyODESolu$. From integration by parts, 
\begin{align*}
\int_{r_0}^\infty \PotlGeneraldrrCoeff (\dr \fnHardyToEstimate)^2 + \PotlGeneralLOCoeff \fnHardyToEstimate^2 \di r
-[\PotlGeneraldrrCoeff\fnHardyODESolu \fnHardyODEAux(\dr(\fnHardyODESolu \fnHardyODEAux))]_{r_0}^\infty
%={}&{} \int_{r_0}^\infty 
%-\fnHardyODESolu \fnHardyODEAux (\dr \PotlGeneraldrrCoeff \dr \fnHardyODESolu \fnHardyODEAux)
%+\PotlGeneralLOCoeff \fnHardyODESolu^2\fnHardyODEAux^2
%dr \\
%={}&{} \int_{r_0}^\infty 
%-\fnHardyODESolu \fnHardyODEAux (\dr \PotlGeneraldrrCoeff ((\dr \fnHardyODESolu) \fnHardyODEAux +\fnHardyODESolu (\dr \fnHardyODEAux) ) )
%+\PotlGeneralLOCoeff \fnHardyODESolu^2\fnHardyODEAux^2
%dr \\
%={}&{} \int_{r_0}^\infty 
%-\fnHardyODESolu \fnHardyODEAux^2 (\dr \PotlGeneraldrrCoeff \dr \fnHardyODESolu)
%+\PotlGeneralLOCoeff \fnHardyODESolu^2\fnHardyODEAux^2
%dr \\
%{}&{}+\int_{r_0}^\infty 
%-2\fnHardyODESolu \fnHardyODEAux \PotlGeneraldrrCoeff (\dr \fnHardyODESolu) (\dr \fnHardyODEAux) 
%-\fnHardyODESolu^2 \fnHardyODEAux (\dr \PotlGeneraldrrCoeff) (\dr \fnHardyODEAux)
%-\fnHardyODESolu^2 \fnHardyODEAux \PotlGeneraldrrCoeff (\dr^2 \fnHardyODEAux) 
%dr \\
%={}&{} \int_{r_0}^\infty 
%-\fnHardyODESolu \fnHardyODEAux^2 (\dr \PotlGeneraldrrCoeff \dr \fnHardyODESolu)
%+\PotlGeneralLOCoeff \fnHardyODESolu^2\fnHardyODEAux^2
%dr \\
%{}&{}+\int_{r_0}^\infty 
%-(\dr \fnHardyODESolu^2) \fnHardyODEAux \PotlGeneraldrrCoeff (\dr \fnHardyODEAux) 
%-\fnHardyODESolu^2 \fnHardyODEAux (\dr \PotlGeneraldrrCoeff) (\dr \fnHardyODEAux)
%-\fnHardyODESolu^2 \fnHardyODEAux \PotlGeneraldrrCoeff (\dr^2 \fnHardyODEAux) 
%dr \\
={}&{} \int_{r_0}^\infty 
\fnHardyODESolu \fnHardyODEAux^2 (-\dr \PotlGeneraldrrCoeff \dr \fnHardyODESolu
+\PotlGeneralLOCoeff \fnHardyODESolu)
\di r \\
{}&{}+\int_{r_0}^\infty 
\fnHardyODESolu^2 \PotlGeneraldrrCoeff (\dr\fnHardyODEAux)^2 
\di r \\
{}&{}-[ \fnHardyODESolu^2\PotlGeneraldrrCoeff \fnHardyODEAux(\dr\fnHardyODEAux)]_{r_0}^\infty .
\end{align*}
Since $\fnHardyODESolu$ satisfies the ODE $-\dr \PotlGeneraldrrCoeff \dr \fnHardyODESolu +\PotlGeneralLOCoeff \fnHardyODESolu =0$, the first term on the right is zero. Cancelling the boundary terms on the right from those on the left leaves the estimate
\begin{align*}
\int_{r_0}^\infty \PotlGeneraldrrCoeff (\dr \fnHardyToEstimate)^2 + \PotlGeneralLOCoeff \fnHardyToEstimate^2 \di r
={}&{}\int_{r_0}^\infty 
\fnHardyODESolu^2 \PotlGeneraldrrCoeff (\dr\fnHardyODEAux)^2 
\di r 
+[ \fnHardyODEAux^2\PotlGeneraldrrCoeff \fnHardyODESolu(\dr\fnHardyODESolu)]_{r_0}^\infty .
\end{align*}
The boundary term vanishes under condition \eqref{eq:HardyBoundaryCondition}, and the integrand on the right is nonnegative, since $\fnHardyToEstimate=\fnHardyODEAux\fnHardyODESolu$. Therefore, 
\begin{align*}
\int_{r_0}^\infty \PotlGeneraldrrCoeff (\dr \fnHardyToEstimate)^2 + \PotlGeneralLOCoeff \fnHardyToEstimate^2 \di r
\geq{}&{} 0 .
\end{align*}

%%%%%%%%%%%%%%%%%%%%%%%%%%%%%%%%%%%%%%%%%%%%%%%%%%%%%%%%%%%%%%%%%%%%%%%%%%%
\newstep{Simplify the estimate to eliminate one of the coefficients.}
For the rest of this proof, we will take
\begin{align*}
\PotlGeneraldrrCoeff={}&{} \frac{\KDelta^2}{r^2(r^2+a^2)} .
\end{align*}
We will consider the function
\begin{align*}
\fnHardyRedToEstimate ={}&{} \PotlGeneraldrrCoeff^{1/2} \fnHardyToEstimate .
\end{align*}
Since $\PotlGeneraldrrCoeff^{1/2}$ is smooth on $[\rp,\infty)$ and vanishes linearly at $\rp $, the new function $\fnHardyRedToEstimate$ is also smooth and vanishes at least linearly at $\rp $. Its derivative satisfies
\begin{align*}
\dr \fnHardyToEstimate ={}&{} \frac{1}{\PotlGeneraldrrCoeff^{1/2}} (\dr \fnHardyRedToEstimate) -\frac12\frac{\dr\PotlGeneraldrrCoeff}{\PotlGeneraldrrCoeff^{3/2}} \fnHardyRedToEstimate . 
\end{align*}
Therefore, the left-hand side of \eqref{eq:AbstractHardy} is given by 
\begin{align*}
\int_{\rp}^\infty{}&{} (\dr\fnHardyRedToEstimate)^2 
-\frac{\dr\PotlGeneraldrrCoeff}{\PotlGeneraldrrCoeff}\fnHardyRedToEstimate(\dr\fnHardyRedToEstimate) 
+\left(\frac14\frac{(\dr\PotlGeneraldrrCoeff)^2}{\PotlGeneraldrrCoeff^2} +\frac{\PotlGeneralLOCoeff}{\PotlGeneraldrrCoeff}\right)\fnHardyRedToEstimate^2 \di r \\
={}&{}\int_{\rp}^\infty
(\dr\fnHardyRedToEstimate)^2 
+\left(\frac{\PotlGeneralLOCoeff}{\PotlGeneraldrrCoeff}+\frac12\frac{\dr^2\PotlGeneraldrrCoeff}{\PotlGeneraldrrCoeff}-\frac14\frac{(\dr\PotlGeneraldrrCoeff)^2}{\PotlGeneraldrrCoeff^2}\right)\fnHardyRedToEstimate^2 \di r 
-\left[ \frac12\frac{\dr\PotlGeneraldrrCoeff}{\PotlGeneraldrrCoeff}\fnHardyRedToEstimate^2\right]_{\rp}^\infty .
\end{align*}
If the original function $\fnHardyToEstimate$ is bounded, then the boundary term in this equality vanishes. The estimate that we shall prove in the subsequent steps of this proof is, for some $\epsilonHardyExtraFactorOnPotentialB>0$, 
\begin{align}
\int_{\rp}^\infty
(\dr\fnHardyRedToEstimate)^2 
+\left(\frac{\PotlGeneralLOCoeff}{\PotlGeneraldrrCoeff}+\frac12\frac{\dr^2\PotlGeneraldrrCoeff}{\PotlGeneraldrrCoeff}-\frac14\frac{(\dr\PotlGeneraldrrCoeff)^2}{\PotlGeneraldrrCoeff^2}\right)\fnHardyRedToEstimate^2 \di r 
\geq{}&{}\epsilonHardyExtraFactorOnPotentialB \int_{\rp}^\infty \frac{1}{\PotlGeneraldrrCoeff r^2} \fnHardyRedToEstimate^2 \di r 
\label{eq:HardyRed}
\end{align}
If $\fnHardyRedToEstimate$ satisfies this, then, by multiplying this estimate by $1-\epsilonHardyExtraFactorOnPotentialC$, we find
\begin{align*}
\int_{\rp}^\infty 
(\dr{}&{}\fnHardyRedToEstimate)^2
+\left(\frac{\PotlGeneralLOCoeff}{\PotlGeneraldrrCoeff}+\frac12\frac{\dr^2\PotlGeneraldrrCoeff}{\PotlGeneraldrrCoeff}-\frac14\frac{(\dr\PotlGeneraldrrCoeff)^2}{\PotlGeneraldrrCoeff^2}\right)\fnHardyRedToEstimate^2 \di r \\
\geq
\int_{\rp}^\infty {}&{}
\epsilonHardyExtraFactorOnPotentialC(\dr\fnHardyRedToEstimate)^2 \\
{}&{}+\left(\epsilonHardyExtraFactorOnPotentialC\left(\frac{\PotlGeneralLOCoeff}{\PotlGeneraldrrCoeff}+\frac12\frac{\dr^2\PotlGeneraldrrCoeff}{\PotlGeneraldrrCoeff}-\frac14\frac{(\dr\PotlGeneraldrrCoeff)^2}{\PotlGeneraldrrCoeff^2}\right)\fnHardyRedToEstimate^2 + (1-\epsilonHardyExtraFactorOnPotentialC)\epsilonHardyExtraFactorOnPotentialB \frac{1}{\PotlGeneraldrrCoeff r^2} \fnHardyRedToEstimate^2 \right) 
\di r .
\end{align*}
By taking $\epsilonHardyExtraFactorOnPotentialC>0$ sufficiently small and substituting back for $\fnHardyToEstimate$, we can conclude inequality \eqref{eq:WeightedHardy} holds.

%%%%%%%%%%%%%%%%%%%%%%%%%%%%%%%%%%%%%%%%%%%%%%%%%%%%%%%%%%%%%%%%%%%%%%%%%%%
\newstep{Construction of the explicit solution for $a=0$ and $\epsilonHardyExtraFactorOnPotential=0$.}
Following the arguments in the first section, we could prove the desired estimate (for $a=0$ and $\epsilonHardyExtraFactorOnPotential=0$) by finding a positive solution to 
\begin{align}
-\dr \PotlGeneraldrrCoeff \dr \fnHardyODESolu + \PotlGeneralLOCoeff \fnHardyODESolu = 0
\label{eq:OriginalUnperturbeODEInHardy}
\end{align}
with
\begin{align*}
\PotlGeneraldrrCoeff    ={}&{} \frac{(r^2-2Mr)^2}{r^4} ,\\
\PotlGeneralLOCoeff     ={}&{} \frac{1}{6}\frac{9r^2-46Mr+54M^2}{r^4}
\end{align*}
on the interval $[2M,\infty)$. However, by using the argument in the previous section, it is easier to use the transformed function
\begin{align}
\fnHardyRedODESolu={}&{} \PotlGeneraldrrCoeff^{1/2}\fnHardyODESolu = \left(\frac{(r^2-2Mr)^2}{r^4}\right)^{1/2}\fnHardyODESolu ,
\label{eq:HardyDefnRedODESolu}\\
x={}&{} r-2M , 
\label{eq:HardyDefnx}
\end{align}
and to solve the ODE \eqref{eq:OriginalUnperturbeODEInHardy}
\begin{align}
-\dx^2 \fnHardyRedODESolu + \PotlRedGeneralLOCoeff \fnHardyRedODESolu ={}&{}0 ,
\label{eq:HardyRedODE} \\
\PotlRedGeneralLOCoeff
={}&{}\frac{\PotlGeneralLOCoeff}{\PotlGeneraldrrCoeff} 
+\frac12\frac{\dx^2\PotlGeneraldrrCoeff}{\PotlGeneraldrrCoeff} 
-\frac14\frac{(\dx\PotlGeneraldrrCoeff)^2}{\PotlGeneraldrrCoeff^2} 
\nonumber\\
={}&{} \frac{9x^2-34Mx-2M^2}{6x^2(x+2M)^2}  
\label{eq:HardyDefnPotlRed}
\end{align}
on the interval $x\in[0,\infty)$. 

We first note the following properties of hypergeometric functions \cite{Abramowitz:SpecialFunctionHandbook,ErdelyiEtAl}. The hypergeometric function is typically written with parameters $F(a,b;c;z)$. This is also referred to as Gauss's hypergeometric function ${}_2F_1(a,b;c;z)$, but we will not use this notation. It should be clear in all cases whether $a$ refers to the first parameter of the hypergeometric function or to the angular momentum parameter of the Kerr spacetime. 
The hypergeometric function $F(a,b;c;z)$ has the following integral representation for $a<0<b<c$ and $z\not\in [1,\infty)$
\begin{align}
F(a,b;c;z)
={}&{} 
\frac{\Gamma(c)}{\Gamma(a)\Gamma(b)} \int_0^1 t^{b-1} (1-t)^{c-b-1} (1-tz)^{-a} \di t .
\label{eq:HardyHypergeometricIntegralRepresentation}
\end{align}
It is not obvious from this representation, but it is true, that $F$ is symmetric in its first two arguments, $F(a,b;c;z)=F(b,a;c;z)$. There are a vast number of further relations. The hypergeometric differential equation is
\begin{align}
z(1-z) \frac{d^2 w}{dz^2} +\left[c-(a+b+1)z\right]\frac{dw}{dz} -abw = 0 .
\end{align}
A pair of solutions to this equation is 
\begin{align*}
F(a,b;c;z) , \\
z^{1-c}F(a-c+1,b-c+1;2-c;z) .
\end{align*}

Returning to the ODE arising from the Hardy estimate, we introduce the parameters $\alpha$, $\beta$, and $d$ (to be chosen later), and the further substitution
\begin{align}
\fnHardyRedODESolu ={}&{} x^\alpha (x+d)^\beta \fnHardyRRedODESolu .
\label{eq:HardyDefnRRedODESolu}
\end{align}
The ODE now becomes 
\begin{align}
\fnHardyRedODESolu''={}&{}\left(\alpha(\alpha-1)x^{\alpha-2}(x+2)^\beta  +2\alpha\beta x^{\alpha-1}(x+d)^{\beta-1} +\beta(\beta-1)x^\alpha(x+d)^{\beta-2}  \right)\fnHardyRRedODESolu \nonumber\\
{}&{}+2\left(\alpha x^{\alpha-1}(x+d)^\beta+\beta x^\alpha(x+d)^{\beta-1}  \right)\fnHardyRRedODESolu' \nonumber\\
{}&{}+x^\alpha(x+d)^\beta \fnHardyRRedODESolu'' ,\nonumber\\
0={}&{}-\fnHardyRedODESolu''+\PotlRedGeneralLOCoeff \fnHardyRedODESolu \nonumber\\
={}&{}x^{\alpha-2}(x+d)^{\beta-2} P,
\label{eq:HardyRRedODE}\\
P={}&{}x^2(x+d)^2 \fnHardyRRedODESolu'' \nonumber\\
{}&{}-2x(x+d)((\alpha+\beta)x+\alpha d)\fnHardyRRedODESolu' \nonumber\\
{}&{}+\bigg(-(\alpha(\alpha-1)(x+d)^2 +2\alpha\beta x(x+d) +\beta(\beta-1)x^2) \nonumber\\
{}&{}\hspace{.5in}+\frac{9x^2 -34Mx-2M^2}{6 x^2 (x+2M)^2} x^2 (x+d)^2 \bigg) \fnHardyRRedODESolu  .
\label{eq:HardyTermP}
\end{align}
We conclude from \eqref{eq:HardyRRedODE}, that $P=0$. If we choose
\begin{align}
d={}&{} 2M , 
\label{eq:HardyDefnd}
\end{align}
then the rational function in the last term on the right reduces to a polynomial. 

The coefficient of $\fnHardyRRedODESolu''$ is $x^2(x+d)^2$, of $\fnHardyRRedODESolu'$ is $x(x+d)$ times a linear function, and of $\fnHardyRRedODESolu$ a quadratic. If we choose the parameters $\alpha$ and $\beta$ so that the coefficient of $\fnHardyRRedODESolu$ is a constant multiple of $x(x+d)$, then an over all factor of $x(x+d)$ can be dropped, leaving the coefficients of $\fnHardyRRedODESolu''$, $\fnHardyRRedODESolu'$, and $\fnHardyRRedODESolu$ as $x(x+d)$, a linear function, and a constant respectively. The substitution $z=-x/d$, then transforms the equation to the hypergeometric differential equation. Our goal is to show that such choices of $\alpha$ and $\beta$ can be made. 

It is now merely a matter of checking by direct calculation that this can be done. The coefficient of $\fnHardyRRedODESolu$ is
\begin{align}
{}&{}-\alpha(\alpha-1)(x^2+2xd+d^2) -2\alpha\beta(x^2+dx) -\beta(\beta-1)x^2 \nonumber\\
{}&{}+(3/2)x^2 -(17/3)Mx -M^2/3 . 
\label{eq:HardyEqnToSolveForAlphaBeta}
\end{align}
In this coefficient, we set the constant order term to zero
\begin{align*}
-\alpha(\alpha-1)d^2 -M^2/3 ={}&{}0 ,\\
\alpha
={}&{}\frac12 \pm \frac{\sqrt{6}}{6} .
\end{align*}
Fortunately, the term $\alpha\beta(x^2+dx)$ is already a multiple of
$x^2+dx$, so we may ignore it when trying to get the coefficient of
$\fnHardyRRedODESolu$ to be a multiple of $x^2+dx$. We set the ratios
of the remaining coefficients of $x^2$ and of $x$ in
\eqref{eq:HardyEqnToSolveForAlphaBeta} to be $d$, so that the
polynomial \eqref{eq:HardyEqnToSolveForAlphaBeta} becomes a multiple
of $x(x+d)$. This condition on the ratio yields
\begin{align*}
d\left((3/2)-\alpha(\alpha-1)-\beta(\beta-1)\right)
={}&{}-2d\alpha(\alpha-1) -(17/3)M .
\end{align*}
We can substitute $-\alpha(\alpha-1)=1/12$ to get
\begin{align*}
2\left((3/2)+(1/12)-\beta(\beta-1)\right)
={}&{}(1/3) -(17/3) , \\
%(19/6) -2\beta^2+2\beta ={}&{} -16/3 ,\\
%-\beta^2 +\beta +51/12={}&{} 0 \\
%-\beta^2 +\beta +17/4={}&{} 0,\\
\beta
%={}&{} (1/2) \pm \sqrt{1+17}/2 \\
={}&{} \frac12 \pm \frac{3\sqrt{2}}{2} .
\end{align*}
The four choices of sign provide four choices of simplified equations to study. For simplicity, we will consider only the equation arising from taking the $+$ sign in $\alpha$ and the $-$ sign in $\beta$.\footnote{This choice simplifies some expressions in the rest of this argument.}

We are left with the differential equation for $\fnHardyRRedODESolu$
\begin{align*}
x(x+2)\fnHardyRRedODESolu''
-2((1+\sqrt{6}/6 +3\sqrt{2}/2)x +1+\sqrt{6}/3) \fnHardyRRedODESolu' \\
+(19/6 -3\sqrt{2}/2 +\sqrt{6}/6 -\sqrt{3}) \fnHardyRRedODESolu{}&{}=0 ,
\end{align*}
Making the substitutions $z=-x/d$ and $\fnHardyRRRedODESolu(z)=\fnHardyRRedODESolu(x)$ gives
\begin{align}
z(1-z) \fnHardyRRRedODESolu''
+\left( (1+\sqrt{6}/3) -(2-3\sqrt{2}+\sqrt{6}/3)z \right) \fnHardyRRRedODESolu' \nonumber\\
+\left( -19/6 +3\sqrt{2}/2+\sqrt{3}-\sqrt{6}/6 \right) \fnHardyRRRedODESolu 
{}&{}=0 .
\label{eq:HardyRRRedODE}
\end{align}
Thus we have a hypergeometric differential equation, with solution $\fnHardyRRedODESolu=F(a,b;c,-x/d)$. We can immediately read off some quantities in terms of the hypergeometric parameters
\begin{align}
c={}&{} 1+\sqrt{6}/3 , 
\label{eq:HardyValueOfC}\\
-a-b-1={}&{} -2 +3\sqrt{2} -\sqrt{6}/3 ,\nonumber\\
-ab={}&{}-19/6 +3\sqrt{2}/2+\sqrt{3}-\sqrt{6}/6 . \nonumber
\end{align}
We can now solve for the remaining two parameters 
\begin{align}
\{a,b\}
%={}&{} \frac{(a+b)\pm\sqrt{(a+b)^2-4ab}}{2} \nonumber\\
%={}&{} \frac12 -\frac32\sqrt{2}+\frac{\sqrt{6}}{6} \pm\frac12\sqrt{(1-3\sqrt{2}+\sqrt{6}/3)^2-4(19/6 -3\sqrt{2}/2+\sqrt{3}-\sqrt{6}/6)} \nonumber\\
%={}&{} \frac12 -\frac32\sqrt{2}+\frac{\sqrt{6}}{6} \pm\frac12\sqrt{ 1-6\sqrt{2}+\frac23\sqrt{6}+18-2\sqrt12+\frac69 -\frac38/3 +6\sqrt{2}+4\sqrt{3}-\frac23\sqrt{6} } \nonumber\\
={}&{} \left\{\frac12 -\frac32\sqrt{2}+\frac{\sqrt{6}}{6} \pm\frac12\sqrt{7} \right\}. 
\label{eq:HardyValueOfAB}
\end{align}
We will make the choice $a<b$ so that 
\begin{align*}
a< -2.5<0<.1<b<.2<1.8< c .
\end{align*}
In particular
\begin{align*}
a < 0 < b < c .
\end{align*}
Thus, the integral representation \eqref{eq:HardyHypergeometricIntegralRepresentation} holds. Dividing by $\Gamma(c)/(\Gamma(a)\Gamma(b))$, we find that $\fnHardyRRRedODESolu(z)$ is positive when $z\leq0$. This means that $\fnHardyRRedODESolu$ is positive when $x\geq0$, $\fnHardyRedODESolu$ is also positive when $x\geq0$, and $\fnHardyODESolu$ is positive when $r>2M$. 

%%%%%%%%%%%%%%%%%%%%%%%%%%%%%%%%%%%%%%%%%%%%%%%%%%%%%%%%%%%%%%%%%%%%%%%%%%%
\newstep{The perturbed estimate for $\fnHardyRedODESolu$.}
In this step, we will prove that there are $0<\epsilonHardySlowRotationInPert$ and $0<\epsilonHardyExtraFactorOnPotentialInPert$ such that for $|a|\leq\epsilonHardySlowRotationInPert$ and all suitable $\fnHardyRedToEstimate$, 
\begin{align*}
\int_{0}^\infty |\dr\fnHardyRedToEstimate|^2 +\PotlRedGeneralLOCoeffPert \fnHardyRedToEstimate^2  \di x \geq 0, 
\end{align*}
for
\begin{align*}
\PotlRedGeneralLOCoeffPert
={}&{} \frac{9x^2 -34Mx-2M^2}{6x^2(x+d)^2} -\epsilonHardyExtraFactorOnPotentialInPert \frac{(M+x)^2}{x^2(x+d)^2} , \\
d={}&{} \rp -r_- \\
r_-={}&{}M-\sqrt{M^2-a^2}.
\end{align*}
This potential is of the form 
\begin{align}
\PotlRedGeneralLOCoeffPert ={}&{} \frac{C_1x^2 +C_2x +C_3}{C_4 x^2 (x+d)^2} , 
\label{eq:HardyDefnPotlRedPert}
\end{align}
with the coefficients $C_1, \ldots, C_4$, and $d$ perturbed from their original values in equation \eqref{eq:HardyDefnPotlRed}.

From the argument in step 1, it is sufficient to find a positive solution to the associated ODE \eqref{eq:HardyRedODE},
\begin{align*}
-\dx^2 \fnHardyRedODESolu + \PotlRedGeneralLOCoeffPert\fnHardyRedODESolu =0 ,
\end{align*}
with the perturbed potential $\PotlRedGeneralLOCoeffPert$. The analysis in step 3 found an explicit, positive solution for $x\in[0,\infty)$ for the parameter values dictated by the potential in equation \eqref{eq:HardyDefnPotlRed}. This step shows that the previous analysis also applies when the coefficients are perturbed. 

The previous analysis began by making the definition of $\fnHardyRRedODESolu$ in equation \eqref{eq:HardyDefnRRedODESolu}, in terms of the parameters $\alpha$ and $\beta$. The analysis then proceeded by choosing values for $\alpha$ and $\beta$ by solving quadratic equations coming from the coefficient in formula \eqref{eq:HardyEqnToSolveForAlphaBeta}, which lead to the new ODE \eqref{eq:HardyRRRedODE}. This ODE could be solved explicitly in terms of a hypergeometric function by solving linear and quadratic equations for the nonzero quantities $a$, $b$, and $c$. Since the coefficients in formula \eqref{eq:HardyEqnToSolveForAlphaBeta} depend continuously on the parameters $C_1$, $C_2$, $C_3$, $C_4$, and $d$ in the potential; since the coefficients in the ODE \eqref{eq:HardyRRRedODE} depend continuously on $\alpha$, $\beta$, and the coefficients in the potential; since all the quadratic equations involved had distinct, real roots; and since solutions to linear and quadratic equations depend continuously on the coefficients; it follows that positive solutions to the ODEs \eqref{eq:HardyRedODE} and \eqref{eq:HardyRRRedODE} can be found explicitly in terms of hypergeometric functions with parameters $a$, $b$, and $c$ depending continuously on the parameters in $\PotlRedGeneralLOCoeffPert$, at least when those parameter values are sufficiently close to the values given in equation \eqref{eq:HardyDefnPotlRed}. Similarly, when the perturbation of the parameter values in the potential $\PotlRedGeneralLOCoeffPert$ is sufficiently small, then the hypergeometric parameters maintain their order $a<0<b<c$. This gives the existence of positive $\epsilonHardySlowRotationInPert$ and $\epsilonHardyExtraFactorOnPotentialInPert$ which give the desired estimate for this step.

%%%%%%%%%%%%%%%%%%%%%%%%%%%%%%%%%%%%%%%%%%%%%%%%%%%%%%%%%%%%%%%%%%%%%%%%
\newstep{The perturbed estimate for the original function $\fnHardyToEstimate$.}
In the previous step, a particular type of perturbation of the potential was considered. In this step, we show that such perturbations are sufficient to control the type of perturbation appearing in our problem. 

From the argument in step 2, we wish to prove that there exist $0<\epsilonHardySlowRotation$ and $0<\epsilonHardyExtraFactorOnPotentialB$ such that for $0\leq|a|<\epsilonHardySlowRotation$ and suitable $\fnHardyRedToEstimate$ estimate \eqref{eq:HardyRed} holds, e.g.{}
\begin{align*}
\int_{\rp}^\infty
(\dr\fnHardyRedToEstimate)^2 
+\left(\frac{\PotlGeneralLOCoeff}{\PotlGeneraldrrCoeff}+\frac12\frac{\dr^2\PotlGeneraldrrCoeff}{\PotlGeneraldrrCoeff}-\frac14\frac{(\dr\PotlGeneraldrrCoeff)^2}{\PotlGeneraldrrCoeff^2}\right)\fnHardyRedToEstimate^2 \di r 
\geq{}&{}\epsilonHardyExtraFactorOnPotentialB \int_{\rp}^\infty  \frac{1}{\PotlGeneraldrrCoeff r^2} \fnHardyRedToEstimate^2 \di r ,
\end{align*}
with
\begin{align*}
\PotlGeneraldrrCoeff ={}&{}  \frac{\KDelta^2}{r^2(r^2+a^2)} ,\\
\PotlGeneralLOCoeff  ={}&{} \frac16 \frac{9r^2-46Mr+54M^2}{r^4}. 
\end{align*}

To simplify the following calculations, we introduce a new rotation parameter\footnote{This is typically denoted $r_-$.}
\begin{align*}
\aAlt={}&{} M-\sqrt{M^2-a^2} .
\end{align*}
When $\aAlt$ is treated as a function of $|a|$ with $M$ fixed, this is a continuous, increasing function on the interval $[0,M]$, which maps the interval $[0,M]$ to $[0,M]$. In addition, since the quantities which appear in our estimate (such as $\KDelta$ and $r^2+a^2$) only have a quadratic dependence on $a$, and since $a^2$ can be solved for as a quadratic expression in $\aAlt$, it follows that the quantities $\PotlGeneraldrrCoeff$ and $\PotlGeneralLOCoeff$ are rational functions in $(r,M,\aAlt)$. 

The new radial coordinate, analogous to the one defined in \eqref{eq:HardyDefnx}, is now defined to be
\begin{align*}
x={}&{} r - \rp  = r -(2M-2\aAlt) .
\end{align*}
Since $r$ can be expressed as a linear function of $(x,M,\aAlt)$, the quantities $\PotlGeneraldrrCoeff$ and $\PotlGeneralLOCoeff$ are rational functions in $(x,M,\aAlt)$. 

The quantity 
\begin{align*}
\PotlRedGeneralLOCoeff
={}&{} \frac{\PotlGeneralLOCoeff}{\PotlGeneraldrrCoeff} +\frac12\frac{\dr^2\PotlGeneraldrrCoeff}{\PotlGeneraldrrCoeff} -\frac14\frac{(\dr\PotlGeneraldrrCoeff)^2}{\PotlGeneraldrrCoeff^2}
\end{align*}
is rational in $(x,M,\aAlt)$; has degree, with respect to $x$, two lower in the numerator than in the denominator; has singularities in $x\in[-d,\infty)$ only at $x\in\{0,-d\}$ for fixed $M$ and $\aAlt$; these are of order at most two; and, for sufficiently small $\aAlt$, has no singularities in $\aAlt$ for fixed $x>0$ and $M$. Thus, we may expand it as
\begin{align*}
\PotlRedGeneralLOCoeff
={}&{} \frac{1}{\KDelta^2}\frac{P_0 +\aAlt P_>}{Q_0+\aAlt Q_>} ,
\end{align*}
where the functions $P_0$ and $Q_0$ are polynomials in $(x,M)$, the functions $P_>$ and $Q_>$ are polynomials in $(x,M,\aAlt)$, and $Q_0$ and $Q_0+\aAlt Q_>$ have no roots in $x\in[-d,\infty)$. Since $P_0/Q_0$ is determined explicitly by equation \eqref{eq:HardyDefnPotlRed}, it follows that 
\begin{align*}
\PotlRedGeneralLOCoeff - \frac{1}{\KDelta^2}\frac{P_0}{Q_0} 
={}&{} \frac{\aAlt}{\KDelta^2}\frac{P_>Q_0-P_0Q_>}{Q_0 (Q_0+\aAlt Q_>)}
\end{align*}
must decay like $r^{-2}$ as $r\rightarrow\infty$ for fixed $\aAlt$ and $M$ and has no singularities in $[-d,\infty)$ except for those coming from $\KDelta^{-2}$. Since this is a rational function, there is a constant $C$ such that
\begin{align*}
\left|\PotlRedGeneralLOCoeff - \frac{1}{\KDelta^2}\frac{P_0}{Q_0} \right|
\leq{}&{} \aAlt C\frac{(M+x)^2}{\KDelta^2} .
\end{align*}
Thus, there are sufficiently small $\epsilonHardySlowRotation$ and $\epsilonHardyExtraFactorOnPotentialB$ such that for $0\leq|a|<\epsilonHardySlowRotation$ 
\begin{align*}
\PotlRedGeneralLOCoeff -\epsilonHardyExtraFactorOnPotentialB \frac{1}{\PotlGeneraldrrCoeff r^2} > \PotlRedGeneralLOCoeffPert , 
\end{align*}
with $\PotlRedGeneralLOCoeffPert$ as in equation \eqref{eq:HardyDefnPotlRedPert}
The smallness of $\epsilonHardySlowRotation$ and $\epsilonHardyExtraFactorOnPotentialB$ is determined by the smallness of $\epsilonHardySlowRotationInPert$ and $\epsilonHardyExtraFactorOnPotentialInPert$. These then give $\epsilonHardySlowRotation$ and $\epsilonHardyExtraFactorOnPotential$ for which the desired estimate holds. 

%%%%%%%%%%%%%%%%%%%%%%%%%%%%%%%%%%%%%%%%%%%%%%%%%%%%%%%%%%%%%%%%%%%%%%%%%%
\newstep{Controlling the boundary terms.}
Since the argument from step 1 was applied to the function $\fnHardyRedToEstimate$, the boundary condition which must be imposed for this argument to hold is that 
\begin{align*}
\fnHardyRedToEstimate^2 \frac{\dr\fnHardyRedODESolu}{\fnHardyRedODESolu}
\end{align*}
vanishes at $\rp $ and at $\infty$. Since the positive solution to the ODE is given by
\begin{align*}
\fnHardyRedODESolu(r) 
={}&{} x^\alpha (x+d)^\beta \fnHardyRRedODESolu
=  x^\alpha (x+d)^\beta F(a,b;c;-z/d) \\
={}&{}  (r-\rp )^\alpha (r-r_-)^\beta F\left(a,b;c;-\frac{r-\rp}{\rp -r_-}\right) ,
\end{align*}
and the hypergeometric function is analytic (in its fourth argument) near zero, the ratio $\dr\fnHardyRedODESolu/\fnHardyRedODESolu$ will diverge at most inverse linearly at $r=\rp $. Thus, it is sufficient that $\fnHardyRedToEstimate$ vanish linearly at $r=\rp $. Since $\fnHardyRedToEstimate = \KDelta(r(r^2+a^2)^{1/2})^{-1} \fnHardyToEstimate$, it is sufficient that $\fnHardyToEstimate$ be smooth near $\rp$. 

To show the vanishing as $x\rightarrow\infty$, we first note that from
the form of the potential $\PotlRedGeneralLOCoeffPert$ in the ODE, the
solution $\fnHardyRedODESolu(r)$ will behave like a polynomial as
$r\rightarrow\infty$, so that
$\dr\fnHardyRedODESolu/\fnHardyRedODESolu$ will decay like a constant
times $1/r$. Thus, it is sufficient that $\fnHardyRedToEstimate$
remains bounded at infinity. Since $\fnHardyRedToEstimate = \KDelta(r(r^2+a^2)^{1/2})^{-1} \fnHardyToEstimate$, it is sufficient that $\fnHardyToEstimate$ be bounded near $\infty$. 

Thus, to obtain the vanishing of $\fnHardyRedToEstimate^2
(\dr\fnHardyRedODESolu)/\fnHardyRedODESolu$ at both $\rp$ and
$\infty$, it is sufficient that $\fnHardyToEstimate$ be smooth and
bounded on $[\rp,\infty)$.
\end{multistep}
\end{proof}
\end{lemma}

%%%%%%%%%%%%%%%%%%%%%%%%%%%%%%%%%%%%%%%%%%%%%%%%%%%%%%%%%%%%%%%%%%%%%%%%%%%%%%
\newpageForSubsection
\subsection{Integrating the Morawetz estimate}
\begin{lemma}
\label{Lemma:IntegratedMorawetzHomo}
There are positive constants $\epsilonMorawetzSlowRotation$, $\epsilonMorawetzPhotonWidth$, 
$\epsilonInvIntegratedMorawetzHomo$, and $\epsilonInvIntegratedMorawetzHomoB$ such that, for all
$|a|\leq\epsilonMorawetzSlowRotation$ and all smooth $\solu$ solving the
wave equation $\gWave \solu = 0$, the estimate
\begin{align}
\epsilonInvIntegratedMorawetzHomo (\GenEnergy{\vecTBlend}[\SOSym\solu](T_2){}&{} +\GenEnergy{\vecTBlend}[\VFTrueSym\solu](T_2) +\GenEnergy{\vecTBlend}[\SOSym\solu](T_1) +\GenEnergy{\vecTBlend}[\VFTrueSym\solu](T_1) )\nonumber\\
\geq \int_{T_1}^{T_2} \int_{\rp}^{\infty} \int_{S^2}
{}&{} \left(\frac{\KDelta^2}{r^4}\right) \normPtwiseTnHomo{2}{\dr \solu}^2 +r^{-2}\normPtwiseTnHomo{2}{\solu}^2  
+ \localiseAwayFromPhotonOrbits \frac{1}{r} \normPtwiseTnHomo{3}{\solu}^2 
\diFour \nonumber\\
-a^2\epsilonInvIntegratedMorawetzHomoB \int_{T_1}^{T_2} {}&{}\int_{\rp}^{\infty} \int_{S^2} \localiseAwayFromPhotonOrbits \frac{1}{r} \normPtwiseTn{2}{\solu}^2 \diFour 
, 
\label{eq:IntegratedMorawetzHomo}
\end{align}
holds, where $\localiseAwayFromPhotonOrbits$ is identically one for $|r-3M|>\epsilonMorawetzPhotonWidth$ and zero otherwise. 

\begin{proof}
We integrate the result of lemma \ref{Lemma:Morawetz} over the coordinate slab $(t,r,\theta,\phi)\in\slabST{T_1}{T_2}$, from which we get the integral of the right-hand side of estimate \eqref{eq:Morawetz}. From the Hardy estimate \eqref{eq:WeightedHardy}, the integral of the first two terms on the right-hand side of \eqref{eq:Morawetz} dominates an absolute constant times
\begin{align*}
\int_{T_1}^{T_2} \int_{S^2}\int_{\rp}^{\infty} \left( \frac{\KDelta^2}{r^2(r^2+a^2)} \normPtwiseTnEpsilon{2}{\epsilonMorawetzTurnOndtSquared}{\dr \solu}^2 + \frac{1}{r^2}\normPtwiseTnEpsilon{2}{\epsilonMorawetzTurnOndtSquared}{\solu}^2 \right) \diFour . 
\end{align*}
By taking $|a|$ sufficiently small relative to
$\epsilonMorawetzTurnOndtSquared$ and
$\epsilonMorawetzTurnOndtSquared$ sufficiently small relative to $1$,
these terms will also dominate the fourth and fifth terms, with a
constant factor left over. Since $\epsilonMorawetzTurnOndtSquared$ can
be chosen independently of $a$, the norms
$\normPtwiseTnEpsilon{2}{\epsilonMorawetzTurnOndtSquared}{\solu}$ can
be replaced by $\normPtwiseTnHomo{2}{\solu}$ at the price of a fixed
constant. The same is true for the norms of $\dr\solu$.

The only term which we still need to estimate is the third, 
\begin{align*}
\int_{T_1}^{T_2} \int_{\rp}^{\infty} \int_{S^2} \frac{(r^2+a^2)^4}{4r (3r^2-a^2)} \OpL^{\alpha\beta}(\partial_\alpha\DiffCurlyRTilde \solu)(\partial_\beta\DiffCurlyRTilde \solu)
\diFour . 
\end{align*}
The integrand can be estimated by
\begin{align*}
\frac{(r^2+a^2)^4}{4r (3r^2-a^2)} {}&{}\OpL^{\alpha\beta}(\partial_\alpha\DiffCurlyRTilde \solu)(\partial_\beta\DiffCurlyRTilde \solu) \\
\geq{}&{} \frac{(r^2+a^2)^4}{4r (3r^2-a^2)} |\VFEffSym \DiffCurlyRTilde\solu|^2 \\
\geq{}&{} \localiseAwayFromPhotonOrbits\frac{(r^2+a^2)^4}{4r (3r^2-a^2)} |\VFEffSym \DiffCurlyRTilde\solu|^2 \\
\geq{}&{} C \localiseAwayFromPhotonOrbits r^{-1} |\VFEffSym\OpLEpsilon{\epsilonMorawetzTurnOndtSquared}\solu|^2 \\
{}&{}+\localiseAwayFromPhotonOrbits O(r^{-5}) \epsilonMorawetzTurnOndtSquared
(|\VFEffSym\OpQ\solu|^2+|\VFEffSym\dphi^2\solu|^2+a^2|\VFEffSym\dt^2\solu|^2)\\
{}&{}+\localiseAwayFromPhotonOrbits O(r^{-3}) (|a|+\epsilonMorawetzTurnOndtSquared^2)|\VFEffSym\SOSym\solu|^2 .
\end{align*}
Recall $\VFEffSym$ is the set defined in section
\ref{SS:FurtherNotation} to consist of $\dt$ and the rotations around
the coordinate axes. To prove a lower bound on the first term in the
integrand, we first commute $\VFEffSym$ derivatives through
$\OpLEpsilon{\epsilonMorawetzTurnOndtSquared}$, and then apply
estimate \eqref{eq:ExplicitLepsilonLExpansionEqualityThirdOrder}:
\begin{align*}
|\VFEffSym\OpLEpsilon{\epsilonMorawetzTurnOndtSquared}\solu|^2
\gtrsim{}&{} |\OpLEpsilon{\epsilonMorawetzTurnOndtSquared}\VFEffSym\solu|^2
-a^2\normPtwiseTn{2}{\solu}^2 \\
\gtrsim{}&{} \normPtwiseTnEpsilondtTurnedOn{3}{\solu}^2 
- a^2 \normPtwiseTnEpsilon{3}{1}{\solu}^2
-a^2\normPtwiseTn{2}{\solu}^2 
-2\epsilonMorawetzTurnOndtSquared\frac{1}{\Svol}\partial_\alpha\left(\Svol\TimeAndAngularDerivativesTwoPtThreea{\VFEffSym\solu}^\alpha\right) .
\end{align*}

To estimate the remaining terms, we note that 
\begin{align*}
\epsilonMorawetzTurnOndtSquared(|\VFEffSym\OpQ\solu|^2{}&{}+|\VFEffSym\dphi^2\solu|^2+a^2|\VFEffSym\dt^2\solu|^2)
+(|a|+\epsilonMorawetzTurnOndtSquared^2)|\VFEffSym\SOSym\solu|^2\\
\lesssim{}&{} (\epsilonMorawetzTurnOndtSquared a^2 +|a|+\epsilonMorawetzTurnOndtSquared^2)|\dt^3\solu|^2 
+ (\epsilonMorawetzTurnOndtSquared a^2 +|a|+\epsilonMorawetzTurnOndtSquared^2)|\dt^2\dAng\solu|^2 \\
{}&{}+ (\epsilonMorawetzTurnOndtSquared +|a|+\epsilonMorawetzTurnOndtSquared^2)|\dt\SLap\solu|^2
+ (\epsilonMorawetzTurnOndtSquared +|a|+\epsilonMorawetzTurnOndtSquared^2)|\dAng\SLap\solu|^2 .
\end{align*}
These terms are dominated by
$\normPtwiseTnEpsilon{3}{\epsilonMorawetzTurnOndtSquared}{\solu}^2$ if
we again impose the conditions that $|a|$ is sufficiently small relative to
$\epsilonMorawetzTurnOndtSquared^2$ and $\epsilonMorawetzTurnOndtSquared$ is sufficiently small relative to
$1$. These smallness conditions are consistent with the one made in
the first paragraph of this proof. Thus,
\begin{align*}
\frac{(r^2+a^2)^4}{2r (3r^2-a^2)} \OpL^{\alpha\beta}(\partial_\alpha\DiffCurlyRTilde \solu)(\partial_\beta\DiffCurlyRTilde \solu)
\gtrsim{}&{} \localiseAwayFromPhotonOrbits r^{-1}\normPtwiseTnEpsilon{3}{\epsilonMorawetzTurnOndtSquared}{\solu}^2 
-Ca^2 \localiseAwayFromPhotonOrbits r^{-1}\normPtwiseTn{2}{\solu}^2\\
{}&{}+\localiseAwayFromPhotonOrbits
O(r^{-1})\frac{1}{\Svol}\partial_\alpha\left(\Svol\TimeAndAngularDerivativesTwoPtThreea{\VFEffSym\solu}\right) .
\end{align*}
Having chosen $\epsilonMorawetzTurnOndtSquared$, we can now make the
estimate
$\normPtwiseTnEpsilon{3}{\epsilonMorawetzTurnOndtSquared}{\solu}^2
\gtrsim\normPtwiseTnEpsilon{3}{1}{\solu}^2 $. 

The time derivative generated in this part of the argument is 
\begin{align*}
\dt (\localiseAwayFromPhotonOrbits O(r^{-1}) (\dt\VFEffSym\solu) (\SLap\VFEffSym\solu)) .
\end{align*}
Thus, the contribution of this time derivative on the boundary of the region of integration is bounded by $\GenMomentum{\vecTBlend}[\VFTrueSym\solu]^t +\GenMomentum{\vecTBlend}[\SOSym\solu]^t$. 

We must now control the integral of the momentum and the boundary terms over the boundary of the slab. All the angular derivative terms vanish, since $S^2$ has no boundary. Similarly, the boundary contributions along $r=\rp $ and $r\rightarrow\infty$ are zero by lemma \ref{Lemma:controlingTheBoundaryTerms}. (Geometrically, one would expect this, since $r=\rp $ is actually a two-dimensional surface, the bifurcation sphere, and not a three-dimensional hypersurface, so it should not contribute any boundary terms.)

We are left to estimate the integral of the momentum and the boundary terms over the hypersurfaces $t=T_1$ and $t=T_2$. From lemma \ref{Lemma:controlingTheBoundaryTerms}, these are estimated, at fixed $t$, by
\begin{align*}
\Big|\int_{\hst{t}} \left(\GenMomentum{\pairM}^\alpha + \MorBoundaryTermA{\pairM}^\alpha + \MorBoundaryTermB{\pairM}^\alpha\right) \diNormal_\alpha \Big| 
{}&{}\lesssim \GenEnergy{\vecTBlend}[\SOSym \solu](t) + \GenEnergy{\vecTBlend}[\dt \solu](t). 
\end{align*}
\end{proof}
\end{lemma} 

The previous lemma alone is insufficient, since it estimates only third derivatives, but the boundary terms involve both the second- and third-order energies. (Certain second-derivative terms are controlled, but these are not the important ones.) In the following lemma, we estimate the lower-order derivatives. 

%%%%%%%%%%%%%%%%%%%%%%%%%%%%%%%%%%%%%%%%%%%%%%%%%%%%%%%%%%%%%%%%%%%%%%%%%%%%%
\begin{lemma}
\label{Lemma:IntegratedMorawetz}
There are positive constants $\epsilonMorawetzSlowRotation$,
$\epsilonMorawetzTurnOndtSquared$, $\epsilonMorawetzPhotonWidth$, and 
$\epsilonInvIntegratedMorawetz$ such that for all
$|a|\leq\epsilonMorawetzSlowRotation$ and all smooth $\solu$ solving the
wave equation $\gWave \solu = 0$, the estimate
\begin{align}
\epsilonInvIntegratedMorawetz (\EnergyThree(T_2){}&{} +\EnergyThree(T_1) )
\label{eq:IntegratedMorawetz-inlemma} \\
\geq \int_{T_1}^{T_2} \int_{\rp}^{\infty} \int_{S^2}
{}&{} \bigg( \left(\frac{\KDelta^2}{r^4}\right) \normPtwiseTn{2}{\dr \solu}^2 +r^{-2}\normPtwiseTn{2}{\solu}^2  
+ \localiseAwayFromPhotonOrbits \frac{1}{r} \left(\normPtwiseTn{2}{\dt \solu}^2+\normPtwiseTn{2}{\dAng \solu}^2\right) \bigg)
\diFour . \nonumber
\end{align}
holds, where $\localiseAwayFromPhotonOrbits$ is identically one for $|r-3M|>\epsilonMorawetzPhotonWidth$ and zero otherwise. 
\end{lemma}
\begin{proof}
The Morawetz estimate, lemma \ref{Lemma:IntegratedMorawetzHomo},
controls the square integral of $\SOSym\solu$ and its first
derivatives. To prove the current lemma, it is sufficient to estimate
the corresponding integrals for $\solu$ and $\VFTrueSym\solu$. 

To treat $\solu$, we prove a Morawetz estimate using a classical,
first-order vector field. The estimate is valid for axi-symmetric
solutions; for nonaxial solutions, there are negative terms that can
be controlled using lemma \ref{Lemma:IntegratedMorawetzHomo}. 

In constructing this classical first-order vector field, we must find
scalar functions to play the roles of quantities previously
constructed from second-order symmetry operators. In particular, 
the role of $\CurlyRTilde$ is played by
$\KDelta$, and the role of $\DiffCurlyRTilde$ is played by a scalar
function $\fnMncClassical$. Since
$\DiffCurlyRTilde=\dr((\fnMna/\KDelta)\CurlyRTilde)$, this leads to
the slightly peculiar expression
$\fnMncClassical=\dr((\fnMna/\KDelta)\KDelta)$. Thus, the
quantities required for the proof of a Morawetz estimate are 
\begin{align*}
\fnMncClassical={}&{} \dr \left(\frac{\fnMna}{\KDelta} \KDelta \right) ,\\
\fnMpclassical ={}&{} \frac12\left(\dr\fnMna\right)\fnMnb\fnMncClassical ,\\
\vecMclassical ={}&{} \fnMna\fnMnb\fnMncClassical \dr ,\\
\fnMclassical  ={}&{} \frac12\left(\dr \vecMclassical^r\right) -\fnMpclassical .
\end{align*}
%Here $\fnMncClassical$ plays the role of $\DiffCurlyRTilde$. The role
%of $\CurlyRTilde$ is played by $\KDelta$, which gives the slightly
%peculiar expression $\dr((\fnMna/\KDelta)\KDelta)$. 
Using the same sort of calculations as before, we can obtain the analogue of \eqref{eq:DivPGenMorawetz}
\begin{align*}
\frac1\Svol\partial_\alpha \left(\Svol\GenMomentum{\pairClassical}^\alpha\right)
={}&{} \CurlyALzZero (\dr\solu)^2 + \CurlyULzZero^{\alpha\beta}(\partial_\alpha\solu)(\partial_\beta\solu) +\CurlyVLzZero|\solu|^2 , 
\end{align*}
with
\begin{align*}
\CurlyALzZero={}&{}\frac12 \fnMna^{1/2} \KDelta^{3/2} \dr\left(\fnMnb \frac{\fnMna^{1/2}}{\KDelta^{1/2}}\fnMncClassical \right) , \\
\CurlyU^{\alpha\beta}={}&{} \frac12 \fnMnb(\dr\fnMncClassical) \DiffCurlyRTilde^{\alpha\beta} ,\\
\CurlyVLzZero={}&{} \frac{1}{4}\dr\KDelta\dr\fnMna(\dr\fnMnb\fnMncClassical) .
\end{align*}
Taking the same choices of $\fnMna$ and $\fnMnb$ as before, we find
\begin{align*}
\CurlyALzZero={}&{} \frac12 \frac{\KDelta^2}{r^2+a^2}\left( \frac{1}{r^2} +|a|O(r^{-3}) +\epsilonMorawetzTurnOndtSquared O(r^{-2}) \right) ,\\
\CurlyU^{\alpha\beta}(\partial_\alpha\solu)(\partial_\beta\solu)
={}&{}\frac12 \fnMnb\fnMncClassical^2\OpQ^{\alpha\beta}(\partial_\alpha\solu)(\partial_\beta\solu) \\
{}&{}+ \frac12 \fnMnb\left(\fnMncClassical^2 +a^2O(r^{-8})\right)(\dphi\solu)^2\\
{}&{}+ a\fnMnb O(r^{-9})(\dt\solu)(\dphi\solu) \\
{}&{}+\frac12\fnMnb\epsilonMorawetzTurnOndtSquared\left(\dr\frac{\KDelta}{(r^2+a^2)^2}\right)^2\left(1-2\epsilonMorawetzTurnOndtSquared\frac{\KDelta}{(r^2+a^2)^2}\right) (\dt\solu)^2 \\
\CurlyVLzZero={}&{} \frac16 \frac{9Mr^2-46M^2r+54M^3}{r^4} +(a+\epsilonMorawetzTurnOndtSquared)O(r^{-4}) .
\end{align*}
From the Hardy estimate, \eqref{eq:WeightedHardy}, it follows that
\begin{align*}
\int_{\rp}^\infty \CurlyALzZero (\dr\solu)^2  +\CurlyVLzZero|\solu|^2 \di r
\gtrsim{}&{} \int_{\rp}^\infty \left(\frac{\KDelta^2}{r^2(r^2+a^2)}|\dr\solu|^2 +\frac{1}{r^2}|\solu|^2 \right) \di r. 
\end{align*}

We now analyse the
$\CurlyU^{\alpha\beta}(\partial_\alpha\solu)(\partial_\beta\solu)$
term. Let $\rPSLzZero$ denote the value of $r$ which maximises
$\KDelta/(r^2+a^2)^2$.  In
$\CurlyU^{\alpha\beta}(\partial_\alpha\solu)(\partial_\beta\solu)$,
the coefficients of $(\dtheta\solu)^2$ and $(\dt\solu)^2$ are
nonnegative and vanish only at $\rPSLzZero$. Since $\fnMncClassical$
decays as $r^{-3}$ and is strictly positive at $r=\rp$, it follows
that the coefficient of $(\dphi\solu)^2$ is positive except in a small
$r$ neighbourhood of $\rPSLzZero$. Similarly, outside a slightly
larger $r$ neighbourhood of $\rPSLzZero$, using the positivity of the
coefficients of $(\dphi\solu)^2$ and $(\dt\phi)^2$, the
$(\dphi\solu)(\dt\solu)$ term can be estimated by the Cauchy-Schwarz
inequality, because of the small parameter $a$ and the faster decay
rate. 

Thus, it is sufficient to estimate the integral of
$|a|((\dphi\solu)^2+(\dphi\solu)(\dt\solu))O(r^{-2})$. Although this
expression does not have a sign, we refer to it as the negative
contribution in this argument. Integrating over the spherical
coordinates, we have
\begin{align*}
\int_{S^2}a((\dphi\solu)^2+(\dphi\solu)(\dt\solu)) \diTwo
={}&{}-\int_{S^2}a((\dphi^2\solu)(\solu)+(\solu)(\dphi\dt\solu)) \diTwo \\
\left|\int_{S^2}a((\dphi\solu)^2+(\dphi\solu)(\dt\solu)) \diTwo\right|
\lesssim{}&{} |a|\int_{S^2} |\solu|^2 \diTwo +|a|\int_{S^2}|\SOSym\solu|^2 \diTwo . 
\end{align*}
The first term on the right can be estimated by the contribution from
$\CurlyALzZero (\dr\solu)^2  +\CurlyVLzZero(\solu)^2$. In the second
term, the integrand can be dominated by
$\normPtwiseTnHomo{2}{\solu}^2$. 
Thus, 
\begin{align}
\int_{T_1}^{T_2}{}&{}\int_{\rp}^\infty\int_{S^2} \frac1\Svol\partial_\alpha
\left(\Svol\GenMomentum{\pairClassical}^\alpha\right) \diFour \nonumber\\
\gtrsim{}&{}  \int_{T_1}^{T_2}\int_{\rp}^\infty\int_{S^2}
\frac{\KDelta^2}{r^2(r^2+a^2)}|\dr\solu|^2 +\frac{1}{r^2}|\solu|^2
 +\localiseAwayFromPhotonOrbits (|\dt\solu|^2 +|\dAng\solu|^2)
\diFour \nonumber\\
{}&{}-|a|\int_{T_1}^{T_2}\int_{\rp}^\infty\int_{S^2}
\frac{1}{r^2}\normPtwiseTnHomo{2}{\solu}^2 \diFour . 
\label{eq:IntegratedMorawetzClassical}
\end{align}

We now treat $\VFTrueSym\solu$, by applying the same argument using a
classical vector field. The only terms in
$\Svol^{-1}\partial_\alpha(\Svol\GenMomentum{\pairClassical}[\VFTrueSym\solu]^\alpha)$
that fail to be nonnegative are those we termed the negative
contribution in the previous paragraph. These can be estimated by
\begin{align*}
\left|\int_{S^2}a((\dphi\VFTrueSym\solu)^2+(\dphi\VFTrueSym\solu)(\dt\VFTrueSym\solu)) \diTwo\right|
\lesssim{}&{} |a|\int_{S^2} |\SOSym\solu|^2 \diTwo 
\lesssim |a|\int_{S^2} \normPtwiseTnHomo{2}{\solu}^2\diTwo.
\end{align*}
This can also be estimated by the second-order terms in lemma
\ref{eq:IntegratedMorawetzClassical}. Thus, the analogue of estimate
\eqref{eq:IntegratedMorawetzClassical} holds with $\solu$ replaced by
$\VFTrueSym\solu$ on the left and in all but the last term on the
right. The last term on the right remains the integral of
$|a|\normPtwiseTnHomo{2}{\solu}^2/r^2$.

We note that the sum of the homogeneous norms $\normPtwiseTnHomo{2}{\solu}^2
+|\VFTrueSym\solu|^2+|\solu|^2$ is uniformly equivalent to the
inhomogeneous norm $\normPtwiseTn{2}{\solu}^2$. The same is obviously
true with $\dr\solu$ or $\dt\solu$ replacing $\solu$. We also note
that
$\normPtwiseTnHomo{3}{\solu}^2+|\dAng\VFTrueSym\solu|^2+|\dAng\solu|^2$
dominates $\normPtwiseTn{2}{\dAng\solu}^2$. 

In analogy with the previous results in lemma \ref{Lemma:controlingTheBoundaryTerms}, there is a constant and an upper bound on $a$ such that
\begin{align*}
|\GenEnergy{\pairClassical}[\solu]| 
\lesssim{}&{} C\GenEnergy{\vecTBlend}[\solu] ,\\
|\GenEnergy{\pairClassical}[\VFTrueSym\solu]| 
\lesssim{}&{} C\GenEnergyOrder{\vecTBlend}{2}[\solu] .
\end{align*}
We can now sum the result of lemma
\ref{Lemma:IntegratedMorawetzHomo}, estimate
\eqref{eq:IntegratedMorawetzClassical}, and its analogue for
$\VFTrueSym\solu$, and use the smallness of
$\epsilonMorawetzSlowRotation\leq|a|$. From this, we obtain the desired
result.
\end{proof}

\newpageForSubsection
%%%%%%%%%%%%%%%%%%%%%%%%%%%%%%%%%%%%%%%%%%%%%%%%%%%%%%%%%%%%%%%%%%%%%%%%%%%%%%
\subsection{Closing the argument}
We are now able to show that the energy associated with $\vecTBlend$ is uniformly bounded by its value on the initial hypersurface. When $a=0$, the energy is conserved. When $a\not=0$, the energy is no longer conserved, but, in the following theorem, we show that the factor by which it can change vanishes linearly in $|a|$. 

\begin{theorem}
\label{Thm:UniformEnergyBound}
There are positive constants $\epsilonBoundedEnergySlowRotation$ and
$\epsilonInvBoundedEnergy$ such that if
$|a|\leq\epsilonBoundedEnergySlowRotation$ and $\solu$ is a solution to the wave
equation 
$\gWave \solu = 0$, 
%\eqref{eq:KerrWave}, 
then for all $t_2\geq t_1\geq 0$: 
\begin{align*}
\EnergyThree(t_2) 
\leq (1+\epsilonInvBoundedEnergy |a|) \EnergyThree(t_1) .
\end{align*}
\begin{proof}
By corollary \ref{Corollary:EnergyGrowthByDerivatives} 
\begin{align*}
\EnergyThree(t_2){}&{}-\EnergyThree(t_1) \\
\leq{}&{} |a| C \int_{\slabST{t_1}{t_2}} \fnAtBlendLocation \left(\normPtwiseTn{2}{\dr \solu}^2 +\normPtwiseTn{3}{\solu}^2\right) \diFour .
\end{align*}
By the Morawetz estimate, lemma \ref{Lemma:IntegratedMorawetz}, for sufficiently small $a$, there is a constant $\epsilonInvBoundedEnergyInner$ such that the integral of the third derivatives is controlled by the energies. Thus, 
\begin{align*}
\EnergyThree(t_2) -\EnergyThree(t_1)
\leq{}&{} |a|\epsilonInvBoundedEnergyInner\left(\EnergyThree(t_2)+\EnergyThree(t_1) \right) . 
\end{align*} 
Thus, for $a$ sufficiently small (by which we mean
$|a|<\epsilonBoundedEnergySlowRotation$, with
$\epsilonBoundedEnergySlowRotation$ defined to be the minimum of the bound on $a$ arising from lemma \ref{Lemma:IntegratedMorawetz} and of the inverse of $\epsilonInvBoundedEnergyInner$) 
\begin{align*}
(1-|a|\epsilonInvBoundedEnergyInner)\EnergyThree(t_2)
\leq{}&{} (1+|a|\epsilonInvBoundedEnergyInner)\EnergyThree(t_1) , \\
\EnergyThree(t_2) \leq{}&{} \frac{1+|a|\epsilonInvBoundedEnergyInner}{1-|a|\epsilonInvBoundedEnergyInner}\EnergyThree(t_1) .
\end{align*}
Since, for sufficiently small $|a|$, the rational function $(1+|a|\epsilonInvBoundedEnergyInner)/(1-|a|\epsilonInvBoundedEnergyInner)$ is bounded above by $1+\epsilonInvBoundedEnergy |a|$ for some $\epsilonInvBoundedEnergy$, the desired result holds. 
\end{proof}
\end{theorem}

Finally, we note that since $\vecTBlend$ and the symmetry operators
are all $t$-translation invariant, the same is true for the set of
quadratic forms they define on each hypersurface of constant $t$,
$\GenEnergyOrder{\vecTBlend}{3}$.

%%%%%%%%%%%%%%%%%%%%%%%%%%%%%%%%%%%%%%%%%%%%%%%%%%%%%%%%%%%%%%%%%
\appendix
\section{Nondegenerate estimates using the Dafermos-Rodnianski
  red-shift vector field}
\label{S:DRVector}
The estimate in theorem \ref{Thm:IntroMorawetz} and the energy,
$\GenEnergy{\vecTBlend}$, which is bounded in theorem
\ref{Thm:UniformEnergyBound} are degenerate in the sense that the
integrands contain terms which vanish as $r\rightarrow\rp$. 
In this section, the degeneracy in the energy and decay estimates are
removed\footnote{We thank one of the referees for suggesting the removal of this degeneracy.} through an application of the red-shift vector field, which was
first used in this context in \cite{DafermosRodnianski:RedShiftSchwarzschild}.

In this section, Greek indices refer to the Kerr coordinates (called Kerr-star coordinates in section 2.5 of \cite{ONeill}) $(\tk,\rk,\hk,\pk)$ given by
\begin{align*}
\tk{}&{}=t+T(r), {}&{}
\rk{}&{}=r, {}&{}
\hk{}&{}=\theta, {}&{}
\pk{}&{}=\phi+A(r), 
\end{align*}
where
\begin{align*}
T(r){}&{}=\int_{3M}^r \frac{x^2+a^2}{x^2-2Mx+a^2} \di x, {}&{}
A(r){}&{}=\int_{3M}^r \frac{a}{x^2-2Mx+a^2} \di x .
\end{align*}
These coordinates are adapted to the future event horizon, but a
similar construction can be made to work in a neighbourhood of the
past horizon. In these coordinates,
\begin{align*}
\dtk{}&{}= \dt, {}&{}
\drk{}&{}= \dr +\frac{1}{\KDelta}\left( (r^2+a^2)\dt +a\dphi\right), {}&{}
\dhk{}&{}= \dtheta, {}&{}
\dpk{}&{}= \dphi ,
\end{align*}
and the metric takes the form
\begin{align*}
g_{tt}\di\tk^2 +2g_{t\phi}\di\tk\di\pk+g_{\phi\phi}\di\pk^2
+\KSigma\di\hk^2
+2\di\tk\di\rk-2a\sin^2\theta\di\pk\di\rk ,
\end{align*}
where $g_{tt}$, etc., refers to the metric in $(t,r,\theta,\phi)$
coordinates. In the $(\tk,\rk,\hk\,\pk)$ coordinates, the metric is
not singular at $r=\rp$. Because $g_{\rk\rk}=0$, the vector field
$\drk$ is null. In these coordinates, the volume form is 
\begin{align*}
\KSigma\sin\theta \di\tk\di\rk\di\hk\di\pk ,
\end{align*}
and the wave equation
becomes\footnote{After multiplying by $\KSigma$, as we have done throughout this paper.}
\begin{align}
0
{}&{}=\drk\KDelta\drk \solu \nonumber\\
{}&{}\quad+ \dtk(r^2+a^2)\drk \solu+\drk(r^2+a^2)\dtk \solu +2a\drk\dpk \solu \nonumber\\
{}&{}\quad+ \frac{1}{\sin\theta}\dhk\sin\theta\dhk \solu
+\frac{1}{\sin^2\theta}\dpk^2 \solu +2a\dpk\dtk \solu +a^2\sin^2\theta \dtk^2
\solu . 
\label{eq:KSWaveEqn}
\end{align}
It is convenient to introduce
\begin{align*}
\opQT^{\alpha\beta}(\partial_\alpha\solu)(\partial_\beta\solu)
=(\dhk\solu)^2+\frac1{\sin^2\theta}
(\dpk\solu)^2+2a(\dpk\solu)(\dtk\solu)+a^2\sin^2\theta(\dtk\solu)^2.
\end{align*}  
The last three terms factor as $((\sin\theta)^{-1}\dpk\solu+a\sin\theta\dtk\solu)^2$, so $\opQT^{\alpha\beta}(\partial_\alpha\solu)(\partial_\beta\solu)\geq0$. 
The associated operator
$\Svol^{-1}\partial_\alpha\Svol\opQT^{\alpha\beta}\partial_\beta$ is a
linear combination of our previous hidden symmetries. The contravariant metric, after rescaling by $\KSigma$, is 
%The Lagrangian for the wave equation is given by
\begin{align*}
%\KSigma\Lagrangian
%{}&{}=\KSigma\gMetric^{\alpha\beta}(\nabla_\alpha\solu)(\nabla_\beta\solu) \\
%{}&{}=\KDelta(\drk\solu)^2 \\
%{}&{}\quad+2(r^2+a^2)(\drk\solu)(\vecTk\solu)+2
%%(\drk\solu)
%a\left(1-\frac{r^2+a^2}{\rp^2+a^2}\right)(\drk\solu)(\dpk\solu)\\
%{}&{}\quad+\opQT^{\alpha\beta}(\partial_\alpha\solu)(\partial_\beta\solu) ,
\KSigma\gMetric^{\alpha\beta}
{}&{}=\KDelta\drk^\alpha\drk^\beta 
+2(r^2+a^2)\drk^{(\alpha}\dtk^{\beta)}+2a\drk^{(\alpha}\dpk^{\beta)}+\opQT^{\alpha\beta}\\
{}&{}=\KDelta\drk^\alpha\drk^\beta 
+2(r^2+a^2)\drk^{(\alpha}\vecTk^{\beta)}
+2 a\left(1-\frac{r^2+a^2}{\rp^2+a^2}\right)\drk^{(\alpha}\dpk^{\beta)}
+\opQT^{\alpha\beta}, 
\end{align*}
where $\vecTk=\dtk+\OmegaH\dpk$. Thus rescaling by $\KSigma$ provides the same simplifications in the Kerr-star coordinates as those described in subsection \ref{SS:SigmaRescaling} in the Boyer-Lindquist coordinates. The energy-momentum tensor, momentum density, and energy on a hypersurface are all covariant quantities, so they can be expressed in the $(\tk,\rk\,\hk,\pk)$ coordinate system. 

\begin{figure}
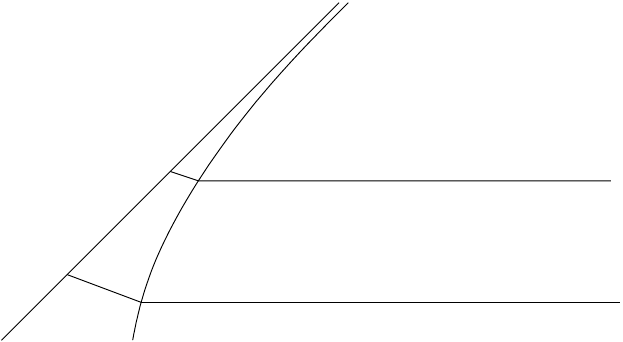
\caption{The hypersurfaces $\hstk{t_i}$, the region $\Omegak{[t_1,t_2]}$, and the hypersurface $\hsrk{[t_1,t_2]}$.}
\label{fig:hstk}
\end{figure}

In a neighbourhood of the horizon, it is convenient to work with
surfaces of constant $\rk-\tk$. The hypersurfaces and regions defined
in this paragraph are illustrated in figure \ref{fig:hstk}. For
$|a|\leq\epsilonInverseAppendix$, let $\epsilonNH$ be a small multiple
of $M$ to be determined later in the argument and let the near-horizon
radius be 
$\rNH=\rp+\epsilonNH$. Define the hypersurfaces $\hstk{\tau}$ as the
union of the hypersurface $\{(\tk,\rk,\hk,\pk):\rk\in[\rp,\rNH],
\rk-\tk=\rNH-\tau-T(\rNH)\}$ (in $(\tk,\rk,\hk,\pk)$ coordinates) with
the hypersurface $\{(t,r,\theta,\phi):r\geq\rNH,t=\tau\}$ (in
$(t,r,\theta,\phi)$ coordinates). This family of hypersurfaces is
continuous and is smooth except at $r=\rk=\rNH$. Define
$\hsrk{[t_1,t_2]}=\{(\tk,\rp,\hk,\pk):\tk\in[\rp-\rNH+t_1-T(\rNH),\rp-\rNH+t_2-T(\rNH)]\}$
and $\Omegak{[t_1,t_2]}$ to be the union of the region
$\{(\tk,\rk,\hk,\pk):\rk\in[\rp,\rNH],\tk\in[\rk-\rNH+t_1-T(\rNH),\rk-\rNH+t_2-T(\rNH)]\}$
with the region $\{(t,r,\theta,\phi):r>\rNH, t\in[t_1,t_2]\}$. Note
that the boundary of $\Omegak{[t_1,t_2]}$ is
$\hstk{t_1}\cup\hstk{t_2}\cup\hsrk{[t_1,t_2]}$.

On the portion of $\hstk{t}$ where $\tk-\rk$ is constant, the
future-directed normal volume form is
\begin{align*}
\diNormal^\alpha
{}&{}=(\KSigma\gMetric^{\tk\alpha}-\KSigma\gMetric^{\rk\alpha})\sin\theta\di\rk\di\hk\di\pk \\
{}&{}=(-\KSigma\dtk^\alpha+2Mr\drk^\alpha)\sin\theta\di\rk\di\hk\di\pk .
\end{align*}
Certain integrals are simplified by noting that for any vector field $\vecX$, 
\begin{align*}
\vecX_\alpha\diNormal^\alpha 
{}&{}=-\vecX_{\alpha}\KSigma\sin\theta(\di\tk^\alpha\di\rk\di\hk\di\pk -\di\rk^\alpha\di\tk\di\hk\di\di\pk) \\
{}&{}=(\vecX^{\tk}-\vecX^{\rk})\KSigma\sin\theta\di\rk\di\hk\di\pk .
\end{align*}
Similarly, on $\hsrk{[t_1,t_2]}$, one finds
\begin{align*}
\diNormal^\alpha{}&{}=\gMetric^{\rk\alpha}\KSigma\sin\theta\di\tk\di\hk\di\pk, \\
\vecX_\alpha\diNormal^\alpha{}&{}=\vecX^{\rk}\KSigma\sin\theta\di\tk\di\hk\di\pk .
\end{align*}

Consider the vector field $\vecY=\chiNH\left(\h\vecTk+\f\drk\right)$
with $\chiNH=\chiNH(r)$ identically $1$ for $r<\rNH=\rp+\epsilonNH$,
decreasing smoothly for $r\in[\rp+\epsilonNH,\rp+2\epsilonNH]$, and
identically zero for $r>\rp+2\epsilonNH$ and with $\h=\h(r)=\h(\rk)$
and $\f=\f(r)=\f(\rk)$ smooth and satisfying, for
$r\in[\rp-\epsilonNH,\rp+\epsilonNH]$, the five conditions (i) $\f<0$,
(ii) $\f'<0$, (iii) $\h>0$, (iv) $\h'>0$, and (v) $\h'>8|\f|/(r-M)$. In
particular, one can choose $\f$ and $\h$ to be linear if one chooses
the values of $\epsilonInverseAppendix$, $\f'(\rp)$, $\f(\rp)$,
$\h'(\rp)$, and $\h(\rp)$ in this order.

We now estimate the energy on $\hstk{0}$ generated by $\vecY$ in terms of the energy on $\hst{0}$ generated by the normal to $\hst{0}$. Let $\normalhst{0}$ be the normal to the hypersurface $\hst{0}$, where $t=0$, $r\geq\rp$. Let $\{\ONBasis_i\}_{i=0}^3$ denote an orthonormal basis at each point on $\hst{0}$ such that $X_0=\normalhst{0}$. The energy $\GenEnergy{\normalhst{0}}(\hst{0})$ is equivalent to the integral of $\sum_{i=0}^3 r^2|\ONBasis_i \solu|^2$. By a standard Hardy estimate, this means it dominates the integral of $|\solu|^2$. In a coordinate system which covers the bifurcation sphere (the limit $r\rightarrow\rp$ with $t=0$) and also covers $\hstk{0}\cap\{r<\rp+2\epsilonNH\}$, because both $\hst{0}$ and $\hstk{0}$ have a timelike normal and $\hstk{0}$ is in the causal future of $\hst{0}$, it follows that the local $H^1$ norm squared on $\hstk{0}$ is bounded by a multiple of the $H^1$ norm squared on $\hst{0}$. Thus, 
\begin{align*}
\GenEnergy{\vecY}[\solu](\hstk{0})\lesssim\GenEnergy{\normalhst{0}}[\solu](\hst{0}).
\end{align*} 
In the terminology of \cite{DafermosRodnianski:RedShiftSchwarzschild}, this is a Cauchy stability argument. Similarly, by the same type of argument the $\GenEnergy{\normalhst{0},3}(\hst{0})$ controls the $L^\infty$ norm in a neighbourhood of the bifurcation sphere and the integral of the derivatives in the spacetime region between $\hst{0}$ and $\hstk{0}$. 

At this stage in the argument, we assume several positivity
conditions. Later, these conditions are shown to hold. First, assume
that the energies defined by $\vecY$ on $\hstk{0}$, $\hstk{T}$, and
$\hsrk{[0,T]}$ are positive. Further assume that, for $r<\rNH$, the
divergence $\nabla_\alpha\GenMomentum{\vecY}^\alpha$ is negative
and that the modulus of the divergence dominates the square integral
of all $(\tk,\rk,\hk,\pk)$ partial derivatives of $\solu$. Recall that for
$\rNH=\rp+\epsilonNH\leq r\leq\rp+2\epsilonNH$, the square integral of
all partial derivatives can be estimated using the Morawetz estimate,
theorem \ref{Thm:IntroMorawetz}. Without loss of generality, we may assume $\rp+2\epsilonNH$ is smaller than $3M-\epsilonMorawetzPhotonWidth$, with $\epsilonMorawetzPhotonWidth$ from lemma \ref{Lemma:IntegratedMorawetz}. The same results apply to
$\GenMomentum{\vecY}[\CQA_\ua\solu]$. Thus,
\begin{align}
\GenEnergy{\vecY,3}[\solu](\hstk{T}) {}&{}+c\int_0^T\int_{\rp}^{\rNH}\int_{S^2} \sum_\alpha \left(|\partial_\alpha\solu|^2+\sum_\ua|\partial_\alpha\CQA_\ua\solu|^2\right)\di^2\omega\di\rk\di\tk\nonumber\\
{}&{}\lesssim \GenEnergy{\vecY,3}[\solu](\hstk{0}) +\GenEnergy{\vecTBlend,3}[\solu](\hst{0}) \nonumber\\
{}&{}\lesssim \GenEnergy{\normalhst{0},3}[\solu](\hst{0}) 
\label{eq:NonDegenerateEnergyBound}.
\end{align}
The control over the spacetime integral appearing in the first line of this equation allows us to replace $(\KDelta\dr\solu)^2$ in the Morawetz estimate \ref{Thm:IntroMorawetz} by $(\drk\solu)^2$, thus removing the degeneracy from that estimate in the region $r\in[\rp,\rNH]$ and $\tk\geq0$. (By using coordinates adapted to the past horizon, a similar control can be obtained near the past horizon and away from the bifurcation sphere. The Cauchy stability argument handles the region near the bifurcation sphere.) 

The crucial positivity and negativity properties arising from
conditions (i)-(v) can be found in the work of Dafermos-Rodnianski
\cite{DafermosRodnianski:RedShiftSchwarzschild,DafermosRodnianski:LectureNotes}. For
the sake of completeness, we calculate the the energy on
$\hstk{t_i}\cap\{r\in\{\rp,\rNH\}\}$, 
\begin{align*}
\GenEnergy{\vecY}{}&{}(\hstk{t_i}\cap\{r\in\{\rp,\rNH\}\})\\
{}&{}=-\int_{[\rp,\rNH]\times S^2} \GenMomentum{\vecY}^\alpha\diNormal_\alpha \\
{}&{}=-\int_{[\rp,\rNH]\times S^2} (\nabla_\alpha\solu)(\nabla_\beta\solu)\vecY^\beta \diNormal^\alpha \\
{}&{}\quad+\frac12\int_{[\rp,\rNH]\times S^2} \gMetric^{\gamma\delta}(\nabla_\gamma\solu)(\nabla_\delta\solu) \vecY_\alpha\diNormal^\alpha \\
%{}&{}=\int_{[\rp,\rNH]\times S^2} (\KSigma\dtk\solu-2Mr\drk\solu)(\h\dtk\solu+\f\drk\solu) \sin\theta \di\rk\di\hk\di\pk \\
%{}&{}\quad+\frac12\int_{[\rp,\rNH]\times S^2}(\h-\f)\left(\KDelta(\drk\solu)^2+2(r^2+a^2)(\drk\solu)(\dtk\solu)+2a(\drk\solu)(\dpk\solu)+2a(\dtk\solu)(\dpk\solu) +\opQT^{\alpha\beta}(\partial_\alpha\solu)(\partial_\beta\solu) \right)\sin\theta\di\rk\di\hk\di\pk ,\\
{}&{}=\int_{[\rp,\rNH]\times S^2} \Big(\left(-2Mr\f-\frac12(\h-\f)\KDelta\right)(\drk\solu)^2  \\
{}&{}\qquad\qquad\qquad\qquad +\left(\KDelta\h-a^2\cos^2\theta\f \right)(\drk\solu)(\dtk\solu) \\
{}&{}\qquad\qquad\qquad\qquad +\left(\KSigma\h \right)(\dtk\solu)^2\\
{}&{}\qquad\qquad\qquad\qquad +\left(\h-\f\right)\opQT^{\alpha\beta}(\partial_\alpha\solu)(\partial_\beta\solu) \\
{}&{}\qquad\qquad\qquad\qquad +a\f(\dpk\solu)(\drk\solu)\\
{}&{}\qquad\qquad\qquad\qquad +a\h\frac{\KSigma}{\rp^2+a^2}(\dpk\solu)(\dtk\solu)
\Big)\sin\theta\di\rk\di\hk\di\pk ,
\end{align*} 
the energy on $\hsrk{[t_1,t_2]}$
\begin{align*}
\GenEnergy{\vecY}(\hsrk{[t_1,t_2]})
%{}&{}=\left(\KDelta\drk\solu+(r^2+a^2)\vecTk\solu+a\left(1-\frac{r^2+a^2}{\rp^2+a^2}\right)\dpk\solu\right)(\h\vecTk\solu+\f\drk\solu)-\frac12\f\Lagrangiank \\
%{}&{}=\frac12\f\KDelta(\drk\solu)^2\\
%{}&{}\quad-(\drk\solu)(\vecTk\solu)(-\KDelta\h) \\
%{}&{}\quad+(\vecTk\solu)^2(\h(r^2+a^2))\\
%{}&{}\quad-(\vecTk\solu)(\dpk\solu)\left(-a\h\left(1-\frac{r^2+a^2}{\rp^2+a^2}\right)\right)\\
%{}&{}\quad-\frac12\f\opQT^{\alpha\beta}(\partial_\alpha\solu)(\partial_\beta\solu)\\
=\int_{[t_1,t_2]\times S^2}
{}&{}\left(\h(r^2+a^2)(\vecTk\solu)^2 -\frac{\f}{2}\opQT^{\gamma\delta}(\partial_\alpha\solu)(\partial_\beta\solu)
\right)\sin\theta\di\tk\di\hk\di\pk  ,
\end{align*}
and the divergence of the momentum, using lemma \ref{Lemma:DivOfP}, 
\begin{align}
-\KSigma\nabla_\alpha\GenMomentum{\vecY}^\alpha
%{}&{}=(\h'\vecTk\solu+\f'\drk\solu)\left(-\KDelta\drk\solu-(r^2+a^2)\vecTk\solu-a\left(1-\frac{r^2+a^2}{\rp^2+a^2}\right)(\dpk\solu)\right)\\
%{}&{}\quad+\frac12f'\Lagrangiank+\frac12\f(\drk\CurlyGk^{\alpha\beta})(\partial_\alpha\solu)(\partial_\beta\solu)\\
%{}&{}=(\drk\solu)^2(-\f'\KDelta+\frac12\f'\KDelta+\frac12\f(2r-2M))\\
%{}&{}\quad+(\drk\solu)(\vecTk\solu)(-\h\KDelta-f'(r^2+a^2)+\frac12\f'2(r^2+a^2)+\frac12\f4r)\\
%{}&{}\quad+(\vecTk\solu)^2(-\h'(r^2+a^2))\\
%{}&{}\quad+(\drk\solu)(\dpk\solu)\left(\frac12\f2a\left(\drk\left(1-\frac{r^2+a^2}{\rp^2+a^2}\right)\right)\right)\\
%{}&{}\quad+(\vecTk\solu)(\dpk\solu)\left(-\h'a\left(1-\frac{r^2+a^2}{\rp^2+a^2}\right)\right)\\
%{}&{}\quad-\frac12\f'\opQT^{\alpha\beta}(\partial_\alpha\solu)(\partial_\beta\solu)
%\\
{}&{}=(\drk\solu)^2(\f(r-M)-\frac12\f'\KDelta)
+(\vecTk\solu)^2(-\h'(r^2+a^2))
+\frac12\f'\opQT^{\alpha\beta}(\partial_\alpha\solu)(\partial_\beta\solu)
\nonumber\\
{}&{}\quad+(\drk\solu)(\vecTk\solu)(-\h\KDelta+2r\f)
+(\drk\solu)(\dpk\solu)\left(\f a\left(\drk\left(1-\frac{r^2+a^2}{\rp^2+a^2}\right)\right)\right)\nonumber\\
{}&{}\quad+(\vecTk\solu)(\dpk\solu)\left(-\h'a\left(1-\frac{r^2+a^2}{\rp^2+a^2}\right)\right) .
\label{eq:NonDegenerateMorawetzBulk}
\end{align}

In considering the positivity or negativity of these terms, it is
convenient to, at first, ignore all factors involving $\KDelta$ or
$a$. One can see that, from conditions (i)-(iv), the coefficients of
$(\drk\solu)^2$, $(\vecTk\solu)^2$, and
$\opQT^{\alpha\beta}(\partial_\alpha\solu)(\partial_\beta\solu)$ all
have the desired sign, except for the $(\drk\solu)^2$ term in
$\GenEnergy{\vecY}(\hsrk{[t_1,t_2]})$ which vanishes. With two exceptions, all
the other terms have either a factor of $a$ or $\KDelta$, so they are
small and can immediately be estimated using the Cauchy-Schwarz inequality, possibly at the expense of introducing a new, smaller upper bound for the rotation parameter, $|a|\leq\epsilonInverseAppendix$. The smallness of the factors involving $\KDelta$ near $r=\rp$ imposes the first smallness condition on $\epsilonNH$. 

Of the two exceptional terms, the first is the term involving
$(\drk\solu)(\vecTk\solu)$ in $\GenEnergy{\vecY}(\hsrk{[t_1,t_2]})$. The
potential problem here is that the coefficient of $(\drk\solu)^2$
vanishes linearly in $r-\rp$, so one must take care in applying the
Cauchy-Schwarz estimate. This term can be estimated by
$|(\drk\solu)(\vecTk\solu)\h\KDelta|\leq
\h\KDelta^{1/2}(\KDelta(\drk\solu)^2+(\vecTk\solu)^2)$ and choosing
$\epsilonNH$ sufficiently small that $\epsilonNH\h<|\f|/10$. The
second exceptional term is the term involving
$(\drk\solu)(\vecTk\solu)$ in the divergence. Using condition (v),
this term can be estimated by the Cauchy-Schwarz estimate. 

Since in the support of $\chiNH$, we have $\vecTk=\vecTBlend$ and adding a positive constant to $\h$ preserves conditions (i)-(v), we have a uniform bound 
\begin{align*}
\GenEnergy{\vecY+\vecTBlend,3}(\hstk{T})
\lesssim \GenEnergy{\normalhst{0},3}(\hst{0}) .
\end{align*}
This provides a nondegenerate energy. The nondegenerate Morawetz estimate follows from combining the degenerate Morawetz estimate \eqref{Thm:IntroMorawetz} and estimate \eqref{eq:NonDegenerateEnergyBound}. 

Since $\GenEnergy{\vecY+\vecTBlend,3}$ dominates the integral of
$r^2|\dr\CQA_\ua\solu|^2 +|\CQA_\ua\solu|^2$ for $r>\rNH$ and of $|\drk\CQA_\ua\solu|^2+|\CQA_\ua\solu|^2$ for $r\in[\rp,\rNH]$, this energy also dominates $\sup_{r>\rp}\int_{S^2} \normPtwiseTn{2}{\solu}^2 \di^2\omega$. From the spherical Sobolev estimate \ref{Lemma:STwoSobolev}, we can conclude that there is a uniform constant $C$ such that $\forall t\in\Reals$, $r>\rp$, $(\theta,\phi)\in S^2$
\begin{align*}
|\solu(t,r,\theta,\phi)|
\leq C \GenEnergy{\normalhst{0},3}[\solu](\hst{0})^{1/2} .
\end{align*}

%%%%%%%%%%%%%%%%%%%%%%%%%%%%%%%%%%%%%%%%%%%%%%%%%%%%%%%%%%%%%%%%%%%%%%%%%%%%%%

\section{The Carter operator and the hidden symmetry in Boyer-Lindquist coordinates}
\label{S:CarterOperatorAndQ}
The purpose of this appendix is to compare the operator $\nabla_\alpha K^{\alpha\beta}\nabla_\beta$ arising from Killing tensor associated to Carter's constant, and the operator $\OpQ$ which it turns out to be  convenient to work with in Boyer-Lindquist coordinates. 
The Killing tensor found by Walker and Penrose
\cite{PenroseWalker:1970} to be associated to Carter's constant, is
given by
\begin{align*}
\TensorKCarter^{\alpha\beta}{}&{}=2 \KSigma l^{(\alpha}n^{\beta)} + r^2\gMetric^{\alpha\beta} ,
\end{align*}
see also \cite[\S 12.3]{wald:MR757180}, 
where $l^\alpha$ and $n^\alpha$ are null vectors with $l^\alpha n_\alpha=-1$ 
and orthogonal to 
$$
\Theta=\dtheta, \quad 
\Phi=\frac{1}{\sin\theta}\left(\dphi+a\sin^2\theta\dt\right). 
$$
Carter's constant $\KCarter$ is given by 
$$
\KCarter = \TensorKCarter_{\alpha\beta} \GeodesicVelocity^\alpha \GeodesicVelocity^\beta. 
$$
The Killing tensor can be written in terms of the vectors $\Theta, \Phi$ as 
\begin{align*}
\TensorKCarter^{\alpha\beta}{}&{}=(-\KSigma+r^2)\gMetric^{\alpha\beta}+\Theta^\alpha\Theta^\beta+\Phi^\alpha\Phi^\beta .
\end{align*}
The operator $\nabla_\alpha \TensorKCarter^{\alpha\beta}\nabla_\beta$, which commutes with the d'Alembertian   $\nabla^\alpha\nabla_\alpha$ can be simplified by using standard formulas for divergences in terms of the volume form $\KSigma\Svol$ and by noting that $-\KSigma+r^2=-a^2\cos^2\theta$ depends only on $\theta$. One finds 
\begin{align*}
\nabla_\alpha \TensorKCarter^{\alpha\beta}\nabla_\beta
%{}&{}=(-a^2\cos^2\theta)\nabla_\alpha g^{\alpha\beta}\nabla_\beta
%+(\nabla_\alpha (-a^2\cos^2\theta))\gMetric^{\alpha\beta}\nabla_\beta
%+\nabla_\alpha(\Theta^\alpha\Theta^\beta+\Phi^\alpha\Phi^\beta)\\
%{}&{}=-a^2\cos^2\theta\nabla^\alpha\nabla_\alpha
%-(\partial_\alpha \KSigma)\KSigma^{-1}\Theta^\alpha\Theta^\beta\partial_\beta
%+\frac{1}{\KSigma\Svol}\partial_\theta %(\KSigma\Svol(\Theta^\alpha\Theta^\beta+\Phi^\alpha\Phi^\beta))\partial_\beta\\
%{}&{}=-a^2\cos^2\theta\nabla^\alpha\nabla_\alpha 
%+\frac{1}%{\Svol}\partial_\alpha(\Svol(\Theta^\alpha\Theta^\beta+\Phi^\alpha\Phi^\beta))\partial_\beta\\
{}&{}=-a^2\cos^2\theta\nabla^\alpha\nabla_\alpha 
+\OpQ +\dphi^2 +2a\dt\dphi .
\end{align*}
where $\OpQ$ is given by
\begin{equation}%\label{eq:OpQdefappendix}
\OpQ = \frac1{\sin\theta}\dtheta\sin\theta\dtheta
+\frac{\cos^2\theta}{\sin^2\theta}\dphi^2 +a^2\sin^2\theta \dt^2 , 
\end{equation}
as in \eqref{eq:OpQdef}.

\subsection*{Acknowledgements} 
A significant part of this work was completed during the ``Geometry, Analysis,
and General Relativity'' programme held at the Mittag-Leffler Institute, 
Djursholm, Sweden, in the fall of 2008. Both authors are grateful to the
institute for hospitality and support during this time, and to many of the
participants in the programme for useful discussions. Both authors are grateful to the anonymous referees for several helpful suggestions. PB thanks the Albert
Einstein Institute for support during the initial phase of this project. LA was
supported in part by the NSF, under contract no. DMS-0707306.

\newcommand{\prd}{Phys. Rev. D} 
\bibliography{HiddenSymmetriesDecayWaveKerr}
\bibliographystyle{abbrv}

\end{document}